\documentclass{birkjour}

\newif\ifignore % when set to true, additional text appears containing
\ignorefalse

\newcommand{\auxproof}[1]{
\ifignore\mbox{}\newline
\textbf{PROOF:} \dotfill\newline
{\it #1}\mbox{}\newline
\textbf{ENDPROOF}\dotfill
\fi}

\usepackage{amsmath}
\usepackage{amssymb}
\usepackage{stmaryrd}
\usepackage{xspace}
\usepackage{xy}
\xyoption{all}
\usepackage{hyperref}

\newdimen\proofrulebreadth \proofrulebreadth=.05em
\newdimen\proofdotseparation \proofdotseparation=1.25ex
\newdimen\proofrulebaseline \proofrulebaseline=2ex
\newcount\proofdotnumber \proofdotnumber=3
\let\then\relax
\def\hfi{\hskip0pt plus.0001fil}
\mathchardef\squigto="3A3B
\newif\ifinsideprooftree\insideprooftreefalse
\newif\ifonleftofproofrule\onleftofproofrulefalse
\newif\ifproofdots\proofdotsfalse
\newif\ifdoubleproof\doubleprooffalse
\let\wereinproofbit\relax
\newdimen\shortenproofleft
\newdimen\shortenproofright
\newdimen\proofbelowshift
\newbox\proofabove
\newbox\proofbelow
\newbox\proofrulename
\def\shiftproofbelow{\let\next\relax\afterassignment\setshiftproofbelow\dimen0 }
\def\shiftproofbelowneg{\def\next{\multiply\dimen0 by-1 }%
\afterassignment\setshiftproofbelow\dimen0 }
\def\setshiftproofbelow{\next\proofbelowshift=\dimen0 }
\def\setproofrulebreadth{\proofrulebreadth}

\def\prooftree{% NESTED ZERO (\ifonleftofproofrule)
\ifnum  \lastpenalty=1
\then   \unpenalty
\else   \onleftofproofrulefalse
\fi
\ifonleftofproofrule
\else   \ifinsideprooftree
        \then   \hskip.5em plus1fil
        \fi
\fi
\bgroup% NESTED ONE (\proofbelow, \proofrulename, \proofabove,
\setbox\proofbelow=\hbox{}\setbox\proofrulename=\hbox{}%
\let\justifies\proofover\let\leadsto\proofoverdots\let\Justifies\proofoverdbl
\let\using\proofusing\let\[\prooftree
\ifinsideprooftree\let\]\endprooftree\fi
\proofdotsfalse\doubleprooffalse
\let\thickness\setproofrulebreadth
\let\shiftright\shiftproofbelow \let\shift\shiftproofbelow
\let\shiftleft\shiftproofbelowneg
\let\ifwasinsideprooftree\ifinsideprooftree
\insideprooftreetrue
\setbox\proofabove=\hbox\bgroup$\displaystyle % NESTED TWO
\let\wereinproofbit\prooftree
\shortenproofleft=0pt \shortenproofright=0pt \proofbelowshift=0pt
\onleftofproofruletrue\penalty1
}

\def\eproofbit{% NESTED TWO
\ifx    \wereinproofbit\prooftree
\then   \ifcase \lastpenalty
        \then   \shortenproofright=0pt  % 0: some other object, no indentation
        \or     \unpenalty\hfil         % 1: empty hypotheses, just glue
        \or     \unpenalty\unskip       % 2: just had a tree, remove glue
        \else   \shortenproofright=0pt  % eh?
        \fi
\fi
\global\dimen0=\shortenproofleft
\global\dimen1=\shortenproofright
\global\dimen2=\proofrulebreadth
\global\dimen3=\proofbelowshift
\global\dimen4=\proofdotseparation
\global\count255=\proofdotnumber
$\egroup  % NESTED ONE
\shortenproofleft=\dimen0
\shortenproofright=\dimen1
\proofrulebreadth=\dimen2
\proofbelowshift=\dimen3
\proofdotseparation=\dimen4
\proofdotnumber=\count255
}

\def\proofover{% NESTED TWO
\eproofbit % NESTED ONE
\setbox\proofbelow=\hbox\bgroup % NESTED TWO
\let\wereinproofbit\proofover
$\displaystyle
}%
\def\proofoverdbl{% NESTED TWO
\eproofbit % NESTED ONE
\doubleprooftrue
\setbox\proofbelow=\hbox\bgroup % NESTED TWO
\let\wereinproofbit\proofoverdbl
$\displaystyle
}%
\def\proofoverdots{% NESTED TWO
\eproofbit % NESTED ONE
\proofdotstrue
\setbox\proofbelow=\hbox\bgroup % NESTED TWO
\let\wereinproofbit\proofoverdots
$\displaystyle
}%
\def\proofusing{% NESTED TWO
\eproofbit % NESTED ONE
\setbox\proofrulename=\hbox\bgroup % NESTED TWO
\let\wereinproofbit\proofusing
\kern0.3em$
}

\def\endprooftree{% NESTED TWO
\eproofbit % NESTED ONE
  \dimen5 =0pt% spread of hypotheses
\dimen0=\wd\proofabove \advance\dimen0-\shortenproofleft
\advance\dimen0-\shortenproofright
\dimen1=.5\dimen0 \advance\dimen1-.5\wd\proofbelow
\dimen4=\dimen1
\advance\dimen1\proofbelowshift \advance\dimen4-\proofbelowshift
\ifdim  \dimen1<0pt
\then   \advance\shortenproofleft\dimen1
        \advance\dimen0-\dimen1
        \dimen1=0pt
        \ifdim  \shortenproofleft<0pt
        \then   \setbox\proofabove=\hbox{%
                        \kern-\shortenproofleft\unhbox\proofabove}%
                \shortenproofleft=0pt
        \fi
\fi
\ifdim  \dimen4<0pt
\then   \advance\shortenproofright\dimen4
        \advance\dimen0-\dimen4
        \dimen4=0pt
\fi
\ifdim  \shortenproofright<\wd\proofrulename
\then   \shortenproofright=\wd\proofrulename
\fi
\dimen2=\shortenproofleft \advance\dimen2 by\dimen1
\dimen3=\shortenproofright\advance\dimen3 by\dimen4
\ifproofdots
\then
        \dimen6=\shortenproofleft \advance\dimen6 .5\dimen0
        \setbox1=\vbox to\proofdotseparation{\vss\hbox{$\cdot$}\vss}%
        \setbox0=\hbox{%
                \advance\dimen6-.5\wd1
                \kern\dimen6
                $\vcenter to\proofdotnumber\proofdotseparation
                        {\leaders\box1\vfill}$%
                \unhbox\proofrulename}%
\else   \dimen6=\fontdimen22\the\textfont2 % height of maths axis
        \dimen7=\dimen6
        \advance\dimen6by.5\proofrulebreadth
        \advance\dimen7by-.5\proofrulebreadth
        \setbox0=\hbox{%
                \kern\shortenproofleft
                \ifdoubleproof
                \then   \hbox to\dimen0{%
                        $\mathsurround0pt\mathord=\mkern-6mu%
                        \cleaders\hbox{$\mkern-2mu=\mkern-2mu$}\hfill
                        \mkern-6mu\mathord=$}%
                \else   \vrule height\dimen6 depth-\dimen7 width\dimen0
                \fi
                \unhbox\proofrulename}%
        \ht0=\dimen6 \dp0=-\dimen7
\fi
\let\doll\relax
\ifwasinsideprooftree
\then   \let\VBOX\vbox
\else   \ifmmode\else$\let\doll=$\fi
        \let\VBOX\vcenter
\fi
\VBOX   {\baselineskip\proofrulebaseline \lineskip.2ex
        \expandafter\lineskiplimit\ifproofdots0ex\else-0.6ex\fi
        \hbox   spread\dimen5   {\hfi\unhbox\proofabove\hfi}%
        \hbox{\box0}%
        \hbox   {\kern\dimen2 \box\proofbelow}}\doll%
\global\dimen2=\dimen2
\global\dimen3=\dimen3
\egroup % NESTED ZERO
\ifonleftofproofrule
\then   \shortenproofleft=\dimen2
\fi
\shortenproofright=\dimen3
\onleftofproofrulefalse
\ifinsideprooftree
\then   \hskip.5em plus 1fil \penalty2
\fi
}

\newtheorem{thm}{Theorem}[section]

\newtheorem{lemma}[thm]{Lemma}
\newtheorem{proposition}[thm]{Proposition}
\theoremstyle{definition}
\newtheorem{definition}[thm]{Definition}
\theoremstyle{remark}
\newtheorem{remark}[thm]{Remark}
\newtheorem{example}[thm]{Example}

\renewcommand{\arraycolsep}{2pt}

\makeatother
 \newdir{ >}{{}*!/-7.5pt/@{>}}
 \newdir{|>}{!/4.5pt/@{|}*:(1,-.2)@^{>}*:(1,+.2)@_{>}}
 \newdir{ |>}{{}*!/-3pt/@{|}*!/-7.5pt/:(1,-.2)@^{>}*!/-7.5pt/:(1,+.2)@_{>}}
\newcommand{\xyline}[2][]{\ensuremath{\smash{\xymatrix@1#1{#2}}}}
\newcommand{\xyinline}[2][]{\ensuremath{\smash{\xymatrix@1#1{#2}}}^{\rule[8.5pt]{0pt}{0pt}}}
\makeatletter

\newenvironment{myproof}[1][Proof]%
   { \begin{trivlist}%
     \item[\hskip \labelsep {\bfseries #1}]%
   }%
   { \end{trivlist}%
   }
\newcommand{\QEDbox}{\square}
\newcommand{\QED}{\hspace*{\fill}$\QEDbox$}

\newcommand{\after}{\mathrel{\circ}}
\newcommand{\relafter}{\mathrel{\bullet}}

\newcommand{\cat}[1]{\ensuremath{\mathbf{#1}}}
\newcommand{\Cat}[1]{\ensuremath{\mathbf{#1}}}

\newcommand{\op}{\ensuremath{^{\mathrm{op}}}}
\newcommand{\idmap}[1][]{\ensuremath{\mathrm{id}_{#1}}}

\newcommand{\support}{\ensuremath{\mathrm{supp}}}

\renewcommand{\ker}{\ensuremath{\mathsf{ker}}}
\newcommand{\Ker}{\ensuremath{\mathsf{Ker}}}
\newcommand{\eq}[1][]{\ensuremath{\mathsf{eq}_{#1}}}

\newcommand{\cp}[1][]{\ensuremath{\mathsf{cp}_{#1}}}
\newcommand{\cpMlt}{\cp[\Mlt]}
\newcommand{\cpPow}{\cp[\Pow]}
\newcommand{\cpLft}{\cp[\Lft]}

\newcommand{\SA}{\ensuremath{\mathcal{S}{\kern-.95ex}\mathcal{A}}}

\newcommand{\Pos}{\ensuremath{\mathcal{P}{\kern-.6ex}o{\kern-.35ex}s}}

\newcommand{\Proj}{\ensuremath{\mathcal{P}{\kern-.45ex}r}}
\newcommand{\DM}{\ensuremath{\mathcal{D}{\kern-.85ex}\mathcal{M}}}

\newcommand{\Hom}{\textsl{Hom}}

\newcommand{\inprod}[2]{\ensuremath{\langle #1\,|\,#2 \rangle}}

\newcommand{\Vect}{\cat{Vect}\xspace}
\newcommand{\Alg}{\textsl{Alg}\xspace}
\newcommand{\Kl}{\mathcal{K}{\kern-.5ex}\ell}
\newcommand{\TRel}{\textsl{TRel}\xspace}

\newcommand{\Mlt}{\ensuremath{\mathcal{M}}}
\newcommand{\Dst}{\ensuremath{\mathcal{D}}}
\newcommand{\Lft}{\ensuremath{\mathcal{L}}}
\newcommand{\UF}{\ensuremath{\mathcal{U}{\kern-.75ex}\mathcal{F}}}
\newcommand{\CS}{\ensuremath{\mathcal{C}{\kern-.75ex}\mathcal{S}}}

\newcommand{\NNO}{\ensuremath{\mathbb{N}}}
\newcommand{\unitR}{\ensuremath{[0,1]}}

\newcommand{\closed}{\ensuremath{\mathcal{C}{\kern-.45ex}\ell}}

\newcommand{\orthogonal}{\mathrel{\bot}}

\newcommand{\powerset}{\mathcal{P}}
\newcommand{\Pow}{\powerset}

\newcommand{\Powfin}{\Pow_{\mathit{fin}}}

\newcommand{\tr}{\ensuremath{\textrm{tr}}\xspace}

\newcommand{\leftScottint}{[{\kern-.3ex}[}
\newcommand{\rightScottint}{]{\kern-.3ex}]}

\newcommand{\BifRel}{\Cat{BifRel}\xspace}
\newcommand{\dBisRel}{\Cat{dBisRel}\xspace}
\newcommand{\BifMRel}{\Cat{BifMRel}\xspace}
\newcommand{\PInj}{\Cat{PInj}\xspace}
\newcommand{\Pfn}{\Cat{Pfn}\xspace}
\newcommand{\Hilb}{\Cat{Hilb}\xspace}

\newcommand{\Sets}{\Cat{Sets}\xspace}
\newcommand{\PoSets}{\Cat{PoSets}\xspace}

\newcommand{\JSL}{\Cat{JSL}\xspace}
\newcommand{\Mod}{\Cat{Mod}\xspace}

\newcommand{\sotimes}{\mathrel{\raisebox{.05pc}{$\scriptstyle \otimes$}}}

\newcommand{\tuple}[1]{\ensuremath{\langle #1 \rangle}}

\newcommand{\set}[2]{\{#1\;|\;#2\}}
\newcommand{\setin}[3]{\{#1\in#2\;|\;#3\}}
\newcommand{\conjun}{\mathrel{\wedge}}

\newcommand{\all}[2]{\forall{#1}.\,#2}
\newcommand{\allin}[3]{\forall{#1\in#2}.\,#3}
\newcommand{\ex}[2]{\exists{#1}.\,#2}
\newcommand{\exin}[3]{\exists{#1\in#2}.\,#3}
\newcommand{\lamin}[3]{\lambda{#1\in#2}.\,#3}

\newcommand{\EndoHom}[1]{{\cal E}{\kern-.5ex}\textit{n}{\kern-.2ex}\textit{d}{\kern-.2ex}\textit{o}(#1)}
\newcommand{\congrightarrow}{\mathrel{\stackrel{
           \raisebox{.5ex}{$\scriptstyle\cong\,$}}{
           \raisebox{0ex}[0ex][0ex]{$\rightarrow$}}}}
\newcommand{\conglongrightarrow}{\mathrel{\stackrel{
           \raisebox{.5ex}{$\scriptstyle\cong\,$}}{
           \raisebox{0ex}[0ex][0ex]{$\longrightarrow$}}}}

\begin{document}

\title{Dagger Categories of Tame Relations}
\author{Bart Jacobs} 
\address{
Institute for Computing and Information Sciences (iCIS), \\
Radboud University Nijmegen, The Netherlands. \\
Webaddress: \url{www.cs.ru.nl/B.Jacobs} \\[.5em]
\today
}

\date{\small \today}

\maketitle

\begin{abstract}
Within the context of an involutive monoidal category the notion of a
comparison relation $\cp\colon \overline{X} \otimes X\rightarrow
\Omega$ is identified. Instances are equality $=$ on sets, inequality
$\leq$ on posets, orthogonality $\orthogonal$ on orthomodular
lattices, non-empty intersection on powersets, and inner product
$\inprod{-}{-}$ on vector or Hilbert spaces. Associated with a
collection of such (symmetric) comparison relations a dagger category
is defined with ``tame'' relations as morphisms. Examples include
familiar categories in the foundations of quantum mechanics, such as
sets with partial injections, or with locally bifinite relations, or
with formal distributions between them, or Hilbert spaces with bounded
(continuous) linear maps. Of one particular example of such a dagger
category of tame relations, involving sets and bifinite multirelations
between them, the categorical structure is investigated in some
detail. It turns out to involve symmetric monoidal dagger structure,
with biproducts, and dagger kernels. This category may form an
appropriate universe for discrete quantum computations, just like
Hilbert spaces form a universe for continuous computation.
\end{abstract}

% Orthomodular lattices: 1998 ACM Subject Classification: F.4.1

%\subjclass{Primary 68Q55; Secondary 18D10, 81P68}
%\keywords{Dagger category, quantum semantics}

% 18D10 Monoidal categories (= multiplicative categories), symmetric monoidal categories, braided categories

% 81P68 Quantum computation

% 68Q55 Semantics

\section{Introduction}\label{IntroSec}

So-called tame relations were introduced in~\cite{BluteP11} in the
construction of a particular (monoidal) dagger category of formal
distributions. The phrase `tame' refers to finiteness restrictions in
two directions, and is best illustrated in the context of
relations. So suppose we have a relation $r\subseteq X\times Y$; it
can be described equivalently as a function $X\rightarrow \Pow(Y)$,
where $\Pow$ is powerset, or via reversal, as a function
$Y\rightarrow\Pow(X)$. The relation is called tame, if both these
functions factorise via the \emph{finite} powerset $\Powfin$, as in
$X\rightarrow\Powfin(Y)$ and $Y\rightarrow\Powfin(X)$. Concretely,
this means that for each $x\in X$ there are only finitely many $y\in
Y$ with $r(x,y)$, and vice-versa. Such relations are often called
(locally) bifinite. They may be used to model finitely
non-determinstic reversible computations.

In~\cite{BluteP11} tameness is used in the context of polynomials. Let
$S[X]$ be the collection of (multivariate) polynomials, with variables
in a set $X$ and coefficients in a semiring $S$; similarly $S[[X]]$ is
used for possibly infinite such polynomials (or power series, or
formal distributions). For certain analogues of relations, giving rise
to mappings $S[X]\rightarrow S[[Y]]$ and $S[Y]\rightarrow S[[X]]$,
tameness means that these mappings factor via finite polynomials, as
in $S[X]\rightarrow S[Y]$ and $S[Y]\rightarrow S[X]$. For details, see
Subsection~\ref{FormDistrSubsec} below.

This paper starts by identifying a general context in which this
notion of tameness makes sense. It involves the notion of a comparison
relation $\cp\colon \overline{X} \otimes X\rightarrow \Omega$.  Such a
relation requires an ambient category with tensors $\otimes$ and
involution $\overline{(-)}$, as described for instance
in~\cite{Egger11,BeggsM09,Jacobs11e}. A relation $r\colon
\overline{X}\otimes Y \rightarrow\Omega$ is then tame, if it factors
via such comparisons, via appropriate maps $r_{*}$ and $r^{*}$. It is
shown that such categories of tame relations give rise to dagger
categories, assuming the underlying comparison relation is
symmetric. Section~\ref{ExSec} illustrates how this general
construction encompasses several known categories that are relevant in
the foundations of quantum mechanics, such as orthomodular lattices
with Galois connections, or sets with partial injections or with
bifinite relations, or with bifinite multirelations, or with
bistochastic relations. Some of these constructions are also described
more abstractly, in terms of the monads involved, namely lift, finite
powerset, multiset and distribution monads, see
Subsection~\ref{MonadSubsec}.

The formal distributions example from~\cite{BluteP11} is re-described
in the present general setting. Additionally, bounded (or continuous)
maps between Hilbert spaces are shown to correspond to tame relations
(see Lemma~\ref{BoundedLem}).

Finally, one particular example category of tame relations, arising
via the multiset monad from the monad construction just mentioned, is
further investigated in Section~\ref{BfMRelSec}. We refer to this as
the category $\BifMRel$ of sets and bifinite multirelations. Morphisms
$X\rightarrow Y$ are functions $r\colon X\times Y\rightarrow S$, into
a semiring $S$, such that for each $x\in X$ there are only finitely
many $y$ with $r(x,y)\neq 0$, and vice-versa. This means that the
relation factors both as $X\rightarrow\Mlt_{S}(Y)$ and as
$Y\rightarrow\Mlt_{S}(X)$, where $\Mlt_{S}$ is the multiset monad
which ``counts in $S$''. Such bifinite multirelations may be used to
model finitely weighted, reversible computations. It is show that this
category $\BifMRel$ has, besides daggers, tensors $\otimes$ and
biproducts $\oplus$. Moreover, it has dagger kernels, as described
in~\cite{HeunenJ10a}. Thus, the category $\BifMRel$ resembles the
category $\Hilb$ of Hilbert spaces. It is suggested that this category
$\BifMRel$ is the discrete analogue of \Hilb, useful for discrete
quantum computations, such as usually occurring in a quantum computer
science context (see \textit{e.g.}~\cite{Mermin07,NielsenC00}). The
quantum walks example from~\cite{Jacobs11a}, formalised in $\BifMRel$,
supports this suggestion, but further evidence is required via more
extensive investigation.

Thus, the contributions of the paper are two-fold: (1)~identifying the
uniformity in various models of quantum computation via a systematic
exposition in terms of comparison relations, and (2)~first
investigation of one particular promising example of such a model for
discrete quantum computation, namely the category of sets and bifinite
multirelations.

\section{Involutive categories, and comparisons therein}\label{InvCatSec}

This section recalls the basics of involutive categories as presented
in~\cite{Jacobs11e} (see also~\cite{BeggsM09,Egger11}). Within such
involutive categories the notion of `comparison' is introduced.

A category $\cat{A}$ will be called \textbf{involutive} if it comes
with a `involution' functor $\cat{A} \rightarrow \cat{A}$, written as
$X\mapsto \overline{X}$, and a natural isomorphism $\iota_{X} \colon
X\conglongrightarrow \overline{\overline{X}}$ satisfying
$\overline{\iota_{X}} = \iota_{\overline{X}} \colon \overline{X}
\rightarrow \overline{\overline{\overline{X}}}$.

%% In this situation it is not hard to show that there is a
%% self-adjointness $\overline{(-)} \dashv \overline{(-)}$, so that
%% the involution functor $\overline{(-)}\colon\cat{A}\rightarrow\cat{A}$
%% preserves all limits and colimits that exist in $\cat{A}$.

Within such an involutive category a \textbf{self-conjugate} is an
object $X$ with a map $j\colon\overline{X}\rightarrow X$ satisfying
$j\after \overline{j} = \iota^{-1} \colon
\overline{\overline{X}}\rightarrow X$.  Such a map $j$ is necessarily
an isomorphism. A self-conjugate is called a star-object
in~\cite{BeggsM09}.

Each category is trivially involutive via the identity functor. The
category \PoSets is involutive via order reversal $(-)\op$. This
applies also to categories of, for instance, distributive lattices or
Boolean algebras. Probably the most standard example of an involutive
category is the category $\Vect_{\mathbb{C}}$ of vector spaces over
the complex numbers $\mathbb{C}$; it is involutive via conjugation:
for a vector space $V\in\Vect_{\mathbb{C}}$ there is the `complex
conjugate' space $\overline{V}\in\Vect_{\mathbb{C}}$ with the same
vectors as $V$, but with adapted scalar multiplication $s
\cdot_{\overline{V}} v = \overline{s} \cdot_{V} v$, for
$s\in\mathbb{C}$ and $v\in V$, where $\overline{s} = a-ib$ is the
conjugate of the complex number $s = a+ib\in\mathbb{C}$. This same
involution exists on categories of Hilbert spaces (over $\mathbb{C}$).

The negation map $\neg \colon B\op \congrightarrow B$ makes each
Boolean algebra $B$ self-conjugate, for the $(-)\op$ involution on the
category of Boolean algebras. The conjugation map $\overline{(-)}$ on
the complex numbers makes $\mathbb{C}$ a self-conjugate
$\overline{\mathbb{C}} \congrightarrow \mathbb{C}$ in the category
of vector (or Hilbert) spaces over $\mathbb{C}$.

\begin{definition}
\label{InvMonCatDef}
An \textbf{involutive (symmetric) monoidal category} is a category
$\cat{A}$ which is both involutive and (symmetric) monodial in which
involution $\overline{(-)}\colon \cat{A} \rightarrow \cat{A}$ is a
(symmetric) monoidal functor---via maps $\zeta\colon I \rightarrow
\overline{I}$ and $\xi\colon \overline{X}\otimes \overline{Y}
\rightarrow \overline{X\otimes Y}$ commuting with the monoidal
isomorphisms---and $\iota \colon \idmap{} \Rightarrow
\overline{\overline{(-)}}$ is a monoidal natural transformation; this
means that the following diagrams commutes.
\begin{equation}
\label{InvMonoidalEqn}
\vcenter{\xymatrix@R1.5pc{
I\ar@{=}[d]\ar@{=}[rr] & & I\ar[d]^{\iota}
& &
X\otimes Y\ar[d]_{\iota\otimes\iota}\ar@{=}[rr] & & 
   X\otimes Y\ar[d]^{\iota} \\
I\ar[r]^-{\zeta} & \overline{I}\ar[r]^-{\overline{\zeta}} &
   \overline{\overline{I}}
& &
\overline{\overline{X}}\otimes\overline{\overline{Y}}\ar[r]^-{\xi} &
   \overline{\overline{X}\otimes\overline{Y}}\ar[r]^-{\overline{\xi}} &
   \overline{\overline{X\otimes Y}}
}}
\end{equation}
\end{definition}

\noindent One can show (see~\cite{Jacobs11e}) that the involution
functor $\overline{(-)}$ is automatically strong monoidal: the maps
$\zeta\colon I \rightarrow \overline{I}$ and $\xi\colon
\overline{X}\otimes \overline{Y} \rightarrow \overline{X\otimes Y}$
are necessarily isomorphisms.

In the symmetric case, with symmetry $\gamma\colon X\otimes Y
\conglongrightarrow Y\otimes X$, we often use the `twist' $\tau$ defined
by:
\begin{equation}
\label{TwistConjugEqn}
\vcenter{\xymatrix@R1.5pc@C+.2pc{
\tau \stackrel{\textrm{def}}{=} \Big(\overline{\overline{X}\otimes Y}
      \ar[r]^-{\xi^{-1}}_-{\cong} &
   \overline{\overline{X}}\otimes\overline{Y}
   \ar[r]^-{\iota^{-1}\otimes\idmap}_-{\cong} &
   X\otimes \overline{Y}\ar[r]^-{\gamma}_-{\cong} & 
   \overline{Y}\otimes X\Big).
}}
\end{equation}

\noindent For $Y=X$ this map makes the object $\overline{X}\otimes X$
self-conjugate.

\auxproof{
This is the same as the (too complicated) map
$$\vcenter{\xymatrix@R1.5pc@C+.2pc{
\tau \stackrel{\textrm{def}}{=} \Big(\overline{\overline{X}\otimes Y}
      \ar[r]^-{\overline{\idmap{}\otimes\iota_{Y}}} &
   \overline{\overline{X}\otimes\overline{\overline{Y}}}
      \ar[r]^-{\overline{\xi}} &
   \overline{\overline{X\otimes \overline{Y}}}\ar[r]^-{\iota^{-1}} &
   X\otimes \overline{Y}\ar[r]^-{\gamma} & \overline{Y}\otimes X\Big)
}}$$

\noindent used in~\cite{Jacobs11e}. Since:
$$\begin{array}{rcl}
\iota^{-1} \after \overline{\xi} \after \overline{\idmap\otimes\iota_{Y}}
& = &
(\iota_{X}^{-1}\otimes\iota_{\overline{Y}}^{-1}) \after \xi^{-1} 
   \after \overline{\idmap\otimes\iota_{Y}} \\
& = &
(\iota_{X}^{-1}\otimes\overline{\iota_{Y}^{-1}}) 
   \after (\overline{\idmap}\otimes\overline{\iota_{Y}}) \after \xi^{-1} \\
& = &
(\iota_{X}^{-1}\otimes\idmap) \after \xi^{-1}.
\end{array}$$

This twist map is indeed a self-conjugate:
$$\xymatrix@R1.5pc@C+1pc{
\overline{\overline{\overline{X}\otimes X}}
      \ar[r]^-{\overline{\xi^{-1}}}\ar@/_8ex/[ddrr]_{\iota^{-1}}
      \ar`u`[rrr]^-{\overline{\tau}}[rrr] &
   \overline{\overline{\overline{X}}\otimes\overline{X}}
      \ar[r]^-{\overline{\iota^{-1}\otimes\idmap}}\ar[d]^{\xi^{-1}} &
   \overline{X\otimes\overline{X}}\ar[r]^-{\overline{\gamma}}
      \ar[d]^{\xi^{-1}} &
   \overline{\overline{X}\otimes X}\ar[d]^{\xi^{-1}}
      \ar`r`[dddl]^-{\tau}[ddd] &
\\
& 
   \overline{\overline{\overline{X}}}\otimes\overline{\overline{X}}
      \ar[r]^-{\overline{\iota^{-1}}\otimes\idmap}
      \ar[dr]_{\iota^{-1}\otimes\iota^{-1}} & 
   \overline{X}\otimes\overline{\overline{X}}\ar[r]^-{\gamma}
      \ar[d]^{\idmap\otimes\iota^{-1}} & 
   \overline{\overline{X}}\otimes \overline{X}\ar[d]^{\iota^{-1}\otimes\idmap}
\\
& & \overline{X}\otimes X\ar[r]^-{\gamma}\ar@{=}[dr] & 
   X\otimes\overline{X}\ar[d]^{\gamma}
\\
& & & \overline{X}\otimes X
}$$
}

\subsection{Comparison relations}\label{ComparRelSubsec}

The equality relation on a set $X$ can be described as a map $= \colon
X\times X\rightarrow 2$ in \Sets, where $2 = \{0,1\}$. We wish to
capture such maps more generally under the name `comparison relation'.

\begin{definition}
\label{ComparRelDef}
Assume an involutive monoidal category with a special object
$\Omega$. A \textbf{comparison relation} is a map of the form
$\cp\colon \overline{X}\otimes X\rightarrow \Omega$ satisfying:
$$f=g \quad\mbox{follows from either}\quad
\left\{\begin{array}{l}
 \cp \after (\idmap\otimes f) = \cp \after (\idmap\otimes g) 
   \quad\mbox{or} \\
 \cp \after (f\otimes\idmap) = \cp \after (g\otimes\idmap)
\end{array}\right.$$

\noindent In presence of exponents $\multimap$, these `mono
requirements' mean that the two associated abstraction maps
$X\rightarrow (\overline{X}\multimap\Omega)$ and
$\overline{X}\rightarrow (X\multimap\Omega)$ are monic.

In a symmetric monoidal setting such a comparison relation is called
\textbf{symmetric} if the following diagram commutes,
\begin{equation}
\label{ComparSymmEqn}
\vcenter{\xymatrix@R1.5pc{
\overline{\overline{X}\otimes X}\ar[d]_{\overline{\cp}}
   \ar[rr]^-{\tau}_-{\cong} & &
   \overline{X}\otimes X\ar[d]^{\cp} \\
\overline{\Omega}\ar[rr]^-{j}_-{\cong} & & \Omega
}}
\end{equation}

\noindent where a self-conjugate structure $\smash{\overline{\Omega}
  \stackrel{j}{\rightarrow} \Omega}$ is assumed, and where $\tau$ is
the twist map from~\eqref{TwistConjugEqn}.
\end{definition}

In the symmetric case the two mono requirements---for each argument
separa\-tely---can be reduced to a single requirement---in one argument
only: if $\cp \after (\idmap\otimes f) = \cp \after (\idmap\otimes g)$
implies $f=g$, then one can deduce that also $\cp \after
(f\otimes\idmap) = \cp \after (g\otimes\idmap)$ implies $f=g$ (and
vice-versa).

\auxproof{
To be really precise:
$$\cp \after (f\otimes\idmap) = \cp \after (g\otimes\idmap)$$

\noindent implies:
$$j \after \overline{\cp \after (\idmap\otimes f)} \after \tau^{-1} = 
  j \after \overline{\cp \after (\idmap\otimes g)} \after \tau^{-1}.$$

\noindent Then one uses:
$$\begin{array}{rcl}
j \after \overline{\cp \after (\idmap\otimes f)} \after \tau^{-1}
& = &
j \after \overline{\cp} \after 
   \overline{(\overline{\idmap}\otimes f)} \after \tau^{-1} \\
& = &
j \after \overline{\cp} \after \tau^{-1} \after 
   (\overline{f}\otimes\idmap) \\
& = &
\cp \after (\overline{f}\otimes\idmap)
\end{array}$$

\noindent so that we can conclude $\overline{f}=\overline{g}$.  But
then also $\iota \after f = \overline{\overline{f}} \after \iota =
\overline{\overline{g}} \after \iota = \iota \after g$. Hence $f=g$
since $\iota$ is an isomorphism.
}

An equality relation $= \colon X\times X\rightarrow 2 = \{0,1\}$ on a
set $X$ is given by $(x = x) = 1$ and $(x = x') = 0$ for $x\neq
x'$. This is a symmetric comparison relation in the category \Sets,
with trivial (identity) involution. More interestingly, for a poset
$(X,\leq)$, the order forms a non-symmetric comparison relation $\leq
\colon X\op\times X\rightarrow 2$ in \PoSets. The involution $(-)\op$
in the type of the first argument is needed for monotonicity, since:
$x\geq x'$ and $x \leq y$ and $y\leq y'$ implies $x'\leq y'$. The mono
requirement translates (in one argument) to: $x=y$ follows from $x\leq
z$ iff $y\leq z$ for all $z$.

A non-trivial symmetric example is the inner product $\inprod{-}{-}
\colon \overline{H}\otimes H\rightarrow \mathbb{C}$ on a Hilbert space
$H$ (over $\mathbb{C}$). The bilinearity and antilinearity
requirements of an inner product are captured via tensor and
conjugation in the input type of the operation: it yields
$\inprod{s\cdot x}{y} = \overline{s}\cdot\inprod{x}{y}$ and
$\inprod{x_{1}+x_{2}}{y} = \inprod{x_1}{y} + \inprod{x_2}{y}$, and
similarly, $\inprod{x}{s\cdot y} = s\cdot\inprod{x}{y}$ and
$\inprod{x}{y_{1}+y_{2}} = \inprod{x}{y_1} + \inprod{x}{y_2}$. The
symmetry requirement for a comparison relation says that
$\overline{\inprod{y}{x}} = \inprod{x}{y}$.  The mono requirement
holds, since if $\inprod{x}{z} = \inprod{y}{z}$ for all $z$, then $0 =
\inprod{x}{z} - \inprod{y}{z} = \inprod{x-y}{z}$.  By taking $z=x-y$
we get $\inprod{x-y}{x-y} = 0$, from which we conclude $x-y=0$ and
thus $x=y$.

\begin{remark}
\label{ZeroRem}
Notice that our notion of comparison does not involve the usual inner
product requirements $\inprod{x}{x} \geq 0$ and $\inprod{x}{x} = 0
\Rightarrow x=0$ for Hilbert spaces. Such requirements are not needed
for what we wish to achieve (in the next section) and involve
additional assumptions, namely the presence of zero objects (or
maps). The kind of inner product that is captured via a comparison
relation corresponds to a Minkowski inner product.

Although we do not pursue this here, we would like to mention that in
presence of such a zero one can introduce complementation with respect
to a comparison relation: for $U\subseteq X$, take $U^{\perp} =
\setin{x}{X}{\allin{x'}{U}{\cp(x,x')=0}}$. For sets this gives
ordinary complement, and for Hilbert spaces it yields
orthocomplementation of closed subsets.

Another point not pursued here is the similarity with
profunctors~\cite{Benabou73,Lawvere73}, commonly understood as
`categorified' relations. Taking opposites $(-)\op$ forms an
involution on the category $\Cat{Cat}$ of (small) categories and
functors between them. One can think of taking homsets $\mathrm{Hom}
\colon \cat{C}\op \times \cat{C} \rightarrow \Sets$ as a comparison
relation in $\Cat{Cat}$. The tame relations discussed in the next
section then correspond to adjunctions. In order to obtain a symmetric
comparison relation we need to replace $\Sets$ by a self-dual
category, like the category $\Cat{Rel}$ of sets and relations, or a
groupoid.
\end{remark}

\section{Tame relations}\label{TameRelSec}

This section introduces the setting in which one can define tameness
for relations, leading to the first result, namely that such tame
relations give rise to a dagger category
(Proposition~\ref{TameRelCatProp}).

\begin{definition}
\label{ComparClusterDef}
A \textbf{comparison cluster} consists of a collection
$\smash{\big(\overline{X_{i}}\otimes
  X_{i}\stackrel{\cp[i]}{\longrightarrow} \Omega\big)_{i}}$ of
comparison maps in an involutive monoidal category (with a shared
target object $\Omega$). This cluster is called \textbf{symmetric} if
each of the comparison relations $\cp[i]$ is symmetric.
\end{definition}

In the category \Sets each object $X$ carries equality $=$ as a
comparison relation $X\times X\rightarrow 2$. But there also
situations where only specific objects in a category carry such a
relation. For instance, in the category \JSL of (finite) join
semilattices the free objects carry such comparisons. Recall that free
semilattices are given by finite powersets $\Powfin(X) =
\set{U\subseteq X}{U\mbox{ is finite}}$.  They carry a comparison
relation $\cp[X] \colon \Powfin(X)\otimes\Powfin(X) \rightarrow 2$ in
\JSL, where $\cp[X](U\sotimes V) = 1$ iff $U\cap V \neq
\emptyset$. The tensor $\otimes$ in \JSL arises because of
bilinearity: $\emptyset\, \cap V \neq \emptyset$ never holds, and
$(U_{1}\cup U_{2}) \cap V \neq \emptyset$ iff either $U_{1}\cap V \neq
\emptyset$ or $U_{2}\cap V \neq \emptyset$. Hence these $\cp[X] \colon
\Powfin(X)\otimes\Powfin(X) \rightarrow 2$ form a comparison cluster
in \JSL, indexed by sets $X$.

Similarly, we may consider the collection of Hilbert spaces with their
inner products $\big(\overline{H}\otimes H
\xrightarrow{\inprod{-}{-}_H} \mathbb{C}\big)_{H\in\Hilb}$ as a
comparison cluster in the category $\Vect_{\mathbb{C}}$ of vector
spaces over $\mathbb{C}$.

More formally, we understand the index elements $i$ in
Definition~\ref{ComparClusterDef} as objects of a discrete category (no
arrows except identities). The mapping $i\mapsto X_{i}$ then forms a
functor, like the finite powerset $\Powfin$ above. We do not need
morphisms between these index elements. This functorial view is
sometimes convenient, so we may describe a comparison cluster in a
category $\cat{A}$ as a collection $\big(\overline{F(X)}\otimes F(X)
\xrightarrow{\cp[X]} \Omega\big)_{X\in\cat{D}}$, where $\cat{D}$ is a
discrete category and $F\colon \cat{D}\rightarrow\cat{A}$ is a
functor.

\begin{definition}
\label{TameRelDef}
Assume a comparison cluster $\smash{\big(\overline{F(X)}\otimes F(X)
  \xrightarrow{\cp[X]} \Omega\big)_{X\in\cat{D}}}$, as described
above, in an involutive monoidal category $\cat{A}$. A map in $\cat{A}$
of the form $r\colon \overline{F(X)}\otimes F(Y) \rightarrow \Omega$
is called a relation. Such a relation is called \textbf{tame} if there
are necessarily unique maps $r_{*}\colon F(X)\rightarrow F(Y)$ and
$r^{*}\colon F(Y)\rightarrow F(X)$ in $\cat{A}$ for which the
following diagram commutes.
\begin{equation}
\label{TameRelDiag}
\vcenter{\xymatrix@R1.5pc@C+1pc{
\overline{F(Y)}\otimes F(Y)\ar[dr]_{\cp[Y]} & 
   \overline{F(X)}\otimes F(Y)\ar[d]_{r}
   \ar[l]_-{\overline{r_{*}}\otimes\idmap}\ar[r]^-{\idmap\otimes r^{*}} & 
   \overline{F(X)}\otimes F(X)\ar[dl]^{\cp[X]} \\
& \Omega &
}}
\end{equation}

\noindent Given another (tame) relation $s\colon
\overline{F(Y)}\otimes F(Z)\rightarrow \Omega$ we define a composition
$s\relafter r$ as:
$$\xymatrix@R1.5pc@C+1pc{
s\relafter r = \Big(\overline{F(X)}\otimes F(Z)
   \ar[r]^-{\overline{r_{*}}\otimes\idmap} &
   \overline{F(Y)}\otimes F(Z)\ar[r]^-{s} & \Omega\Big).
}$$
\end{definition}

Notice that if $r$ is a tame relation, $r_{*}$ and $r^{*}$ determine
each other: $r_{*}$ determines $r$, as $r = \cp \after
(\overline{r_{*}}\otimes\idmap)$, and thus $r^{*}$ via the
mono-property of $\cp$. As we shall see in the examples below,
commutation of the triangles~\eqref{TameRelDiag} amounts to an
adjointness correspondence.

Recall the symmetric comparison cluster $\big(X\times X\xrightarrow{=}
2\big)_{X\in\Sets}$ given by equality. A relation $r\colon X\times Y
\rightarrow 2$ is tame (wrt.\ this cluster) if there are functions
$r_{*}\colon X\rightarrow Y$ and $r^{*}\colon Y\rightarrow X$ such
that, for all $x\in X$ and $y\in Y$,
$$\begin{array}{rcccl}
r_{*}(x) = y
& \Longleftrightarrow &
r(x,y) = 1
& \Longleftrightarrow &
x = r^{*}(y).
\end{array}$$

\noindent This means that $r_{*}$ and $r^{*}$ are each other's
inverses. One can interpret this as: set-theoretic reversible
computation requires isomorphisms (bijections).

Before we can form a category of tame relations, we need the following
results.

\begin{lemma}
\label{TameRelCompLem}
In the context of the previous definition:
\begin{enumerate}
\item comparison relations are tame, with $(\cp[X])_{*} = \idmap[F(X)]
  = (\cp[X])^{*}$;

\item for tame relations $r\colon F(X)\otimes F(Y)\rightarrow\Omega$
  and $s\colon F(Y)\otimes F(Z)\rightarrow\Omega$, the relation
  composition $s\relafter r$ is tame, with $(s\relafter r)_{*} = s_{*}
  \after r_{*}$ and $(s\relafter r)^{*} = r^{*} \after s^{*}$.
\end{enumerate}
\end{lemma}

\begin{myproof}
The first point is immediate, and for the second point we show
that the maps $s_{*} \after r_{*}$ and $r^{*} \after s^{*}$ satisfy
the appropriate equations, making the relation $s\relafter r$ tame:
$$\begin{array}[b]{rcl}
\cp \after (\overline{(s_{*}\after r_{*})}\otimes\idmap) 
& = &
\cp \after (\overline{s_{*}}\otimes\idmap) \after 
   (\overline{r_{*}}\otimes\idmap) \\
& = &
s \after (\overline{r_{*}}\otimes\idmap) \\
& = &
s \relafter r \\
\cp \after (\idmap\otimes (r^{*} \after s^{*}))
& = &
\cp \after (\idmap\otimes r^{*}) \after (\idmap\otimes s^{*}) \\
& = &
r \after (\idmap\otimes s^{*}) \\
& = &
\cp \after (\overline{r_{*}}\otimes\idmap) \after (\idmap\otimes s^{*}) \\
& = &
\cp \after (\idmap\otimes s^{*}) \after (\overline{r_{*}}\otimes\idmap) \\
& = &
s \after (\overline{r_{*}}\otimes\idmap) \\
& = &
s \relafter r.
\end{array}\eqno{\QEDbox}$$
\end{myproof}

The comparison cluster $\big(X\op\times X\xrightarrow{\leq}
  2\big)_{X\in\PoSets}$ from the previous section is non-symmetric. A
relation $r\colon X\times Y \rightarrow 2$ in \PoSets is tame if there
are monotone functions $r_{*}\colon X \rightarrow Y$ and $r^{*}\colon
Y\rightarrow X$ such that for all $x\in X$ and $y\in Y$,
$$\begin{array}{rcccl}
r_{*}(x) \leq y
& \Longleftrightarrow &
r(x,y) = 1
& \Longleftrightarrow &
x \leq r^{*}(y).
\end{array}$$

\noindent Thus a tame relation $r$ comes from a Galois connection.  As
is well-known, Galois connections are closed under composition, in the
obvious manner. But exchanging $r_{*}$ and $r^{*}$ does (in general)
not yield another Galois connection---but see
Subsection~\ref{OMSubsec} for a remedy for orthomodular lattices. In
the next result we shall use symmetry to obtain such reversals, in the
form of daggers.

\begin{lemma}
\label{TameRelDagLem}
For a tame relation $r\colon X\otimes Y\rightarrow \Omega$ we define a
swapped version:
$$\xymatrix@R1.5pc@C+.5pc{
r^{\dag} = \Big(\overline{F(Y)}\otimes F(X)
   \ar[r]^-{\overline{r^{*}}\otimes\idmap} &
   \overline{F(X)}\otimes F(X)\ar[r]^-{\cp} & \Omega\Big).
}$$

\noindent Assuming that the comparison cluster is symmetric, we get:
\begin{enumerate}
\item $r^{\dag}$ is the same as the composite:
$$\xymatrix@R1.5pc{
\overline{F(Y)}\otimes F(X)\ar[r]^-{\tau^{-1}}_-{\cong} &
   \overline{\overline{F(X)}\otimes F(Y)}\ar[r]^-{\overline{r}} &
   \overline{\Omega}\ar[r]^-{j}_-{\cong} & \Omega,
}$$

\noindent where $\tau$ is the twist map from~\eqref{TwistConjugEqn};

\item $(r^{\dag})_{*} = r^{*}$ and $(r^{\dag})^{*} = r_{*}$, making
also $r^{\dag}$ tame.
\end{enumerate}
\end{lemma}

\begin{myproof}
For the first point we obtain, by symmetry~\eqref{ComparSymmEqn}:
$$\begin{array}[b]{rcl}
j \after \overline{r} \after \tau^{-1}
& = &
j \after \overline{\cp} \after \overline{(\idmap\otimes r^{*})} 
   \after \tau^{-1} \\
& = &
j \after \overline{\cp} \after \tau^{-1} \after 
   (\overline{r^{*}}\otimes\idmap) 
   \qquad \mbox{by naturality of }\tau \\
& \smash{\stackrel{\eqref{ComparSymmEqn}}{=}} &
\cp \after (\overline{r^{*}}\otimes\idmap) \\
& = &
r^{\dag}.
\end{array}$$

\noindent By construction of $r^{\dag}$ as $\cp \after
(\overline{r^{*}}\otimes\idmap)$, the map $r^{*}$ plays the role of
$(r^{\dag})_{*}$. It is easy to see that $r_{*}$ acts as
$(r^{\dag})^{*}$:
$$\begin{array}[b]{rcl}
\cp \after (\idmap\otimes r_{*})
& \smash{\stackrel{\eqref{ComparSymmEqn}}{=}} &
j \after \overline{\cp} \after \tau^{-1} \after 
   (\overline{\idmap}\otimes r_{*}) \\
& = &
j \after \overline{\cp} \after 
   \overline{(\overline{r_{*}}\otimes\idmap)} \after \tau^{-1} \\
& = &
j \after \overline{r} \after \tau^{-1} \\
& = &
r^{\dag}, \qquad \mbox{as just shown.}
\end{array}\eqno{\QEDbox}$$
\end{myproof}

We summarise the situation.

\begin{proposition}
\label{TameRelCatProp}
A comparison cluster $\smash{\cp = \big(\overline{F(X)}\otimes F(X)
  \xrightarrow{\cp[X]} \Omega\big)_{X}}$ in a category $\cat{A}$ gives
rise a category $\TRel(\cat{A},\cp)$ of tame relations; it has indices
$X$ as objects, and its morphisms $X\rightarrow Y$ are tame relations
$\overline{F(X)}\otimes F(Y)\rightarrow \Omega$. Comparison relations
$\cp[X]$ form identity maps on $X$, and composition is given by
$\relafter$, as in Definition~\ref{TameRelDef}.

In case the comparison cluster is symmetric, $\TRel(\cat{A},\cp)$ is a
dagger category, with dagger $(-)^\dag$ as in
Lemma~\ref{TameRelDagLem}.
\end{proposition}

\begin{myproof}
We briefly check the basic properties, using Lemma~\ref{TameRelCompLem}
and~\ref{TameRelDagLem}.
$$\begin{array}[b]{rclcrcl}
\cp \relafter r
& = &
\cp \after (\overline{r_{*}}\otimes\idmap)
& &
\cp^{\dag}
& = &
\cp \after (\overline{\cp^{*}}\otimes\idmap) \\
& = &
r
& &
& = &
\cp \after (\overline{\idmap}\otimes\idmap) \\
s \relafter \cp
& = &
s \after (\overline{{\cp}_{*}}\otimes\idmap)
& &
& = &
\cp \\
& = &
s \after (\overline{\idmap}\otimes\idmap)
& &
(s\relafter r)^{\dag}
& = &
\cp \after (\overline{(s\relafter r)^{*}}\otimes\idmap) \\
& = &
s 
& &
& = &
\cp \after (\overline{r^{*}}\otimes\idmap) \after 
   \rlap{$(\overline{s^{*}}\otimes\idmap)$}\qquad \\
t \relafter (s \relafter r)
& = &
t \after (\overline{(s\relafter r)_{*}}\otimes\idmap)
& &
& = &
r^{\dag} \after (\overline{(s^{\dag})_{*}}\otimes\idmap) \\
& = &
t \after (\overline{(s_{*}\after r_{*})}\otimes\idmap)
& &
& = &
r^{\dag} \relafter s^{\dag} \\
& = &
t \after (\overline{s_{*}}\otimes\idmap) \after 
   (\overline{r_{*}}\otimes\idmap)
& &
r^{\dag\dag}
& = &
\cp \after (\overline{(r^{\dag})^{*}}\otimes\idmap) \\
& = &
(t \relafter s) \after (\overline{r_{*}}\otimes\idmap)
& &
& = &
\cp \after (\overline{r_{*}}\otimes\idmap) \\
& = &
(t \relafter s) \relafter r
& &
& = &
r.
\end{array}\eqno{\QEDbox}$$
\end{myproof}

In the sequel we focus on symmetric comparison clusters.
We end this section with some easy but useful observation.

\begin{lemma}
\label{TameRelDagMonoEpiLem}
For a map $r\colon X\rightarrow Y$ in the dagger category
$\TRel(\cat{A},\cp)$ of a symmetric comparison cluster one has:
$$\begin{array}{rcl}
r\mbox{ is a dagger mono, \textit{i.e.} } r^{\dag} \relafter r = \idmap
& \Longleftrightarrow &
r^{*} \after r_{*} = \idmap \\
r\mbox{ is a dagger epi, \textit{i.e.} } r \relafter r^{\dag} = \idmap
& \Longleftrightarrow &
r_{*} \after r^{*} = \idmap.
\end{array}$$

\noindent As a result we can characterise dagger isomorphisms (or:
unitary maps) as:
$$\begin{array}{rcccccl}
r\mbox{ is a dagger iso} 
& \smash{\stackrel{\textrm{def}}{\Longleftrightarrow}} &
r^{\dag} = r^{-1}  
& \Longleftrightarrow &
\left\{\begin{array}{rcl}
r_{*} \after r^{*} & = & \idmap \\
r^{*} \after r_{*} & = & \idmap
\end{array}\right. 
& \Longleftrightarrow &
\begin{array}{rcl} 
(r^{\dag})_{*} & = & (r_{*})^{-1}.
\end{array}
\end{array}$$
\end{lemma}

\begin{myproof}
Assume $r^{\dag} \relafter r = \idmap$. Then, using
Lemma~\ref{TameRelCompLem} and~\ref{TameRelDagLem}, $r^{*} \after
r_{*} = r^{*} \after (r^{\dag})^{*} = (r^{\dag} \relafter r)^{*} =
{\cp}^{*} = \idmap$. Conversely, if $r^{*} \after r_{*} = \idmap$,
then $r^{\dag} \relafter r = r^{\dag} \after (r_{*}\otimes\idmap) =
\cp \after (r^{*}\otimes\idmap) \after (r_{*}\otimes\idmap) = \cp =
\idmap$. The dagger epi case is handled similarly, and the result for
dagger isos follows by combining these two cases. \QED

\auxproof{ 
Just to be sure that there are no surprises.  Assume $r \relafter
r^{\dag} = \idmap$. Then, using Lemma~\ref{TameRelCompLem}
and~\ref{TameRelDagLem}, $r_{*} \after r^{*} = r_{*} \after
(r^{\dag})_{*} = (r \relafter r^{\dag})_{*} = {\cp}_{*} =
\idmap$. Conversely, if $r_{*} \after r^{*} = \idmap$, then $r
\relafter r^{\dag} = r \after ((r^{\dag})_{*}\otimes\idmap) = r \after
(r^{*}\otimes\idmap) = \cp \after (r_{*}\otimes\idmap) \after
(r^{*}\otimes\idmap) = \cp = \idmap$.  
}
\end{myproof}

Later on, in Section~\ref{BfMRelSec}, we shall see examples of
dagger monos in a category of tame relation (see especially
Lemma~\ref{BfMRelDagMonoLem}).

\begin{lemma}
In the same context as the previous lemma, an endomap $r\colon
X\rightarrow X$ is self-adjoint (\textit{i.e.}~$r^{\dag} = r$) iff
$r_{*} = r^{*}$.

It is a projection (\textit{i.e.}~$r \relafter r = r = r^{\dag}$) iff
$r_{*} = r^{*}$ and $r_{*} \after r_{*} = r_{*}$.
\end{lemma}

\begin{myproof}
If $r = r^{\dag}$ then $r_{*} = (r^{\dag})_{*} = r^{*}$. Conversely,
if $r_{*} = r^{*}$, then $r^{\dag} = \cp \after
(\overline{r^{*}}\otimes \idmap) = \cp \after (\overline{r_{*}}\otimes
\idmap) = r$.

If $r$ is a projection, then it is a self-adjoint and so $r_{*} =
r^{*}$. Further, $\cp \after (\overline{r_{*}}\otimes \idmap) = r = r
\relafter r = r \after (\overline{r_{*}}\otimes \idmap) = \cp \after
(\overline{r_{*}}\otimes \idmap) \after (\overline{r_{*}}\otimes
\idmap) = \cp \after (\overline{r_{*}\after r_{*}}\otimes \idmap)$.
Hence $\overline{r_{*}} = \overline{r_{*} \after r_{*}}$, by the
mono-requirement for $\cp$, and thus $r_{*} = r_{*} \after r_{*}$.
The converse is obvious. \QED
\end{myproof}

%% \begin{lemma}
%% Assume a symmetric comparison cluster $\smash{\cp =
%%   \big(\overline{F(X)}\otimes F(X) \xrightarrow{\cp[X]}
%%   \Omega\big)_{X}}$ in a category $\cat{A}$. An isomorphism
%% $\varphi\colon F(X) \rightarrow F(Y)$ in $\cat{A}$ gives rise to a
%% dagger isomorphism $\cp \after (\varphi\otimes\idmap)\colon
%% X\congrightarrow Y$ in $\TRel(\cat{A},\cp)$.
%% \end{lemma}

%% \begin{proof}
%% For convenience we abbreviate $\widehat{\varphi} = \cp \after
%% (\varphi\otimes\idmap)\colon X\otimes Y \rightarrow \Omega$ in
%% $\cat{A}$. It forms a tame relation with $(\widehat{\varphi})_{*} =
%% \varphi$ and $(\widehat{\varphi})^{*} = \varphi^{-1}$

%% \auxproof{
%% $$\xymatrix@R1.5pc{
%% \overline{F(Y)}\otimes F(Y)\ar[ddr]_{\cp} &
%%    \overline{F(X)}\otimes F(Y)\ar[l]_-{\overline{\varphi}\otimes\idmap}
%%       \ar[r]^-{\idmap\otimes\varphi^{-1}}
%%       \ar[d]^{\overline{\varphi}\otimes\idmap} &
%%    \overline{F(X)}\otimes F(X)\ar[ddl]^{\cp} \\
%% & \overline{F(Y)}\otimes F(Y)\ar[d]^{\cp} \\
%% & \Omega
%% }$$
%% }

%% \end{proof}

\section{Examples of categories of tame relations}\label{ExSec}

All the illustrations of comparison clusters in this section will be
symmetric---resulting in dagger categories of tame relations. In many
of the examples below we have closed structure---with an exponent
$\multimap$ for $\otimes$. Thus we can equivalently describe such
relations $\overline{F(X)}\otimes F(Y) \rightarrow \Omega$ as maps
$F(Y)\rightarrow \big(\overline{F(X)}\multimap\Omega\big)$. This is
often more convenient, since it avoids tensors.

\subsection{Orthomodular lattices and Galois connections}\label{OMSubsec}

In Subsection~\ref{ComparRelSubsec} we have seen that the order on a
poset $X$ forms a non-symmetric comparison relation $\leq \colon
\overline{X}\times X\rightarrow 2$ in \PoSets, where $\overline{(-)}$
is order-reversal. Now assume that $X$ is an orthomodular lattice
(see~\cite{Kalmbach83} for details), with orthocomplement $(-)^{\perp}
\colon X\rightarrow \overline{X}$. It satisfies, among other things,
$x^{\perp\perp} = x$ and: $x^{\perp} \leq y$ iff $y^{\perp} \leq x$.
When $x \leq y^{\perp}$ one calls $x,y$ orthogonal, which is also
written as $x\orthogonal y$. We obtain a comparison relation
$\cp[\perp] \colon X\times X\rightarrow 2$ in \PoSets (with identity
involution), via $\cp[\perp](x,y)=1$ iff $x^{\perp} \leq y$. By using
orthocomplement in the first coordinate the contravariance
disappears. This relation is the same as $(x,y) \mapsto x^{\perp}
\orthogonal y^{\perp}$, that is, as orthogonality of
orthocomplements. It forms a symmetric comparison relation, since
orthogonality is symmetric. The resulting category of tame relations
is known from~\cite{Crown75,Jacobs10c}.

\auxproof{
This $\cp[\perp]$ is a map in \PoSets, since if $x\leq x'$ and 
$y\leq y'$, then:
$$\cp[\perp](x,y) = 1
\Rightarrow 
x^{\perp} \leq y
\Rightarrow
x'^{\perp} \leq x^{\perp} \leq y \leq y'
\Rightarrow
\cp[\perp](x',y') = 1.$$

The mono-requirement holds: assume $\cp[\perp](x,z)$ iff $\cp[\perp](y,z)$,
for all $z$. This means $x^{\perp} \leq z$ iff $y^{\perp} \leq z$,
and thus $x^{\perp} = y^{\perp}$; hence $x=y$.
}

\begin{proposition}
\label{OMProp}
The category of tame relations $\TRel(\PoSets, \cp[\perp])$ for the
symmetric comparison cluster $\big(X\times X\xrightarrow{\cp[\perp]}
2\big)_{X\in\Cat{OrthMod}}$ given by orthogonality of
orthocomplements, is the category \Cat{OMLatGal} of orthomodular
lattices and (antitone) Galois connections between them.
\end{proposition}

\begin{myproof}
A tame relations $r\colon X\rightarrow Y$, for $X,Y$ orthomodular
lattices is determined by monotone functions $r_{*}\colon X\rightarrow Y$
and $r^{*}\colon Y\rightarrow X$ satisfying:
$$\begin{array}{rcl}
r_{*}(x)^{\perp} \leq y
\hspace*{\arraycolsep} \Longleftrightarrow \hspace*{\arraycolsep}
\cp[\perp](r_{*}(x), y) = 1
& \Longleftrightarrow &
r(x,y) = 1 \\
& \Longleftrightarrow &
\cp[\perp](x, r^{*}(y))
\hspace*{\arraycolsep} \Longleftrightarrow \hspace*{\arraycolsep}
x^{\perp} \leq r^{*}(y).
\end{array}$$

\noindent These $r_{*}$ and $r^{*}$ are completely determined by
monotone functions $r_{\#} = r_{*} \after (-)^{\perp} \colon
\overline{X} \rightarrow Y$ and $r^{\#} = r^{*} \after (-)^{\perp}
\colon Y \rightarrow \overline{X}$ satisfying:
$$\begin{array}{rcccccl}
x = x^{\perp\perp} \leq r^{\#}(y) = r^{*}(y^{\perp})
& \Longleftrightarrow &
r_{*}(x^{\perp})^{\perp} \leq y^{\perp}
& \Longleftrightarrow &
y \leq r_{*}(x^{\perp}) = r_{\#}(x).
\end{array}$$

\noindent This precisely says that $r_{\#}, r^{\#}$ form an antitone
Galois connection---or an adjunction $r^{\#} \dashv r_{\#}$. \QED

\auxproof{
Just to be sure, we also check the other direction, namely given
an antitone Galois connection $r_{\#}, r^{\#}$, defining
$r_{*} \after (-)^{\perp}$ and $r^{\#} = r^{*} \after (-)^{\perp}$
yields a tame relation, since:
$$\begin{array}{rcl}
\cp[\perp](r_{*}(x),y) = \cp[\perp](r_{\#}(x^{\perp}), y) = 1
& \Longleftrightarrow &
r_{\#}(x^{\perp})^{\perp} \leq y \\
& \Longleftrightarrow &
y^{\perp} \leq r_{\#}(x^{\perp}) \\
& \Longleftrightarrow &
x^{\perp} \leq r^{\#}(y^{\perp}) = r^{*}(y) \\
& \Longleftrightarrow &
\cp[\perp](x, r^{*}(y)) = 1.
\end{array}$$
}
\end{myproof}

In~\cite{Jacobs10c} it is shown that $\Cat{OMLatGal}$ is a dagger
kernel category with (dagger) biproducts, and that every dagger
kernel category maps into it.

\subsection{Locally bifinite relations and partial injections}

We have already seen the finite powerset $\Powfin(X) = \set{U\subseteq
  X}{U\mbox{ is finite}}$ as free functor $\Powfin\colon \Sets
\rightarrow \JSL$, left adjoint to the forgetful functor from the
category of join semi-lattices (finite joins only).  This category
$\JSL$ is in fact the category of (Eilenberg-Moore) algebras of the
commutative (symmetric monoidal) monad $\Powfin$. Hence \JSL is
symmetric monoidal closed, following the constructions
in~\cite{Kock71a,Kock71b}, where $\Powfin$ preserves the monoidal
structure: $\Powfin(1) = 2$ is tensor unit and $\Powfin(X\times Y)
\cong \Powfin(X)\otimes\Powfin(Y)$. We first review the comparison
structure on free semilattices $\Powfin(X)$, with respect to the
trivial (identity) involution on \JSL.

As $\Omega\in\JSL$ we take the two-element lattice $2 =
\Powfin(1)$. Then we have correspondences between `abstract' relations
and ordinary relations, in:
\begin{equation}
\label{PowfinRelEquiv}
\begin{prooftree}
\begin{prooftree}
\begin{prooftree}
X\times Y\longrightarrow 2 \rlap{\hspace*{5em} in \Sets} 
\Justifies
\Powfin(X\times Y)\longrightarrow 2 \rlap{\hspace*{3.8em} in \JSL} 
\end{prooftree}
\Justifies
\Powfin(X)\otimes \Powfin(Y) \longrightarrow 2 
   \rlap{\hspace*{2.6em} in \JSL} 
\end{prooftree}
\Justifies
\Powfin(Y)\longrightarrow \big(\Powfin(X) \multimap 2\big)
   \rlap{\hspace*{2em} in \JSL} 
\end{prooftree}
\end{equation}

\noindent Starting from the equality relation $=\colon X\times X
\rightarrow 2$ in \Sets this correspondence yields a comparison
relation $\cpPow\colon \Powfin(X) \rightarrow (\Powfin(X)\multimap 2)$
given by:
\begin{equation}
\label{PowfinComparEqn}
\begin{array}{rcccl}
\cpPow(U)(U')
& = &
\displaystyle\bigvee_{(x,x')\in U\times U'}(x=x')
& = &
\left\{\begin{array}{ll}
1\; & \mbox{if }U\cap U'\neq \emptyset \\
0 & \mbox{otherwise.}
\end{array}\right.
\end{array}
\end{equation}

\noindent Clearly, this relation $\cpPow$ is symmetric; it is also
monic: if $\cpPow(U) = \cpPow(V)$, then:
$$x\in U
\Longleftrightarrow
\cpPow(U)(\{x\}) = 1
\Longleftrightarrow
\cpPow(V)(\{x\}) = 1
\Longleftrightarrow
x\in V.$$

\noindent Hence $U=V$.

\begin{proposition}
\label{LocBifinProp}
The dagger category $\TRel(\JSL,\cpPow)$ of tame relations for the
symmetric comparison cluster $\cpPow\colon
\Powfin(X)\otimes\Powfin(X)\rightarrow 2$ determined
by~\eqref{PowfinComparEqn} is the category of sets with bifinite
relations between them, \textit{i.e.}~with those relations $r\subseteq
X\times Y$ where for each $x\in X$ and $y\in Y$ both the sets
$$\setin{z}{Y}{r(x,z)}
\qquad\mbox{and}\qquad
\setin{w}{X}{r(w,y)}$$

\noindent are finite. Such a relation factors in two directions as
$X\rightarrow\Powfin(Y)$ and as $Y\rightarrow\Powfin(X)$. Thus we also
write $\BifRel = \TRel(\JSL,\cpPow)$ for this category of sets and
bifinite relations.
\end{proposition}

\begin{myproof}
Assume $r\subseteq X\times Y$, which corresponds to $\widehat{r}\colon
\Powfin(Y) \rightarrow (\Powfin(X)\multimap 2)$ in \JSL like
in~\eqref{PowfinRelEquiv}, given by $\widehat{r}(V)(U) = 1$ iff
$r(x,y)$ holds for some $x\in U$ and $y\in V$. We shall prove the
equivalence of:
\begin{enumerate}
\item[(a)] for each $y\in Y$, the set $\set{x}{R(x,y)} \subseteq X$
is finite;

\item[(b)] there is a necessarily unique map $r^{*}\colon
\Powfin(Y)\rightarrow\Powfin(X)$ in \JSL in the diagram:
$$\xymatrix@R1.5pc@C-1pc{
\Powfin(Y)\ar[rr]^-{\widehat{r}}\ar@{-->}[dr]_{r^{*}} & & 
   \Powfin(X)\multimap 2 \\
& \Powfin(X)\ar@{ >->}[ur]_{\cpPow}
}$$
\end{enumerate}

\noindent This diagram corresponds to the triangle on the right
in~\eqref{TameRelDiag}. The analogous statement for $r_{*}$ is left to
the reader.

So assume~(a) holds. Then we can define $r^{*}(V) \in \Powfin(X)$, for
$V\in\Powfin(Y)$, as finite union of finite sets, namely as $r^{*}(V)
=\bigcup_{y\in V}\set{x}{r(x,y)}$. It makes the triangle in~(b)
commute:
$$\begin{array}{rcl}
\big(\cpPow \after r^{*}\big)(V)(U) = 1
& \Longleftrightarrow &
\cpPow(r^{*}(V))(U) = 1 \\
& \Longleftrightarrow &
U \cap r^{*}(V) \neq \emptyset \\
& \Longleftrightarrow &
\exin{x}{U}{\exin{y}{V}{r(x,y)}} \\
& \Longleftrightarrow &
\widehat{r}(V)(U) = 1.
\end{array}$$

Conversely, assume~(b) holds, so that we have a map $r^{*}\colon
\Powfin(Y)\rightarrow \Powfin(X)$ in \JSL in the above triangle. Then:
$$\begin{array}{rcl}
r(x,y)
& \Longleftrightarrow &
\widehat{r}(\{y\})(\{x\}) = 1 \\
& \Longleftrightarrow &
\cpPow(r^{*}(\{y\}))(\{x\}) = 1 \\
& \Longleftrightarrow &
\{x\} \cap r^{*}(\{y\}) \neq \emptyset \\
& \Longleftrightarrow &
x\in r^{*}(\{y\}).
\end{array}$$

\noindent Since $r^{*}(\{y\})\in\Powfin(X)$ there are at most
finitely many $x$ that satisfy $R(x,y)$.

\auxproof{
\begin{enumerate}
\item[(a)] for each $x\in X$, the set $\set{y}{r(x,y)} \subseteq Y$
is finite;

\item[(b)] there is a necessarily unique map $r_{*}\colon
\Powfin(X)\rightarrow\Powfin(Y)$ in \JSL in the diagram:
$$\xymatrix@R1.5pc@C-1pc{
\Powfin(Y)\ar[rr]^-{\widehat{r}}\ar@{ >->}[dr]_{\cpPow} & & 
   \Powfin(X)\multimap 2 \\
& \Powfin(Y)\multimap 2\ar[ur]_{\qquad r_{*}\multimap 2 = (-)\after r_{*}}
}$$
\end{enumerate}

\noindent This diagram corresponds to the triangle on the left
in~\eqref{TameRelDiag}. 

So assume~(a) holds. Then we can define $r_{*}(U) \in \Powfin(Y)$, for
$U\in\Powfin(X)$, as finite union of finite sets, namely as $r_{*}(U)
=\bigcup_{x\in U}\set{y}{r(x,y)}$. It makes the triangle in~(b)
commute:
$$\begin{array}{rcl}
\big((r_{*}\multimap 2) \after \cpPow\big)(V)(U) = 1
& \Longleftrightarrow &
\big(\cpPow(V) \after r_{*}\big)(U) = 1 \\
& \Longleftrightarrow &
\cpPow(V)(r_{*}(U)) = 1 \\
& \Longleftrightarrow &
r_{*}(U) \cap V \neq \emptyset \\
& \Longleftrightarrow &
\exin{x}{U}{\exin{y}{V}{r(x,y)}} \\
& \Longleftrightarrow &
\widehat{r}(V)(U) = 1.
\end{array}$$

Conversely, assume~(b) holds, so that we have a map $r_{*}\colon
\Powfin(X)\rightarrow \Powfin(Y)$ in \JSL in the above triangle. Then:
$$\begin{array}{rcl}
r(x,y)
& \Longleftrightarrow &
\widehat{r}(\{y\})(\{x\}) = 1 \\
& \Longleftrightarrow &
\cpPow(\{y\})(r_{*}(\{x\})) = 1 \\
& \Longleftrightarrow &
r_{*}(\{x\}) \cap \{y\} \neq \emptyset \\
& \Longleftrightarrow &
y\in r_{*}(\{x\}).
\end{array}$$

\noindent Since $r_{*}(\{x\})\in\Powfin(Y)$ there are at most
finitely many $y$ that satisfy $r(x,y)$. 
}

Finally, it is easy to see that composition in the category
$\BifRel = \TRel(\JSL,\cpPow)$ is just relational composition, and
that the dagger is relational converse. \QED

\auxproof{
Assume $r\subseteq X\times Y$ and $s\subseteq Y\times Z$ are
considered as maps in $\TRel(\JSL,\cpPow)$. Their composition, 
as such relations is:
$$\begin{array}{rcl}
(\widehat{s} \relafter \widehat{r})(\{z\})(\{x\}) = 1
& \Longleftrightarrow &
\widehat{s}(\{z\})(\widehat{r}(\{x\})) = 1 \\
& \Longleftrightarrow &
\exin{y}{\widehat{r}(\{x\})}{s(y,z)} \\
& \Longleftrightarrow &
\ex{y}{r(x,y) \conjun s(y,z)}.
\end{array}$$

Similarly,
$$\begin{array}{rcl}
\widehat{r}^{\dag}(\{x\})(\{y\}) = 1
& \Longleftrightarrow &
\cpPow(\{y\})(r^{*}(\{x\})) = 1 \\
& \Longleftrightarrow &
x \in r^{*}(\{y\}) \\
& \Longleftrightarrow &
<<<<<<< tame-relations.tex
\cp(R^{*}(\{y\}))(\{x\}) = 1 \\
=======
r(x,y) \\
>>>>>>> 1.33
& \Longleftrightarrow &
y \in r_{*}(\{x\}) \\
& \Longleftrightarrow &
\widehat{r}(\{y\})(\{x\}) = 1.
\end{array}$$
}
\end{myproof}

For a map $r\colon X\times Y\rightarrow 2$, as morphism in
$\BifRel = \TRel(\JSL,\cpPow)$, the `adjointness'
correspondence~\eqref{TameRelDiag} takes the form:
$$\begin{array}{rcccl}
r_{*}(U) \cap V \neq \emptyset
& \Longleftrightarrow &
\exin{x}{U}{\exin{y}{V}{r(x,y)}}
& \Longleftrightarrow &
U \cap r^{*}(V) \neq \emptyset,
\end{array}$$

\noindent for $U\in\Powfin(X)$ and $V\in\Powfin(Y)$. Moreover, such
a map $r\colon X\rightarrow Y$ is unitary if and only it is given
by an isomorphism of sets $X\cong Y$.

\auxproof{
Assume $r$ is unitary. Then:
$$\textstyle 1
=
(r^{\dag} \relafter r)(x,x)
=
\bigvee_{y}r(x,y)\cdot r(x,y)
=
\bigvee_{y}r(x,y)$$

\noindent and if $x\neq x'$, then:
$$\textstyle 0
=
(r^{\dag} \relafter r)(x,x')
=
\bigvee_{y}r(x,y)\cdot r(x',y).$$

\noindent Hence for each $y$ there is an $x$ with $r(x,y) = 1$.
This $x$ is unique, since if $r(x,y) = 1 = r(x',y)$, then the
second equation above is violated. Similarly the other way around.
}

Our next example is fairly similar to the previous one. Below in
Subsection~\ref{MonadSubsec} we shall capture this similarity in terms
of certain monads. But we prefer to describe this second example
concretely, because it leads to a well-known category, namely the
category \PInj of sets and partial injections between them (see
\textit{e.g.}~\cite{HaghverdiS06,HeunenJ10a}). We start with the
category $\Sets_{\bullet}$ of pointed sets. Objects are sets $X$
containing a distinguished base point $\bullet\in X$. Morphisms are
ordinary functions that preserve this base point. This category
$\Sets_{\bullet}$ is equivalent to the category \Pfn of sets and
partial functions between them.

There is a ``lift'' functor $\Lft = 1+(-)\colon \Sets \rightarrow
\Sets_{\bullet}$ that adds such a base point to set; it is left
adjoint to the fogetful functor $\Sets_{\bullet} \rightarrow
\Sets$. An element $z\in\Lft(X) = 1+X$ is either of the form
$z=\bullet\in 1$ or $z=x\in X$, for a unique $x\in X$. Thus one can
see $z\in\Lft(X)$ as a subset of $X$ with at most one element (a
`subsingleton').  This category $\Sets_{\bullet}$ is the category of
algebras of $\Lft$, as monad on $\Sets$; thus, $\Sets_{\bullet}$ is
symmetric monoidal closed, following the constructions
in~\cite{Kock71a,Kock71b}. If we take $\Omega = 2 =
\Lft(1)\in\Sets_{\bullet}$, then we have a bijective correspondence
between abstract relations $\Lft(X)\otimes\Lft(Y)\rightarrow 2$ and
ordinary relations $X\times Y \rightarrow 2$ in \Sets, like
in~\eqref{PowfinRelEquiv}.

The comparison relation $\cpLft$ we use here for $\Lft$ is the same as
before---for $\Powfin$ in~\eqref{PowfinComparEqn}, when considered as
relation $=\colon X\times X\rightarrow 2$. But when we translate it
into a map $\cpLft\colon \Lft(X) \rightarrow (\Lft(X)\multimap 2)$ it
becomes:
\begin{equation}
\label{LiftComparEqn}
\begin{array}{rcl}
\cpLft(z)(z')
& = &
\left\{\begin{array}{ll}
1\; & \mbox{if $z = x = z'$ for some (necessarily unique) $x\in X$} \\
0 & \mbox{otherwise.}
\end{array}\right.
\end{array}
\end{equation}

\noindent Again this relation is symmetric, and satisfies the mono
requirement from Definition~\ref{ComparRelDef}: if $\cpLft(z) = \cpLft(w)$,
then for each $x\in X$,
$$z=x
\Longleftrightarrow
\cpLft(z)(x) = 1
\Longleftrightarrow
\cpLft(w)(x) = 1
\Longleftrightarrow
w=x.$$

\noindent Hence $z=w$.

\begin{proposition}
\label{PInjProp}
The dagger category $\TRel(\Sets_{\bullet},\cpLft)$ of tame relations for
the comparison relations~\eqref{LiftComparEqn} is the category \PInj
of sets with partial injections between them: relations $r\subseteq
X\times Y$ satisfying both:
$$\begin{array}{rclcrcl}
r(x,y) \mbox{ and } r(x,y')
& \Longrightarrow &
y=y'
& \qquad &
r(x,y) \mbox{ and } r(x',y)
& \Longrightarrow &
x=x'.
\end{array}$$

\noindent That is: $r$ factors both as $X\rightarrow \Lft(Y)$ and as
$Y\rightarrow\Lft(Y)$.
\end{proposition}

\begin{myproof}
We prove the equivalence of:
\begin{enumerate}
\item[(a)] $r(x,y) \mbox{ and } r(x',y) \Longrightarrow x=x'$;

\item[(b)] there is a necessarily unique map $r^{*}\colon
  \Lft(Y)\rightarrow\Lft(X)$ in $\Sets_\bullet$ in the diagram:
$$\xymatrix@R1.5pc@C-1pc{
\Lft(Y)\ar[rr]^-{\widehat{r}}\ar@{-->}[dr]_{r^{*}} & & 
   \Lft(X)\multimap 2 \\
& \Lft(X)\ar@{ >->}[ur]_{\cpLft}
}$$

\noindent where $\widehat{r}(z)(w)=1$ iff $w=x\in X$ and $z=y\in Y$
and $r(x,y)$.
\end{enumerate}

\noindent Assuming~(a) we define:
$$\begin{array}{rcl}
r^{*}(z)
& = &
\left\{\begin{array}{ll}
x\; & \mbox{if $z=y\in Y$ 
   and there is a (necessarily unique) $x$ with $r(x,y)$  } \\
\bullet & \mbox{otherwise.}
\end{array}\right.
\end{array}$$

\noindent Then:
$$\begin{array}{rcl}
\big(\cpLft \after r^{*})(z)(w) = 1
& \Longleftrightarrow &
\cpLft(r^{*}(z))(w) = 1 \\
& \Longleftrightarrow &
r^{*}(z) = x = w \in X \\
& \Longleftrightarrow &
z=y\in Y \mbox{ and } w=x\in X \mbox{ and }r(x,y) \\
& \Longleftrightarrow &
\widehat{r}(z)(w).
\end{array}$$

Conversely, assume $r^{*}\colon \Lft(Y)\rightarrow \Lft(X)$ as in~(b).
Then:
$$\begin{array}{rcl}
r(x,y)
& \Longleftrightarrow &
\widehat{r}(y)(x) = 1 \\
& \Longleftrightarrow &
\cpLft(r^{*}(y))(x) = 1 \\
& \Longleftrightarrow &
r^{*}(y) = x.
\end{array}$$

\noindent There is thus at most one such $x$. \QED
\end{myproof}

In the end we note that there is an obvious inclusion of categories:
$$\xymatrix@R1.5pc{
\PInj = \TRel(\Sets_{\bullet},\cpLft) \ar[r] &
   \TRel(\JSL, \cpPow) = \BifRel
}$$

\noindent since a relation $X\times Y\rightarrow 2$ that is
`bi-singlevalued' is also `bifinite'.

\subsection{Monad-based examples}\label{MonadSubsec}

The previous two examples of functors with equality arise from certain
monads, namely finite powerset $\Powfin$ and lift $\Lft$. The
constructions involved will be generalised next. Subsequently, in the
next subsection, the multiset monad $\Mlt$ and the distribution monad
$\Dst$ will be used as additional examples.

So let $\cat{B}$ be an involutive symmetric monoidal category (SMC)
carrying an involutive monad $T = (T,\eta,\mu, \sigma)$ which is
symmetric monoidal (or `commutative'), via maps $I\rightarrow T(I)$
and $T(X)\otimes T(Y)\rightarrow T(X\otimes Y)$, and with its
involution described via a distributive law $\nu_{X}\colon
T(\overline{X})\Rightarrow \overline{T(X)}$, commuting appropriately
with these two maps and with $\eta$ and $\mu$, and satisfying
$\overline{\nu} \after \nu \after T(\iota) = \iota$. These
requirements imply that $\nu$ is an isomorphism, see~\cite{Jacobs11e}
for further details.

In case the category $\Alg(T)$ of (Eilenberg-Moore) algebras has
enough coequalisers, it is also involutive symmetric monoidal, and the
free functor $F\colon\cat{B}\rightarrow\Alg(T)$ is strong
monoidal. The monoidal construction is described
in~\cite{Kock71a,Kock71b} and the involution structure
in~\cite{Jacobs11e}.  Additionally, exponents $\multimap$ in $\Alg(T)$
can be obtained from exponents in the underlying category $\cat{B}$,
via equalisers.

This situation applies to (involutive) commutative monads $T$ on
$\Sets$. The resulting category of algebras $\Alg(T)$ is always
monoidal closed.  The finite powerset $\Powfin$ and the lift monad
$\Lft$ are instances, with identity involutions; the multiset and
distribution monad form other examples below. In the rest of this
subsection we restrict to \Sets as base category.

The candidate comparison relations are defined on free objects, given
by the free functor $F\colon\Sets\rightarrow\Alg(T)$. We assume an
object $\Theta\in\Sets$ for which the free algebra $\Omega =
T(\Theta)\in\Alg(T)$ contains two different objects $0,1\in\Omega$. In
our examples it is usually obvious what these elements $0,1$ are, for
instance, for $\Omega = 2 = \{0,1\}$, or for $\Omega = [0,1]$, or for
$\Omega$ a semiring $S$, with $0$ as additive unit, and $1$ as
multiplicative unit. 

%% In all these examples, $\Omega = T(\Theta)$ is a
%% (multiplicative) commutative monoid $(\Omega, 1, \cdot)$ with zero
%% element $0\in\Omega$ satisfying $0\cdot x = 0 = x\cdot 0$.

Since we use the identity involution on \Sets there is a map
$$\xymatrix@R1.5pc@C+1pc{
\overline{\Omega} = \overline{T(\Theta)}
   \ar[r]^-{\nu^{-1}_{\Theta}}_-{\cong} &
   T(\overline{\Theta}) = T(\Theta) = \Omega
}$$

\noindent that makes this $\Omega$, like any free algebra, into a
self-conjugate object in $\Alg(T)$.

\auxproof{
We first check that this $\nu_{\Theta}^{-1}$ is self-conjugate,
using that the $\iota$ in \Sets is the identity.
$$\xymatrix@C+1pc@R1.5pc{
\overline{\overline{T(\Theta)}}\ar[r]^-{\overline{\nu^{-1}}}
   \ar[d]_{\iota^{-1}} &
   \overline{T(\overline{\Theta})}\ar[d]^{\nu^{-1}} \\
T(\Theta)\ar@{=}[r]_-{T(\iota)} & T(\overline{\overline{\Theta}})
}$$

\noindent Further, this $\nu_{\Theta}^{-1}$ is a map of algebras, from
$\overline{F(\Theta)} = \big(\overline{\mu_{\Theta}} \after \nu_{\Omega}
\colon T(\overline{\Omega}) = T(\overline{T(\Theta)}) \rightarrow
\overline{T^{2}(\Theta)}\rightarrow \overline{T(\Theta)} =
\overline{\Omega}\big)$ to $F(\Theta) = \big(\mu_{\Theta} \colon 
T(\Omega) = T^{2}(\Theta) \rightarrow T(\Theta) = \Omega\big)$, in:
$$\xymatrix@C+1pc@R1.5pc{
T(\overline{T(\Theta)})\ar[d]_{\nu_\Omega}\ar@{=}[dr]
\\
\overline{T^{2}(\Theta)}\ar[d]_{\overline{\mu}}
   \ar[r]_-{\nu^{-1}_{T(\Theta)}} & 
   T(\overline{T(\Theta)})\ar[r]_-{T(\nu^{-1}_{\Theta})} &
   T^{2}(\overline{\Theta})\ar[d]^{\mu_{\Theta}}
\\
\overline{T(\Theta)}\ar[rr]_-{\nu^{-1}_{\Theta}} & & T(\overline{\Theta})
}$$
}

In this situation we can define an equality function:
\begin{equation}
\label{EqFunEqn}
\eq[X]\colon \overline{X}\times X\longrightarrow \Omega
\qquad\mbox{by}\qquad
(x,x')\longmapsto\left\{\begin{array}{ll}
   1 \; & \mbox{if }x=x' \\ 0 & \mbox{otherwise,} \end{array}\right.
\end{equation}

\noindent where $\overline{X} = X$ is the trivial involution on \Sets
(written only for formal reasons). This equality map in $\Sets$ gives
rise to a comparison relation $\cp$ in the category $\Alg(T)$ on free
algebras, via:
$$\xymatrix@R1.5pc@C-1pc{
\cp \stackrel{\textrm{def}}{=} \Big(\overline{F(X)}\otimes F(X)
   \ar[rr]^-{\nu^{-1}\otimes\idmap}_-{\cong} & &
   F(\overline{X})\otimes F(X)\ar[r]^-{\xi}_-{\cong} & 
   F(\overline{X}\times X) \ar[r]^-{F(\eq)} & F(\Omega)\ar[r]^-{\mu} &
   \Omega\Big)
}$$

\noindent where $\mu$ is the monad's multiplication $T^{2}(\Theta)
\rightarrow T(\Theta) = \Omega$. It is not hard to see that this $\cp$
is automatically symmetric.
% (when $\nu$ is the identity, in this set-theoretic case). 
The mono-requirements from Definition~\ref{ComparRelDef} have to
checked explicitly in specific situations.

\auxproof{
For symmetry we have to check that $\cp = \mu \after T(\eq)
\after \xi^{T} \after (\nu^{-1}\otimes\idmap)$ satisfies
$\cp \after \tau = \nu^{-1}_{\Theta} \after \overline{\cp}$,
where we use the map $\nu^{-1}_{\Theta} \colon \overline{\Omega}
\rightarrow \Omega$ as self-conjugate.
{\small$$\hspace*{-3em}\xymatrix@R1.5pc{
\overline{\overline{TX}\otimes TX}\ar[r]^-{\xi^{-1}}
      \ar[d]^{\overline{\nu^{-1}\otimes\idmap}}
      \ar`u`[rrrr]^-{\tau}[rrrr]\ar`l`[ddddr]_-{\overline{\cp}}[dddd] &
   \overline{\overline{TX}}\otimes\overline{TX}
      \ar[rr]^-{\iota^{-1}\otimes\idmap}
      \ar[d]_{\overline{\nu^{-1}}\otimes\idmap} & &
   TX\otimes \overline{TX}\ar[r]^-{\gamma}
      \ar[d]_{\idmap\otimes\nu^{-1}} &
   \overline{TX}\otimes TX\ar[d]_{\nu^{-1}\otimes\idmap}
      \ar`r`[ddddl]^-{\cp}[dddd]
\\
\overline{T(\overline{X})\otimes TX}\ar[r]^-{\xi^{-1}}
      \ar[d]^{\overline{\xi^T}} &
   \overline{T(\overline{X})}\otimes\overline{TX}
      \ar[r]^-{\nu^{-1}\otimes\nu^{-1}} &
   T(\overline{\overline{X}})\otimes T(\overline{X})\ar[d]_{\xi^T}
      \ar@{=}[r]^-{T(\iota^{-1})\otimes\idmap} &
   TX\otimes T(\overline{X})\ar[r]^-{\gamma}\ar[d]_{\xi^T}
       &
   T(\overline{X})\otimes TX\ar[d]_{\xi^T}
\\
\overline{T(\overline{X}\times X)}\ar[r]^-{\nu^{-1}}
      \ar[d]^{\overline{T(\eq)}} &
   T(\overline{\overline{X}\times X})\ar@{=}[r]^-{T(\xi^{-1})}
      \ar[d]^{T(\overline{\eq})} &
   T(\overline{\overline{X}}\times\overline{X})
      \ar@{=}[r]^-{T(\iota^{-1}\times\idmap)} &
   T(X\times\overline{X})\ar[r]^-{\gamma}\ar[dr]_{T(\eq)} &
   T(\overline{X}\times X)\ar[d]_{T(\eq)}
\\ 
\overline{T^{2}(\Theta)}\ar[d]^{\overline{\mu}}\ar[r]^-{\nu^{-1}} &
   T(\overline{T(\Theta)})\ar[r]^-{T(\nu^{-1})} &
   T^{2}(\overline{\Theta})\ar[d]_{\mu}\ar@{=}[rr] & &
   T^{2}(\Theta)\ar[d]_{\mu}
\\
\overline{T(\Theta)}\ar[rr]_-{\nu^{-1}} & &
   T(\overline{\Theta})\ar@{=}[rr] & &
   T(\Theta)
}$$}

\noindent The rectangle in the middle with $\eq$ commutes if the
involution involution $\nu$ on $T$ is the identity.  
}

This general form of comparison, obtained by lifting
equality~\eqref{EqFunEqn} to a category of algebras, turns out to be
appropriate in many situations of interest. For instance, for the
finite powerset monad $\Powfin$, with $\nu=\idmap$, we get the earlier
comparison relation~\eqref{PowfinComparEqn}, since for
$U,U'\in\Powfin(X)$ this description yields:
$$\begin{array}{rcl}
\cp(U,U')
\hspace*{\arraycolsep}=\hspace*{\arraycolsep}
\bigvee\Powfin(\eq)(\xi(U,U'))
& = &
\bigvee\Powfin(\eq)(U\times U') \\
& = &
\bigvee\set{\eq(x,x')}{x\in U, x'\in U'} \\
& = &
\left\{\begin{array}{ll}
1 \; & \mbox{if }\;\exin{x}{U}{\exin{x'}{U'}{x=x'}} \\
0 & \mbox{otherwise}
\end{array}\right. \\
& = &
\left\{\begin{array}{ll}
1 \; & \mbox{if }\;U\cap U'\neq\emptyset \\
0 & \mbox{otherwise}
\end{array}\right.
\end{array}$$

\noindent The comparison relations~\eqref{LiftComparEqn} for the lift
monad $\Lft$ are also of this kind. We shall see more examples in
Subsections~\ref{MltDstSubsec} and~\ref{FormDistrSubsec} below.

In our set-theoretic examples we often take
$\Theta = 1$---but not always, see the distribution monad example
below. There are now several ways to describe `relations':
$$\begin{prooftree}
\begin{prooftree}
X\times Y \longrightarrow \Omega = T(\Theta) 
   \rlap{\hspace*{1.8em} in \Sets} 
\Justifies
\llap{$\overline{F(X)}\otimes F(Y) \;\cong\;\;$}
F(X\times Y)\longrightarrow \Omega 
   \rlap{\hspace*{2.9em} in $\Alg(T)$, by freeness} 
\end{prooftree}
\Justifies
F(Y) \longrightarrow \big(\overline{F(X)}\multimap \Omega\big)
   \rlap{\hspace*{1.5em} in $\Alg(T)$, with exponents} 
\end{prooftree}\hspace*{5em}$$

\noindent We often use such correspondences implicitly and freely
switch between different (Curry-ied or non-Curry-ied) notations for
comparison.

%% The resulting category $\TRel(\Alg(T),\cp)$ of tame
%% relations comes with two obvious functors to the Kleisli category
%% $\Kl(T)$, in a situation:
%% $$\xymatrix@R1.5pc{
%% \Kl(T) & & \TRel(\Alg(T),\cp)\ar[ll]_-{(-)_*}\ar[rr]^-{(-)^*} & & \Kl(T)\op
%% }$$

\subsection{Multiset and distribution monads}\label{MltDstSubsec}

We sketch two more applications of the monad-based construction
described above, involving the multiset monad $\Mlt_S$ and the
distribution monad $\Dst$. We shall use the multiset monad in full
generality, over a (commutative) involutive semiring $S$, like the
complex numbers $\mathbb{C}$. Such a semiring consists of a
commutative additive monoid $(S,+,0)$ and a (commutative)
multiplicative monoid $(S,\cdot,1)$, where multiplication distributes
over addition, together with an involution $\overline{(-)}\colon
S\rightarrow S$ satisfying $\overline{\overline{s}} = s$, and forming
a map of semirings. One can define a ``multiset'' functor
$\Mlt_{S}\colon\Sets\rightarrow\Sets$ by:
$$\begin{array}{rcl}
\Mlt_{S}(X)
& = &
\set{\varphi\colon X\rightarrow S}{\support(\varphi)\mbox{ is finite}},
\end{array}$$

\noindent where $\support(\varphi) = \setin{x}{X}{\varphi(x) \neq 0}$
is the support of $\varphi$. For a function $f\colon X\rightarrow Y$
one defines $\Mlt_{S}(f) \colon \Mlt_{S}(X) \rightarrow \Mlt_{S}(Y)$ by:
$$\begin{array}{rcl}
\Mlt_{S}(f)(\varphi)(y)
& = &
\sum_{x\in f^{-1}(y)}\varphi(x).
\end{array}$$

\noindent Such a multiset $\varphi\in \Mlt_{S}(X)$ may be written as
formal sum $s_{1}x_{1}+\cdots+s_{k}x_{k}$ where $\support(\varphi) =
\{x_{1}, \ldots, x_{k}\}$ and $s_{i} = \varphi(x_{i})\in S$ describes
the ``multiplicity'' of the element $x_{i}$. This formal sum notation
might suggest an order $1,2,\ldots k$ among the summands, but this sum
is considered, up-to-permutation of the summands.  Also, the same
element $x\in X$ may be counted multiple times, but $s_{1}x + s_{2}x$
is considered to be the same as $(s_{1}+s_{2})x$ within such
expressions. With this formal sum notation one can write the
application of $\Mlt_{S}$ on a map $f$ as
$\Mlt_{S}(f)(\sum_{i}s_{i}x_{i}) = \sum_{i}s_{i}f(x_{i})$.

This multiset functor is a monad, whose unit $\eta\colon X\rightarrow
\Mlt_{S}(X)$ is $\eta(x) = 1x$, and multiplication $\mu\colon
\Mlt_{S}(\Mlt_{S}(X)) \rightarrow \Mlt_{S}(X)$ is
$\mu(\sum_{i}s_{i}\varphi_{i})(x) = \sum_{i}s_{i}\cdot\varphi_{i}(x)$.
There is also an involution $\nu\colon\Mlt_{S}(X)\rightarrow\Mlt_{S}(X)$
given by $\nu(\sum_{i}s_{i}x_{i}) = \sum_{i}\overline{s_i}x_{i}$.

For the semiring $S=\NNO$ one gets the free commutative monoid
$\Mlt_{\NNO}(X)$ on a set $X$. The monad $\Mlt_{\NNO}$ is also known
as the `bag' monad, containing ordinary ($\NNO$-valued) multisets. If
$S=\mathbb{Z}$ one obtains the free Abelian group
$\Mlt_{\mathbb{Z}}(X)$ on $X$. The Boolean semiring $2 = \{0,1\}$
yields the finite powerset monad $\Powfin = \Mlt_{2}$. By taking the
complex numbers $\mathbb{C}$ as semiring one obtains the free vector
space $\Mlt_{\mathbb{C}}(X)$ on $X$ over $\mathbb{C}$.

An (Eilenberg-Moore) algebra $a\colon\Mlt_{S}(X)\rightarrow X$ for the
multiset monad corresponds to a monoid structure on $X$---given by
$x+y = a(1x + 1y)$---together with a scalar multiplication $\bullet
\colon S\times X\rightarrow X$ given by $s\mathrel{\bullet} x =
a(sx)$. It preserves the additive structure (of $S$ and of $X$) in
each coordinate separately. This makes $X$ a module, over the semiring
$S$. Conversely, such an $S$-module structure on a commutative monoid
$M$ yields an algebra $\Mlt_{S}(M)\rightarrow M$ by
$\sum_{i}s_{i}x_{i} \mapsto \sum_{i}s_{i}\mathrel{\bullet}x_{i}$.
Thus the category of algebras $\Alg(\Mlt_{S})$ is isomorphic to the
category $\Mod_{S}$ of $S$-modules. When $S$ happens to be a field,
this category $\Mod_{S}$ is the category of vector spaces $\Vect_{S}$
over $S$. It carries an involution in case $S$ is involutive,
see~\cite{Jacobs11e}.

We show that free modules $\Mlt_{S}(X)$ carry a comparison
relation. We take $\Theta = 1\in\Sets$, so that $\Omega = \Mlt_{S}(1)
= S\in\Mod_{S}$. We shall call maps $X\times Y\rightarrow S$
multirelations, in analogy with multisets; they may be seen as fuzzy
relations, assigning a possibly more general value than 0,1 to a
pair of elements. Such multirelations can thus also be described as
module maps $\Mlt_{S}(Y) \rightarrow
\big(\overline{\Mlt_{S}(X)}\multimap S\big)$, like
in~\eqref{PowfinRelEquiv}. The comparison relation, as a map
$\cp\colon \Mlt_{S}(X)\rightarrow \big(\overline{\Mlt_{S}(X)}\multimap
S\big)$ is given by (finite) sums:
\begin{equation}
\label{MltComparEqn}
\begin{array}{rcl}
\cpMlt(\varphi)(\varphi')
& = &
\sum_{x}\overline{\varphi'(x)}\cdot\varphi(x).
\end{array}
\end{equation}

\noindent This comparison captures the usual inner product (or `dot'
product) for vectors wrt.\ a basis. Symmetry amounts to
$\overline{\cpMlt(\varphi)(\varphi')} = \cpMlt(\varphi')(\varphi)$,
and thus clearly holds. In order to see that $\cpMlt$ is injective,
assume $\cpMlt(\varphi) = \cpMlt(\psi)$. Then, for each $x\in X$,
$$\varphi(x)
=
\cpMlt(\varphi)(1x)
=
\cpMlt(\psi)(1x)
=
\psi(x).$$

\noindent Hence $\varphi=\psi$, as functions $X\rightarrow S$.

The following result is no surprise anymore. The proof proceeds along
the lines of Propositions~\ref{LocBifinProp} and~\ref{PInjProp};
details are left to the interested reader.

\begin{proposition}
\label{MltTameRelProp}
Let $S$ be an involutive commutative semiring. The dagger category
$\TRel(\Mod_{S},\cpMlt)$ of tame relations for the comparison relation
$\cpMlt\colon\overline{\Mlt_{S}(X)}\otimes \Mlt_{S}(X) \rightarrow S$
corresponding to~\eqref{MltComparEqn} contains sets as objects and
morphisms $X\rightarrow Y$ are `multirelations' $r\colon X\times
Y\rightarrow S$ for which the two obvious maps obtained by abstraction
(or Curry-ing) satisfy:
$$\left\{\begin{array}{l}
\mbox{$\Lambda(r) \colon Y\longrightarrow S^{X}$ factors through
$Y\longrightarrow \Mlt_{S}(X)$} \\
\mbox{$\Lambda(r^{\dag}) \colon X\longrightarrow S^{Y}$ factors through
$X\longrightarrow \Mlt_{S}(Y)$.} 
\end{array}\right.$$

\noindent More concretely, this means that for each $x\in X$ there are
only finitely many $y\in Y$ with $r(x,y)\neq 0$, and vice-versa.

We shall also write $\BifMRel_{S} =
\TRel(\Mod_{S},\cpMlt)$ for this category of sets and bifinite
multirelations.
\end{proposition}

A bifinite multirelation $r\colon X\times Y\rightarrow S$, as morphism
$X\rightarrow Y$ in the category $\BifMRel_{S} =
\TRel(\Mod_{S},\cpMlt)$, satisfies the `adjointness' correspondence
(from~\eqref{TameRelDiag}):
$$\begin{array}{rcccl}
\cpMlt\big(r_{*}(\varphi), \psi\big)
& = &
\sum_{x,y}\overline{\varphi(x)}\cdot r(x,y) \cdot \psi(y)
& = &
\cpMlt\big(\varphi, r^{*}(\psi)\big),
\end{array}$$

\noindent for $\varphi\in\Mlt_{S}(X)$ and $\psi\in\Mlt_{S}(Y)$.
This category of bifinite multirelations will be investigated
more closely in Section~\ref{BfMRelSec}. Here we only mention that
there is an inclusion of categories:
$$\xymatrix@R1.5pc{
\BifRel = \TRel(\JSL, \cpPow) \ar[r] &
   \TRel(\Sets_{\bullet},\cpLft) = \BifMRel
}$$

\noindent since we can turn a bifinite relation $X\times Y\rightarrow
2$ into a bifinite multirelation $X\times Y\rightarrow S$ via the
inclusion $\{0,1\} \hookrightarrow S$.

\auxproof{ Assume $r\colon X\times Y\rightarrow S$, which corresponds
  to $\widehat{r}\colon \Mlt_{S}(Y) \rightarrow
  (\overline{\Mlt_{S}(X)}\multimap S)$ in $\Mod_{S}$, given by
  $\widehat{r}(\psi)(\varphi) = \sum_{x,y} \overline{\varphi(x)}\cdot
  \psi(y) \cdot r(x,y)$. We shall prove the equivalence of:
\begin{enumerate}
\item[(a)] $\Lambda(r)\colon Y\rightarrow S^{X}$ factors through
$Y\rightarrow \Mlt_{S}(X)$;

\item[(b)] there is a necessarily unique map $r^{*}\colon
\Mlt_{S}(Y)\rightarrow\Mlt_{S}(X)$ in $\Mod_{S}$ in the diagram:
$$\xymatrix@R1.5pc@C-1pc{
\Mlt_{S}(Y)\ar[rr]^-{\widehat{r}}\ar@{-->}[dr]_{r^{*}} & & 
   \overline{\Mlt_{S}(X)}\multimap S \\
& \Mlt_{S}(X)\ar@{ >->}[ur]_{\cpMlt}
}$$
\end{enumerate}

So assume~(a) holds. Then we can define $r^{*}(\psi) \in
\Mlt_{S}(X)$ as finite sum of multisets, namely as $r^{*}(\psi)
=\sum_{y}\psi(y)\cdot\Lambda(r)(y)$. It makes the triangle in~(b)
commute:
$$\begin{array}{rcl}
\big(\cpMlt \after r^{*}\big)(\psi)(\varphi)
& = &
\cpMlt(r^{*}(\psi))(\varphi) \\
& = &
\sum_{x}\overline{\varphi(x)}\cdot r^{*}(\psi)(x) \\
& = &
\sum_{x}\overline{\varphi(x)}\cdot 
   \Big(\sum_{y}\psi(y)\cdot\Lambda(r)(y)\Big)(x) \\
& = &
\sum_{x,y}\overline{\varphi(x)}\cdot \psi(y) \cdot r(x, y) \\
& = &
\widehat{r}(\psi)(\varphi).
\end{array}$$

Conversely, assume~(b) holds, so that we have a map $r^{*}\colon
\Mlt_{S}(Y)\rightarrow \Mlt_{S}(X)$ in $\Mod_{S}$ in the above
triangle. Then:
$$\begin{array}{rcl}
r(x,y)
& = &
\widehat{r}(1y)(1x) \\
& = &
\cpMlt(r^{*}(1y))(1x)  \\
& = &
r^{*}(1y)(x)
\end{array}$$

\noindent Since $r^{*}(1y)\in\Mlt_{S}(X)$ there are at most
finitely many $x$ that satisfy $r(x,y)\neq 0$.

Dually, we have:
\begin{enumerate}
\item[(a)] $\Lambda(r^{\dag})\colon X\rightarrow S^{Y}$ factors
  through $X\rightarrow \Mlt_{S}(Y)$;

\item[(b)] there is a necessarily unique map $r_{*}\colon
\Mlt_{S}(X)\rightarrow\Mlt_{S}(Y)$ in $\Mod_{S}$ in the diagram:
$$\xymatrix@R1.5pc@C-1pc{
\Mlt_{S}(Y)\ar[rr]^-{\widehat{r}}\ar[dr]_{\cpMlt} & & 
   \overline{\Mlt_{S}(X)}\multimap S \\
& \overline{\Mlt_{S}(Y)}\multimap S
      \ar[ur]_{\qquad\overline{r_*}\multimap S = (-) \after \overline{r_*}}
}$$
\end{enumerate}

So assume~(a) holds. Then we can define $r_{*}(\varphi) \in
\Mlt_{S}(Y)$ as finite sum of multisets, namely as $r_{*}(\varphi)
=\sum_{x}\varphi(x)\cdot\Lambda(r^{\dag})(x)$. It makes the
triangle in~(b) commute:
$$\begin{array}{rcl}
\big((\overline{r_*}\multimap S) \after \cpMlt\big)(\psi)(\varphi)
& = &
\big(\cpMlt(\psi) \after \overline{r_*}\big)(\varphi) \\
& = &
\cpMlt(\psi)(\overline{r_*}(\varphi)) \\
& = &
\sum_{y}\overline{\overline{r_{*}}(\varphi)(y)}\cdot \psi(y) \\
& = &
\sum_{y}\overline{\Big(\sum_{x}
   \varphi(x)\cdot\Lambda(r^{\dag})(x)\Big)(y)} \cdot \psi(y) \\
& = &
\sum_{x,y}\overline{\varphi(x) \cdot \overline{r(x, y)}} \cdot \psi(y) \\
& = &
\sum_{x,y}\overline{\varphi(x)} \cdot r(x, y) \cdot \psi(y) \\
& = &
\sum_{x,y}\overline{\varphi(x)} \cdot \psi(y) \cdot r(x, y) \\
& = &
\widehat{r}(\psi)(\varphi).
\end{array}$$

Conversely, assume~(b) holds, so that we have a map $r^{*}\colon
\Mlt_{S}(Y)\rightarrow \Mlt_{S}(X)$ in $\Mod_{S}$ in the above
triangle. Then:
$$\begin{array}{rcl}
r(x,y)
& = &
\widehat{r}(1y)(1x) \\
& = &
\cpMlt(1y)(r_{*}(1x))  \\
& = &
r_{*}(1x)(y)
\end{array}$$

\noindent Since $r_{*}(1x)\in\Mlt_{S}(Y)$ there are at most
finitely many $y$ that satisfy $r(x,y)\neq 0$. 
}

Analogously to the multiset monad the distribution monad
$\Dst\colon\Sets\rightarrow\Sets$ is defined as:
\begin{equation}
\label{DstUnitEqn}
\begin{array}{rcl}
\Dst(X)
& = &
\set{\varphi\colon X\rightarrow \unitR}{\support(\varphi)
   \mbox{ is finite and }\sum_{x\in X}\varphi(x) = 1}.
\end{array}
\end{equation}

\noindent Elements of $\Dst(X)$ are convex combinations
$s_{1}x_{1}+\cdots+s_{k}x_{k}$, where the probabilities
$s_{i}\in\unitR$ satisfy $\sum_{i}s_{i} = 1$.  Unit and multiplication
making $\Dst$ a monad can be defined as for $\Mlt_S$. The distribution
monad $\Dst$ is always symmetric monoidal (commutative) and its
category of algebras is the category \Cat{Conv} of convex sets with
affine maps between them, see
also~\cite{Keimel08,Doberkat06,Jacobs10e}.

The functor $\Dst\colon\Sets\rightarrow\Cat{Conv}$ also comes with
equality. We now choose $\Theta = 2 \in \Sets$, so that $\Omega =
\Dst(2) = [0,1]\in\Cat{Conv}$. Comparison $\cp\colon \Dst(Y)
\rightarrow \big(\Dst(X)\multimap[0,1]\big)$ can be defined as
in~\eqref{MltComparEqn} (but without conjugation).

\begin{proposition}
\label{DistProp}
The dagger category of tame relations $\TRel(\Cat{Conv},\cp)$ has
morphisms $X\rightarrow Y$ given by discrete `bistochastic' relations
$X\times Y\rightarrow [0,1]$ satisfying both:
$$\left\{\begin{array}{l}
\mbox{$X\rightarrow [0,1]^{Y}$ factors through
$X\rightarrow \Dst(Y)$} \\
\mbox{$Y\rightarrow [0,1]^{X}$ factors through
$Y\rightarrow \Dst(X)$.}
\end{array}\right.$$

\noindent We also write $\dBisRel = \TRel(\Cat{Conv},\cp)$ for
this category of discrete bistochastic relations.
\end{proposition}

\auxproof{ 
Assume $R\colon X\times Y\rightarrow [0,1]$, which corresponds to
$\widehat{R}\colon \Dst(X) \rightarrow (\Dst(Y)\multimap [0,1])$
in $\Cat{Conv}$, given by $\widehat{R}(\varphi)(\psi) = \sum_{x,y}
  \varphi(x)\cdot R(x,y) \cdot \psi(y)$. We shall prove the
  equivalence of:
\begin{enumerate}
\item[(a)] $\Lambda(R)\colon X\rightarrow [0,1]^{Y}$ factors through
$X\rightarrow \Dst(Y)$;

\item[(b)] there is a necessarily unique map $R_{*}\colon
\Dst(X)\rightarrow\Dst(Y)$ in \Cat{Conv} in the diagram:
$$\xymatrix@R1.5pc@C-1pc{
\Dst(X)\ar[rr]^-{\widehat{R}}\ar@{-->}[dr]_{R_{*}} & & 
   \Dst(Y)\multimap 2 \\
& \Dst(Y)\ar@{ >->}[ur]_{\cp}
}$$
\end{enumerate}

So assume~(a) holds. Then we can define $R_{*}(\varphi) \in
\Dst(Y)$ as finite sum of multisets, namely as $R_{*}(\varphi)
=\sum_{x}\varphi(x)\cdot\Lambda(R)(x)$. It makes the triangle in~(b)
commute:
$$\begin{array}{rcl}
\big(\cp \after R_{*}\big)(\varphi)(\psi) 
& = &
\cp(R_{*}(\varphi))(\psi) \\
& = &
\sum_{y}R_{*}(\varphi)(y)\cdot\psi(y) \\
& = &
\sum_{y}\Big(\sum_{x}\varphi(x)\cdot \Lambda(R)(x)(y)\Big)\cdot \psi(y) \\
& = &
\sum_{x,y}\varphi(x)\cdot R(x, y)\cdot \psi(y) \\
& = &
\widehat{R}(\varphi)(\psi).
\end{array}$$

Conversely, assume~(b) holds, so that we have a map $R_{*}\colon
\Dst(X)\rightarrow \Dst(Y)$ in $\Cat{Conv}$ in the above
triangle. Then:
$$\begin{array}{rcl}
R(x,y)
& = &
\widehat{R}(1x)(1y) \\
& = &
\cp(R_{*}(1x))(1y)  \\
& = &
R_{*}(1x)(y)
\end{array}$$

\noindent Since $R_{*}(1x)\in\Dst(Y)$ there are at most
finitely many $y$ that satisfy $R(x,y)\neq 0$. \QED
}

These bistochastic relations are reversible by
definition. Reversibility of arbitrary stochastic relations is studied
for instance in~\cite{Doberkat03}.

%% In some applications one may wish to replace the monads $\Mlt_S$ and
%% $\Dst$ that we have used in this subsection by non-finite versions
%% that allow countable support, so that `eventual' behaviour obtained
%% via limits can be incorporated.

\subsection{Hilbert spaces}\label{HilbertSubsec}

The so-called $\ell^2$-construction can be seen as an infinite
version of the multiset monad $\Mlt_{\mathbb{C}}$. For a set $X$ one
takes the square-summable sequences indexed by $X$, as in:
$$\begin{array}{rcl}
\ell^{2}(X)
& = &
\set{\varphi\colon X\rightarrow\mathbb{C}}{\sum_{x\in X}
   \|\varphi(x)\|^{2} < \infty},
\end{array}$$

\noindent where $\|\varphi(x)\|^{2} =
\varphi(x)\cdot\overline{\varphi(x)}$. As is well-known, the
$\ell^2$-construction forms a functor $\ell^{2} \colon \PInj
\rightarrow \Hilb$, but not a functor $\Sets \to \Hilb$, see
\textit{e.g.}~\cite{Barr92,HaghverdiS06}. However, in the present
setting we do not need functoriality for the indices of comparison
relations. Thus we have the (standard) inner products
$$\xymatrix@R1.5pc{
\overline{\ell^{2}(X)}\otimes\ell^{2}(X)
   \ar[r]^-{\inprod{-}{-}} & \mathbb{C}
\qquad\mbox{given by}\qquad
\inprod{\varphi}{\varphi'} = \sum_{x}\overline{\varphi(x)}\cdot\varphi'(x)
}$$

\noindent forming a symmetric cluster of comparison relations, much
like in~\eqref{MltComparEqn} for multisets. As we show below, it does
not matter if we consider these inner products as morphisms in
$\Vect_{\mathbb{C}}$ or in \Hilb. The resulting category of tame
relations has sets as objects and continuous linear functions
$\ell^{2}(X)\rightarrow\ell^{2}(Y)$ as morphisms $X\rightarrow
Y$. This follows from the lemma below.

Given an arbitrary Hilbert space $H$, we can consider its inner
product $\overline{H}\otimes H \xrightarrow{\inprod{-}{-}} \mathbb{C}$
as a comparison relation in the category $\Vect_{\mathbb{C}}$, as
already mentioned in Section~\ref{TameRelSec}.  It is well-known that
a linear map between Hilbert spaces is continuous if and only if it is
bounded.  Jorik Mandemaker suggested the next result (and proof),
which shows that boundedness/continuity can be captured in terms of
tameness.

\begin{lemma}
\label{BoundedLem}
Consider two Hilbert spaces $H_{1},H_{2}$, with their inner products
$\overline{H_i}\otimes H_{i} \xrightarrow{\inprod{-}{-}} \mathbb{C}$
as comparisons in $\Vect_{\mathbb{C}}$. There is a bijective
correspondence between:
$$\begin{prooftree}
{\xymatrix@R1.5pc{\mbox{ tame }\; \overline{H_{1}}\otimes H_{2}
   \ar[r]^-{r} & \mathbb{C}}}
\Justifies
{\xymatrix@R1.5pc{\mbox{ bounded } \; H_{1}\ar[r]_-{f} & H_{2}}}
\end{prooftree}$$

\noindent As a result, $\Hilb = \TRel(\Vect_{\mathbb{C}},
\inprod{-}{-})$, where $\inprod{-}{-}$ is the comparison cluster
$\big(\overline{H}\otimes H \xrightarrow{\inprod{-}{-}_H}
\mathbb{C}\big)_{H\in\Hilb}$ indexed by Hilbert spaces.
\end{lemma}

Thus, morphisms between Hilbert spaces can also be understood as
(tame) relations, like morphisms in many other categories of interest
in quantum foundations.

\begin{myproof}
If a linear map $f\colon H_{1}\rightarrow H_{2}$ is bounded, then it
has a dagger $f^{\dag} \colon H_{2}\rightarrow H_{1}$ satisfying
$\inprod{f(x)}{y} = \inprod{x}{f^{\dag}(y)}$, for all $x\in H_{1}$ and
$y\in H_{2}$. Thus, by construction, the relation $r(x,y) =
\inprod{f(x)}{y} = \inprod{x}{f^{\dag}(y)}$ is tame, with $r_{*} = f$
and $r^{*} = f^{\dag}$.

Conversely, given a tame relation $r\colon \overline{H_{1}}\otimes
H_{2} \rightarrow \mathbb{C}$ we use the Closed Graph Theorem in order
to show that $r_{*}\colon H_{1} \rightarrow H_{2}$ is
continuous. Assume we have a Cauchy sequence $(x_{n})_{n\in\NNO}$ in
$H_1$ with limit $x$, and let the sequence $(r_{*}(x_{n}))_{n\in\NNO}$
in $H_{2}$ have limit $z$. It suffices to show $r_{*}(x) = z$. We use
that the inner product is continuous (which follows from
Cauchy-Schwarz), in:
$$\begin{array}{rcl}
\inprod{r_{*}(x)}{y}
\hspace*{\arraycolsep} = \hspace*{\arraycolsep}
\inprod{r_{*}(\lim_{n}x_{n})}{y}
& = &
\inprod{\lim_{n}x_{n}}{r^{*}(y)} \\
& = &
\lim_{n}\inprod{x_{n}}{r^{*}(y)} \\
& = &
\lim_{n}\inprod{r_{*}(x_{n})}{y} \\
& = &
\inprod{\lim_{n}r_{*}(x_{n})}{y}
\hspace*{\arraycolsep} = \hspace*{\arraycolsep}
\inprod{z}{y}.
\end{array}$$

\noindent Since this holds for each $y$, we get $r_{*}(x) = y$ by the
mono-property of comparisons (or inner products). \QED
\end{myproof}

\subsection{Formal distributions}\label{FormDistrSubsec}

We now use the present framework of comparisons for re-describing the
dagger category of formal distributions introduced
in~\cite{BluteP11}. First we show how to capture polynomials via
multiset monads (from Subsection~\ref{MltDstSubsec}). Laurent
polynomials, with negative powers $x^{-1}$, are used
in~\cite{BluteP11}, but here we stick to ordinary polynomials.

As described in the previous subsection, a multiset
$\varphi\in\Mlt_{\NNO}(X)$ can be described as a formal sum
$n_{1}x_{1} + n_{2}x_{2} + \cdots + n_{k}x_{k}$, with
$n_{i}\in\NNO$. We might as well write $\varphi$ multiplicatively, as
in $x_{1}^{n_1}x_{2}^{n_2} \cdots x_{k}^{n_{k}}$. This is convenient,
because we can now describe a (multivariate) polynomial as a
`multiset of multisets' $p\in\Mlt_{S}(\Mlt_{\NNO}(X))$. If we use
additive notation for the outer multiset $\Mlt_{S}$ we can write $p$
as formal sum:
$$\textstyle p = \sum_{i}s_{i}\varphi_{i}
\qquad\mbox{where}\qquad
\varphi_{i} = x_{i1}^{n_{i1}}\cdots x_{ik_{i}}^{n_{ik_{i}}}\in \Mlt_{\NNO}(X).$$

\noindent The univariate polynomials, with only one variable, appear
by taking $X=1$, namely as $p\in\Mlt_{S}(\Mlt_{\NNO}(1)) =
\Mlt_{S}(\NNO)$.  Such a $p$ can be written as $\sum_{i}s_{i}n_{i}$,
or as polynomial $\sum_{i}s_{i}x^{n_i}$ for some variable $x$.

We write $S[X] = \Mlt_{S}\Mlt_{\NNO}(X)$ for the set of (multivariate)
polynomials with variables from an arbitrary set $X$ and coefficients
from the commutative semiring $S$. These polynomials are finite, by
construction. Possibly infinite polynomials---also known as power
series or as formal distributions---are obtained via the function
space $S[[X]] = S^{\Mlt_{\NNO}(X)}$. Thus $q\in S[[X]]$ can be written
as (possibly infinite) formal sum
$\sum_{\varphi\in\Mlt_{\NNO}(X)}q(\varphi)\varphi$, where
$q(\varphi)\in S$ gives the coefficient. There is an obvious inclusion
$S[X] \hookrightarrow S[[X]]$ that will play the role of comparison
below. But first we need to relate finite and infinite polynomials
more closely.

\begin{lemma}
\label{FinInfPolLem}
For a commutative semiring $S$ and set $X$ of `variables', modules
$S[X],S[[X]]\in\Cat{Mod}_{S}$ of finite and infinite multivariate
polynomials, with variables from $X$, are defined as:
$$S[X] = \Mlt_{S}\Mlt_{\NNO}(X)
\qquad\mbox{and}\qquad
S[[X]] = S^{\Mlt_{\NNO}(X)} = \Hom_{\Sets}\big(\Mlt_{\NNO}(X), S\big).$$

\noindent Then: $S[[X]] \cong \big(S[X]\multimap S\big)$, where
$\multimap$ is exponent in $\Cat{Mod}_{S}$.
\end{lemma}

\begin{myproof}
We use the following chain of isomorphisms, which exploits that
multiset $\Mlt_{S}\colon \Sets \rightarrow\Mod_{S}$ is the free
functor and that the category $\Mod_{S}$ of modules is monoidal
closed.
$$\begin{array}[b]{rcl}
S[[X]] 
\hspace*{\arraycolsep} = \hspace*{\arraycolsep}
S^{\Mlt_{\NNO}(X)} 
& = &
\Hom_{\Sets}\big(\Mlt_{\NNO}(X), \; S\big) \\
& \cong &
\Hom_{\Mod_{S}}\big(\Mlt_{S}\Mlt_{\NNO}(X), \; S\big) \\
& \cong &
\Hom_{\Mod_{S}}\big(S[X] \otimes \Mlt_{S}(1), \; S\big) \\
& & \qquad \mbox{since $\Mlt_{S}(1) \cong S$ is the tensor unit} \\
& \cong &
\Hom_{\Mod_{S}}\big(\Mlt_{S}(1), \; S[X] \multimap S\big) \\
& \cong &
\Hom_{\Sets}\big(1, \; S[X] \multimap S\big) \\
& \cong &
\big(S[X] \multimap S\big).
\end{array}\eqno{\QEDbox}$$
\end{myproof}

We shall introduce comparisons $\cp \colon S[X]\otimes S[X]\rightarrow
S$ in the category $\Cat{Mod}_{S}$ of $S$-modules, following the
recipe from Subsection~\ref{MonadSubsec}. We start from the equality
relation $\eq\colon \Mlt_{\NNO}(X)\times\Mlt_{\NNO}(X)\rightarrow S$,
following~\eqref{EqFunEqn}, which gives rise to $\cp$ as composite:
\begin{equation}
\label{FormDistrComparEqn}
\vcenter{\xymatrix@R-1.4pc@C-1.6pc{
S[X]\otimes S[X]\ar@{=}[r] &
   \Mlt_{S}\big(\Mlt_{\NNO}(X)\big) \otimes \Mlt_{S}\big(\Mlt_{\NNO}(X)\big)
      \ar@{=}[d]_-{\wr} \\
& \Mlt_{S}\Big(\Mlt_{\NNO}(X) \times \Mlt_{\NNO}(X)\Big)
   \ar[rrr]^-{\Mlt_{S}(\eq)} & & & \Mlt_{S}(S)\ar[rr]^-{\mu} & & S
}}
\end{equation}

\noindent Concretely, $\cp(p,p') = \sum_{\varphi\in\Mlt_{\NNO}(X)}p(\varphi)
\cdot p'(\varphi)$, like in~\eqref{MltComparEqn}.

More generally, relations in this setting will be module maps of the
form $S[X]\otimes S[Y]\rightarrow S$. Obviously, by Curry-ing they can
also be described as maps $S[Y] \rightarrow \big(S[X]\multimap S\big)
\cong S[[X]]$. It is not hard to see that the map $S[X] \rightarrow
S[[X]]$ corresponding to comparison $\cp$
in~\eqref{FormDistrComparEqn} is inclusion. In particular, this shows
that the mono-requirement from Definition~\ref{ComparRelDef} is
satisfied.

\auxproof{
We note that the map $\big(S[X] \multimap S) \conglongrightarrow S[[X]]$
is given by 
$$f\mapsto \lamin{\varphi}{\Mlt_{\NNO}(X)}{f(1\varphi)}.$$

\noindent Thus, starting from $p\in S[X]$ we apply this to
$\cp(p,-)\in (S[X]\multimap S)$ and get in $S[[X]]$,
$$\begin{array}{rcl}
\lefteqn{\lamin{\varphi}{\Mlt_{\NNO}(X)}{\cp(p, 1\varphi)}} \\
& = &
\lamin{\varphi}{\Mlt_{\NNO}(X)}{p(\varphi)} \\
& = &
p.
\end{array}$$
}

There is one further observation that we need to make.

\begin{lemma}
\label{AdditiveLem}
Each multiset monad $\Mlt_{S}$ is an `additive'
monad~\cite{CoumansJ12}: it maps finite coproducts to products, in a
canonical way: $\Mlt_{S}(0) \cong 1$ and $\Mlt_{S}(X+Y) \cong
\Mlt_{S}(X)\times \Mlt_{S}(Y)$. The latter isomorphism will be written
explicitly as:
$${\def\labelstyle{\textstyle}\xymatrix@C+2pc{
\Mlt_{S}(X+Y)\ar@/^2ex/[rr]^-{\chi \longmapsto \tuple{\chi(\kappa_{1}-),
   \chi(\kappa_{2}-)}}
& \cong &
\Mlt_{S}(X)\times \rlap{$\Mlt_{S}(Y)$}\ar@/^2ex/[ll]^-{\varphi\star\psi\longmapsfrom(\varphi,\psi)}
}}$$

\noindent where the operation $\star$ multiplies $\varphi,\psi$, after
appropriate relabeling has put them in the same set of multisets:
$$\begin{array}{rcl}
\varphi\star\psi
& = &
\Mlt_{S}(\kappa_{1})(\varphi)\cdot\Mlt_{S}(\kappa_{2})(\psi).
\end{array}$$

\noindent (We use multiplicative notation $\cdot$ in the definition of
$\star$ for multiset addition because later on we use $\star$ when we
read multisets multiplicatively; the $\kappa_i$ are the coprojections
associated with the coproduct.) \QED
\end{lemma}

Using this additivity of the multiset monad we show that relations can
be described in another way as formal distributions.

\begin{proposition}
\label{FormDistrRelProp}
In the setting described above, there is an isomorphism of modules
between formal distributions in the coproduct $X+Y$ and relations on
$X$ and $Y$, as in:
$$\begin{array}{rcl}
S[[X+Y]]
& \cong &
\big(S[X]\otimes S[Y] \multimap S\big).
\end{array}$$

\noindent The formal distribution in $S[[X+X]]$ corresponding to the
comparison relation $\cp\colon S[X]\otimes S[X]\rightarrow S$ is
the function:
$$\lamin{\varphi}{\Mlt_{\NNO}(X+X)}{\left\{\begin{array}{ll}
   1 \; & \mbox{if }\varphi(\kappa_{1}-) = \varphi(\kappa_{2}-) 
      \mbox{ in }\Mlt_{S}(X) \\
   0 & \mbox{otherwise.} \end{array}\right.}.$$
\end{proposition}

\begin{myproof}
Because multiset monads are additive and free functors we have:
$$\begin{array}{rcl}
S[X+Y]
& = &
\Mlt_{S}\big(\Mlt_{\NNO}(X+Y)\big) \\
& \cong &
\Mlt_{S}\big(\Mlt_{\NNO}(X) \times \Mlt_{\NNO}(Y)\big) \\
& \cong &
\Mlt_{S}\big(\Mlt_{\NNO}(X)\big) \otimes \Mlt_{S}\big(\Mlt_{\NNO}(Y)\big) \\
& = &
S[X] \otimes S[Y].
\end{array}$$

\noindent Hence Lemma~\ref{FinInfPolLem} gives:
$$\begin{array}{rcl}
S[[X+Y]] 
& \cong & 
\big(S[X+Y] \multimap S\big) \\
& \cong &
\big(S[X] \otimes S[Y] \multimap S\big) \\
& \cong &
\big(S[X] \multimap (S[Y] \multimap S)\big).
\end{array}$$

\noindent The formal power series in $S[[X+X]]$ can be obtained by
following these isomorphisms backwards. \QED
\end{myproof}

In~\cite{BluteP11} a category of formal distributions is defined with
(finite) sets as objects and morphisms $X\rightarrow Y$ given by
``tame'' formal distributions $p\in S[[X+Y]]$. Here we re-describe
them in the current framework, namely as category
$\TRel(\Cat{Mod}_{S},\cp)$ for the comparison
cluster~\eqref{FormDistrComparEqn}. Indeed, for a morphism $r\colon
X\rightarrow Y$ in this category, considered as a map of modules
$r\colon S[Y]\rightarrow \big(S[X]\multimap S\big)$, tameness means
the existence of a map $r^{*}\colon S[Y]\rightarrow S[X]$ in
$\Mod_{S}$, as indicated:
$$\xymatrix@R1.5pc@C-1pc{
S[Y]\ar[rr]^-{r}\ar@{-->}[dr]_{r^{*}} & & 
   \big(S[X]\multimap S\big) \rlap{$\;=S[[X]]$} \\
& S[X]\ar@{^(->}[ur]_{\cp}
}$$

\noindent (And similarly for $r_*$.) We can translate this condition
to formal distributions as morphisms, using
Lemma~\ref{FormDistrRelProp}. Indeed, a formal distribution $p\in
S[[X+Y]]$ gives rise to a map $\widehat{p} \colon S[Y] \rightarrow
S[[X]]$, namely:
$$\begin{array}{rcl}
\widehat{p}(q)
& = &
\lamin{\varphi}{\Mlt_{\NNO}(X)}{\displaystyle\sum_{\psi\in\Mlt_{\NNO}(Y)}
   p\big(\varphi\star\psi\big)}
\end{array}$$

\noindent where $\star$ is the operation for additivity from
Lemma~\ref{AdditiveLem}. Tameness says that $\widehat{p}(q)$ is a
finite polynomial, for each $q\in S[Y]$; it means that $\widehat{p}$
factors as $S[Y]\rightarrow S[X]$ (and vice-versa).

For completeness we include formulations of composition and dagger 
for formal distributions. Given $p\in S[[X+Y]]$ and $q\in
S[[Y+Z]]$ we have:
$$\begin{array}{rcl}
p \relafter q
& = &
\lamin{\chi}{\Mlt_{\NNO}(X+Z)}{\displaystyle\sum_{\psi\in\Mlt_{\NNO}(Y)}
   p\Big(\chi(\kappa_{1}-)\star \psi\Big) \cdot 
   q\Big(\psi \star \chi(\kappa_{2}-)\Big)} \\
\big(p\big)^{\dag}
& = &
\lamin{\chi}{\Mlt_{\NNO}(Y+X)}{p\Big(\chi(\kappa_{2}-)\star
   \chi(\kappa_{1}-)\Big)}.
\end{array}$$

\noindent The tameness requirement ensures that these sums $\sum$
exist. It is not hard to see that the formal distribution described at
the end of Proposition~\ref{FormDistrRelProp} is the identity map.

In the end we see that this formal distribution example fits in the
general recipe for monads $T$ from Subsection~\ref{MonadSubsec},
except that we start with an (additional) additive monad $R$. Equality
is used on $R$, in the form of maps $\eq\colon R(X)\times
R(X)\rightarrow T(\Theta)=\Omega$, and is lifted to comparisons $\cp
\colon TR(X)\otimes TR(X)\rightarrow \Omega$. Additivity of $R$ allows
us to translate between coproducts and products to make the machinery
work (via the $\star$'s above). Hence one may construct other examples
of dagger categories of this kind.

\auxproof{
$$\begin{array}{rcl}
(\idmap \relafter p_{1})(\chi)
& = &
\sum_{\psi}p_{1}\big(\chi(\kappa_{1}-)\star\psi\big) \cdot
   \idmap\big(\psi\star\chi(\kappa_{2}-)\big) \\
& = &
\sum\set{p_{1}\big(\chi(\kappa_{1}-)\star\psi\big)}{
   \big(\psi\star\chi(\kappa_{2}-)\big)(\kappa_{1}-) =
   \big(\psi\star\chi(\kappa_{2}-)\big)(\kappa_{2}-)} \\
& = &
\sum\set{p_{1}\big(\chi(\kappa_{1}-)\star\psi\big)}{
   \psi = \chi(\kappa_{2}-)} \\
& = &
p_{1}\big(\chi(\kappa_{1}-)\star\chi(\kappa_{2}-)\big) \\
& = &
p_{1}(\chi) \\
(p_{2} \relafter \idmap)(\chi)
& = &
\sum_{\psi}\idmap\big(\chi(\kappa_{1}-)\star\psi\big) \cdot
   p_{2}\big(\psi\star\chi(\kappa_{2}-)\big) \\
& = &
p_{2}\big(\chi(\kappa_{1}-)\star\chi(\kappa_{2}-)\big) \\
& = &
p_{2}(\chi).
\end{array}$$

\noindent Associativity is more work. Assume further $p_{3}\in S[[Z+W]]$.
$$\begin{array}{rcl}
\lefteqn{\big(p_{3} \relafter (p_{2}\relafter p_{1})\big)(\chi)} \\
& = &
\sum_{\psi}(p_{2}\relafter p_{1})\big(\chi(\kappa_{1}-)\star\psi\big) \cdot
   p_{2}\big(\psi\star\chi(\kappa_{2}-)\big) \\
& = &
\sum_{\psi}\Big(\sum_{\varphi}
  p_{1}\big((\chi(\kappa_{1}-)\star\psi)(\kappa_{1}-)\star \varphi\big) 
  \cdot 
  p_{2}\big(\varphi\star (\chi(\kappa_{1}-)\star\psi)(\kappa_{2}-)\big)\Big)
  \cdot \\
& & \qquad
   p_{2}\big(\psi\star\chi(\kappa_{2}-)\big) \\
& = &
\sum_{\varphi,\psi} p_{1}\big(\chi(\kappa_{1}-)\star\varphi\big) \cdot
   p_{2}\big(\varphi\star\psi\big) \cdot
   p_{3}\big(\psi \star \chi(\kappa_{2}-)\big) \\
& = &
\sum_{\varphi}p_{1}\big(\chi(\kappa_{1}-)\star\varphi\big) \cdot \\
& & \qquad
\Big(\sum_{\psi}p_{2}\big((\varphi\star\chi(\kappa_{2}-))(\kappa_{1}-)\star\psi\big) \cdot
   p_{3}\big(\psi\star(\varphi\star\chi(\kappa_{2}-))(\kappa_{2}-)\big)\Big) \\
& = &
\sum_{\varphi}p_{1}\big(\chi(\kappa_{1}-)\star\varphi\big) \cdot
   (p_{3} \relafter p_{2})\big(\varphi\star\chi(\kappa_{2}-)\big) \\
& = &
\big((p_{3} \relafter p_{2})\relafter p_{1}\big)(\chi).
\end{array}$$

Finally, we check the dagger properties.
$$\begin{array}{rcl}
\big(\idmap\big)^{\dag}(\chi)
& = &
\idmap\big(\chi(\kappa_{2}-)\star\chi(\kappa_{1}-)\big) \\
& = &
\left\{\begin{array}{ll}
   1 \; & \mbox{if }
   \big(\chi(\kappa_{2}-)\star\chi(\kappa_{1}-)\big)(\kappa_{1}-) =
   \big(\chi(\kappa_{2}-)\star\chi(\kappa_{1}-)\big)(\kappa_{2}-) \\
   0 & \mbox{otherwise} \end{array}\right. \\
& = &
\left\{\begin{array}{ll}
   1 \; & \mbox{if }
   \big(\chi(\kappa_{2}-) = \chi(\kappa_{1}-) \\
   0 & \mbox{otherwise} \end{array}\right. \\
& = &
\left\{\begin{array}{ll}
   1 \; & \mbox{if }
   \big(\chi(\kappa_{1}-) = \chi(\kappa_{2}-) \\
   0 & \mbox{otherwise} \end{array}\right. \\
& = &
\idmap(\chi) \\
\big(p_{2} \relafter p_{1}\big)^{\dag}(\chi)
& = &
\big(p_{2} \relafter p_{1}\big)
   \big(\chi(\kappa_{2}-)\star\chi(\kappa_{1}-)\big) \\
& = &
\sum_{\psi}p_{1}\big(\big(\chi(\kappa_{2}-)\star\chi(\kappa_{1}-)\big)(\kappa_{1}-)\star\psi\big) \cdot
   p_{2}\big(\psi\star\big(\chi(\kappa_{2}-)\star\chi(\kappa_{1}-)\big)(\kappa_{2}-)\big) \\
& = &
\sum_{\psi}p_{1}\big(\chi(\kappa_{2}-)\star\psi\big) \cdot
   p_{2}\big(\psi\star\chi(\kappa_{1}-)\big) \\
& = &
\sum_{\psi}\big(p_{2}\big)^{\dag}\big(\chi(\kappa_{1}-)\star\psi\big) 
   \cdot \big(p_{1}\big)^{\dag}\big(\psi\star\chi(\kappa_{2}-)\big) \\
& = &
\Big(\big(p_{1}\big)^{\dag} \relafter \big(p_{2}\big)^{\dag}\Big)(\chi) \\
\big(p_{1}\big)^{\dag\dag}(\chi)
& = &
\big(p_{1})^{\dag}\big(\chi(\kappa_{2}-)\star\chi(\kappa_{1}-)\big) \\
& = &
p_{1}\big(\big(\chi(\kappa_{2}-)\star\chi(\kappa_{1}-)\big)(\kappa_{2}-)
   \star\big(\chi(\kappa_{2}-)\star\chi(\kappa_{1}-)\big)(\kappa_{1}-)\big) \\
& = &
p_{1}\big(\chi(\kappa_{1}-)\star\chi(\kappa_{2}-)\big) \\
& = &
p_{1}(\chi).
\end{array}$$

}

%% \begin{remark}
%% In the theory of formal distributions an important role is played by
%% the so-called residue operation $\res$. In its simplest form it takes
%% of a univariant polynomial $p = \sum_{n\in\NNO}s_{n}x^{n}$ the
%% coefficient at $0$\footnote{In~\cite{BluteP11} the coefficient is
%%   taken at -1, for Laurent polynomials, but that is irrelevant here.},
%% so $\res(p) = s_{0}$. This can be generalised to multivariate
%% polynomials. For instance
%% $\res_{x}\big(\sum_{n,m\in\NNO}s_{n,m}x^{n}y^{m}\big) =
%% \sum_{m}s_{0,m}y^{m}$.

%% We have avoided this residue operation via the additivity of the
%% multiset monad. For instance, the unique map
%% $\smash{0\stackrel{!}{\rightarrow}1}$ gives rise to:
%% $$\xymatrix@C+.5pc{
%% \Mlt_{\NNO}(0) = 1\ar[r]^-{0=\Mlt_{\NNO}(!)} & \NNO = \Mlt_{\NNO}(1)
%% \quad\mbox{and}\quad
%% S^{\mathbb{N}} = S[[1]] \ar[r]^-{\res = S^{0}} &S[[0]] = S^{1}
%% }$$

%% \noindent Similarly, the `right coprojection' $\kappa_{2}\colon
%% 1\rightarrow 2 = 1+1$ gives rise to:
%% $$\xymatrix@C+1.5pc{
%% \Mlt_{\NNO}(1) = \NNO\ar[r]^-{\tuple{0,-}=\Mlt_{\NNO}(\kappa_{2})} & 
%%    \NNO^{2} = \Mlt_{\NNO}(2)
%% \quad\mbox{and}\quad
%% S^{\mathbb{N}^{2}} = S[[2]] \ar[r]^-{\res_{x} = S^{\tuple{0,-}}} &
%%    S[[1]] = S^{\NNO}
%% }$$

%% \noindent This map on the left can be understood as $y^{m}\mapsto
%% x^{0}y^{m}$.
%% \end{remark}

\section{The category of bifinite multirelations}\label{BfMRelSec}

Subsection~\ref{MltDstSubsec} introduced the category $\BifMRel_{S} =
\TRel(\Mod_{S},\cpMlt)$ of sets and bifinite multirelations, with
values in an involutive semiring $S$ (such as $\mathbb{C}$). Here we
shall investigate its categorical structure in more detail.

(Describing the categorical structure of categories
$\TRel(\cat{A},\cp)$ in full generality turns out to be rather
involved. In contrast, for several examples, this structure is
essentially straightforward. That is why we prefer this more concrete
approach.)

There is a special reason why we concentrate on $\BifMRel_{S}$---and
not on other categories of tame relatons.  The category $\BifMRel_{S}$
may be seen a universe for `discrete' quantum computation (like
in~\cite{Jacobs11a}), just like the category of Hilbert spaces may be
used for `continuous' computation. We shall illustrate this in a
moment, but first we describe the category $\BifMRel_{S}$ concretely,
and state an elementary result.

Objects in the category $\BifMRel_{S}$ are sets; it is important that
infinite sets are allowed as objects, so that computations with
infinitely many (orthogonal) states can be covered---unlike in
finite-dimensional vector (or Hilbert) spaces. A morphism $r\colon
X\rightarrow Y$ in $\BifMRel_{S}$ is a multirelation $r\colon X\times
Y\rightarrow S$ such that for each $x\in X$ the subset
$\set{y}{r(x,y)\neq 0} = \support(r(x,-))$ is finite, and similarly,
for each $y\in Y$ the set $\set{x}{r(x,y)\neq 0} = \support(r(-,y))$
is finite.  Composition of $r\colon X\rightarrow Y$ with $s\colon
Y\rightarrow Z$ can be described as matrix compositon: $(s\relafter
r)(x,z) = \sum_{y} r(x,y)\cdot s(y,z)$. The dagger $r^{\dag}\colon
Y\rightarrow X$ is given by the adjoint matrix: $r^{\dag}(y,x) =
\overline{r(x,y)}$, obtained by mirroring and conjugation in $S$.
Notice that the special case $S=2=\{0,1\}$ covers the category
$\BifRel = \BifMRel_{2}$ of bifinite relations.

We show how unitary maps give rise to bistochastic relations (for the
standard semiring examples in this context).

\begin{lemma}
\label{BifMRelUnitaryLem}
Assume an involutive semiring $S$ like $2, \mathbb{R},
\mathbb{R}_{\geq 0}$ or $\mathbb{C}$, for which the mapping $a\mapsto
a\cdot\overline{a}$ yields a function $S\rightarrow \mathbb{R}_{\geq
  0}$, which we write as squared norm $\|-\|^{2}$. A unitary map
$r\colon X\rightarrow Y$ in $\BifMRel_{S}$ then yields a discrete
bistochastic relation, $\|r\|^{2}\colon X\rightarrow Y$,
\textit{i.e.}~a morphism in the category $\dBisRel$ from
Proposition~\ref{DistProp}, given by $\|r\|^{2}(x,y) =
\|r(x,y)\|^{2}$.
\end{lemma}

\begin{myproof}
Suppose $r\colon X\rightarrow Y$ in \BifMRel is unitary,
\textit{i.e.}~$r^{\dag}$ is $r$'s inverse. Then, for each $x\in X$,
$$\textstyle 1 
= 
\idmap[X](x,x) 
= 
(r^{\dag} \relafter r)(x,x)
=
\sum_{y}r(x,y)\cdot \overline{r(y,x)}
=
\sum_{y}\|r(x,y)\|^{2}.$$

\noindent And similarly for $y\in Y$. Hence, post-composition with the
squared norm $\|-\|^{2}\colon S \rightarrow \mathbb{R}_{\geq 0}$ turns
the unitary bifinite multirelation $r\colon X\times Y\rightarrow S$
into a bistochastic relation $X\times Y\rightarrow [0,1]$. \QED

\auxproof{
If we use the \textit{ad hoc} notation $D(r) = \|-\|^{2} \after r
\colon X\times Y\rightarrow \mathbb{R}_{\geq 0}$, then we
can see where functoriality is problematic:
$$\begin{array}{rcl}
\big(D(s) \relafter D(s)\big)(x,z)
& = &
\sum_{y}D(r)(x,y)\cdot D(s)(y,z) \\
& = &
\sum_{y}\|r(x,y\|^{2}\cdot \|s(y,z)\|^{2} \\
& = &
\sum_{y} r(x,y)\cdot\overline{r(x,y)}\cdot s(y,z) \cdot
   \overline{s(y,z)} \\
& = &
\ldots \\
& = &
\|\sum_{y}r(x,y)\cdot s(y,z)\|^{2} \\
& = &
\|(s\relafter r)(x,z)\|^{2} \\
& = &
D(s\relafter r)(x,z)
\end{array}$$
}
\end{myproof}

Notice that an arbitrary morphism $q\colon 1\rightarrow 2$ corresponds
to a map $q\colon 1\times 2 \rightarrow S$, and thus to two scalars $a
= q(*,0)\in S$ and $b = q(*,1)\in S$, where we use $1 = \{*\}$ and $2
= \{0,1\}$. One can call such a $q$ a \emph{unit} if $\|q\|^{2} = 1$,
\textit{i.e.}~if $(q^{\dag}\relafter q)(*,*) = \|a\|^{2} + \|b\|^{2} =
1$ in $\mathbb{R}_{\geq 0}$. Such a unit is a \emph{quantum bit} for
$S=\mathbb{C}$ and a \emph{classical bit} for $S=2$.

We briefly illustrate the use of the category $\BifMRel_{\mathbb{C}}$
to model discrete quantum computations (on an infinite state space).
In~\cite{Jacobs11a} quantum walks (see
also~\cite{Kempe03,VenegasAndraca08}) are investigated in relation to
possibilistic and probabilistic walks.  Such walks involves discrete
steps on an infinite line, given by the integers $\mathbb{Z}$. In a
single move, left or right steps can be made, described as $-1$
decrements or $+1$ increments. The walks are steered by Hadamard's
matrix acting on a qubit. They can be described via a function
$\mathbb{C}^{2}\otimes\Mlt_{\mathbb{C}}(\mathbb{Z}) \rightarrow
\mathbb{C}^{2}\otimes\Mlt_{\mathbb{C}}(\mathbb{Z})$, where
$\mathbb{C}^2$ represents the qubit,
see~\cite{Jacobs11a}. Alternatively, they can be described via a
bifinite multirelation on $\mathbb{Z}+\mathbb{Z}$. We write $\kappa_1$
and $\kappa_2$ as left and right coprojection for this coproduct,
corresponding to the up and down orientations of the qubit that steers
the movement. This kind of quantum walk can now be given as an endomap
$q\colon \mathbb{Z}+\mathbb{Z} \rightarrow \mathbb{Z}+\mathbb{Z}$ in
$\BifMRel_{\mathbb{C}}$, which we describe by listing only the
non-zero values of $q$, as multirelation:
$$\big(\mathbb{Z}+\mathbb{Z}\big)\times\big(\mathbb{Z}+\mathbb{Z}\big)
\stackrel{q}{\longrightarrow} \mathbb{C}
\qquad\mbox{where}\qquad
\left\{
\begin{array}{rcl}
q(\kappa_{1}n,\kappa_{1}(n-1)) & = & \frac{1}{\sqrt{2}} \\
q(\kappa_{1}n,\kappa_{2}(n+1)) & = & \frac{1}{\sqrt{2}} \\
q(\kappa_{2}n,\kappa_{1}(n-1)) & = & \frac{1}{\sqrt{2}} \\
q(\kappa_{2}n,\kappa_{2}(n+1)) & = & -\frac{1}{\sqrt{2}}
\end{array}\right.
%% (\kappa_{i}n,\kappa_{j}m) \longmapsto \left\{\begin{array}{ll}
%% \frac{1}{\sqrt{2}}  & \mbox{if }i=0,j=0\mbox{ and }m=n-1 \\
%% \frac{1}{\sqrt{2}}  & \mbox{if }i=0,j=1\mbox{ and }m=n+1 \\
%% \frac{1}{\sqrt{2}}  & \mbox{if }i=1,j=0\mbox{ and }m=n+1 \\
%% -\frac{1}{\sqrt{2}} \; & \mbox{if }i=1,j=1\mbox{ and }m=n-1 \\
%% 0 & \mbox{otherwise.}
%%\end{array}\right.
$$

\noindent The $n\in\mathbb{Z}$ in the first argument of $q$ represents
the current position; the second argument describes the successor
position, which is either a step left or right. The labels $\kappa_i$
capture orientations.

It is not hard to see that this map $q$ is unitary.  By iterating the map
in \BifMRel, like in $q^{2} = q \relafter q, q^{3} = q \relafter q
\relafter q, \ldots$, and subsequently taking the resulting
bistochastic relation (see Lemma~\ref{BifMRelUnitaryLem}), one can
compute the iterated distributions of the original quantum walk (and
the stationary distribution as suitable limit).

\auxproof{
Notice that we can write the non-zero case of $q^{\dag}$ as:
$$\begin{array}{rclcrcl}
q^{\dag}(\kappa_{1}(n-1),\kappa_{1}n) & = & \frac{1}{\sqrt{2}}
& \qquad &
q^{\dag}(\kappa_{2}(n+1),\kappa_{1}n) & = & \frac{1}{\sqrt{2}} \\
q^{\dag}(\kappa_{1}(n-1),\kappa_{2}n) & = & \frac{1}{\sqrt{2}}
& &
q^{\dag}(\kappa_{2}(n+1),\kappa_{2}n) & = & -\frac{1}{\sqrt{2}}
\end{array}$$

\noindent \textit{i.e.}~as:
$$\begin{array}{rclcrcl}
q^{\dag}(\kappa_{1}m,\kappa_{1}(m+1)) & = & \frac{1}{\sqrt{2}}
& \qquad &
q^{\dag}(\kappa_{2}m,\kappa_{1}(m-1)) & = & \frac{1}{\sqrt{2}} \\
q^{\dag}(\kappa_{1}m,\kappa_{2}(m+1)) & = & \frac{1}{\sqrt{2}}
& &
q^{\dag}(\kappa_{2}m,\kappa_{2}(m-1)) & = & -\frac{1}{\sqrt{2}}
\end{array}$$

\noindent Then:
$$\begin{array}{rcl}
\lefteqn{\big(q \relafter q^{\dag}\big)(\kappa_{1}m, \kappa_{1}n)} \\
& = &
\sum_{x}q^{\dag}(\kappa_{1}m,x)\cdot q(x,\kappa_{1}n) \\
& = &
q^{\dag}(\kappa_{1}m,\kappa_{1}(m+1))\cdot q(\kappa_{1}(m+1),\kappa_{1}n) +
q^{\dag}(\kappa_{1}m,\kappa_{2}(m+1))\cdot q(\kappa_{2}(m+1),\kappa_{1}n) \\
& = &
\left\{\begin{array}{ll}
\frac{1}{\sqrt{2}}\cdot \frac{1}{\sqrt{2}} +
   \frac{1}{\sqrt{2}} \cdot \frac{1}{\sqrt{2}} \quad & \mbox{if }n=m \\
0 & \mbox{otherwise} \end{array}\right. \\
& = &
\idmap(\kappa_{1}m,\kappa_{1}n) \\
\lefteqn{\big(q \relafter q^{\dag}\big)(\kappa_{1}m, \kappa_{2}n)} \\
& = &
\sum_{x}q^{\dag}(\kappa_{1}m,x)\cdot q(x,\kappa_{2}n) \\
& = &
q^{\dag}(\kappa_{1}m,\kappa_{1}(m+1))\cdot q(\kappa_{1}(m+1),\kappa_{2}n) +
q^{\dag}(\kappa_{1}m,\kappa_{2}(m+1))\cdot q(\kappa_{2}(m+1),\kappa_{2}n) \\
& = &
\left\{\begin{array}{ll}
\frac{1}{\sqrt{2}}\cdot \frac{1}{\sqrt{2}} +
   \frac{1}{\sqrt{2}} \cdot -\frac{1}{\sqrt{2}} \quad & \mbox{if }n=m+2 \\
0 & \mbox{otherwise} \end{array}\right. \\
& = &
0 \\
& = &
\idmap(\kappa_{1}m,\kappa_{2}n) \\
\lefteqn{\big(q \relafter q^{\dag}\big)(\kappa_{2}m, \kappa_{1}n)} \\
& = &
\sum_{x}q^{\dag}(\kappa_{2}m,x)\cdot q(x,\kappa_{1}n) \\
& = &
q^{\dag}(\kappa_{2}m,\kappa_{1}(m-1))\cdot q(\kappa_{1}(m-1),\kappa_{1}n) +
q^{\dag}(\kappa_{2}m,\kappa_{2}(m-1))\cdot q(\kappa_{2}(m-1),\kappa_{1}n) \\
& = &
\left\{\begin{array}{ll}
\frac{1}{\sqrt{2}}\cdot \frac{1}{\sqrt{2}} +
   -\frac{1}{\sqrt{2}} \cdot \frac{1}{\sqrt{2}} \quad & \mbox{if }n=m-2 \\
0 & \mbox{otherwise} \end{array}\right. \\
& = &
0 \\
& = &
\idmap(\kappa_{2}m,\kappa_{1}n) \\
\lefteqn{\big(q \relafter q^{\dag}\big)(\kappa_{2}m, \kappa_{2}n)} \\
& = &
\sum_{x}q^{\dag}(\kappa_{2}m,x)\cdot q(x,\kappa_{1}n) \\
& = &
q^{\dag}(\kappa_{2}m,\kappa_{1}(m-1))\cdot q(\kappa_{1}(m-1),\kappa_{2}n) +
q^{\dag}(\kappa_{2}m,\kappa_{2}(m-1))\cdot q(\kappa_{2}(m-1),\kappa_{2}n) \\
& = &
\left\{\begin{array}{ll}
\frac{1}{\sqrt{2}}\cdot \frac{1}{\sqrt{2}} +
   \frac{1}{\sqrt{2}} \cdot \frac{1}{\sqrt{2}} \quad & \mbox{if }n=m \\
0 & \mbox{otherwise} \end{array}\right. \\
& = &
\idmap(\kappa_{2}m,\kappa_{1}n).
\end{array}$$

\noindent Similarly, the other way around:
$$\begin{array}{rcl}
\lefteqn{\big(q^{\dag} \relafter q\big)(\kappa_{1}n, \kappa_{1}m)} \\
& = &
\sum_{x}q(\kappa_{1}n,x)\cdot q^{\dag}(x,\kappa_{1}m) \\
& = &
q(\kappa_{1}n,\kappa_{1}(n-1))\cdot q^{\dag}(\kappa_{1}(n-1),\kappa_{1}m) +
q(\kappa_{1}n,\kappa_{2}(n+1))\cdot q^{\dag}(\kappa_{2}(n+1),\kappa_{1}m) \\
& = &
\left\{\begin{array}{ll}
\frac{1}{\sqrt{2}}\cdot \frac{1}{\sqrt{2}} +
   \frac{1}{\sqrt{2}} \cdot \frac{1}{\sqrt{2}} \quad & \mbox{if }n=m \\
0 & \mbox{otherwise} \end{array}\right. \\
& = &
\idmap(\kappa_{1}n,\kappa_{1}m) \\
\lefteqn{\big(q^{\dag} \relafter q\big)(\kappa_{1}n, \kappa_{2}m)} \\
& = &
\sum_{x}q(\kappa_{1}n,x)\cdot q^{\dag}(x,\kappa_{2}m) \\
& = &
q(\kappa_{1}n,\kappa_{1}(n-1))\cdot q^{\dag}(\kappa_{1}(n-1),\kappa_{2}m) +
q(\kappa_{1}n,\kappa_{2}(n+1))\cdot q^{\dag}(\kappa_{2}(n+1),\kappa_{2}m) \\
& = &
\left\{\begin{array}{ll}
\frac{1}{\sqrt{2}}\cdot \frac{1}{\sqrt{2}} +
   \frac{1}{\sqrt{2}} \cdot -\frac{1}{\sqrt{2}} \quad & \mbox{if }n=m \\
0 & \mbox{otherwise} \end{array}\right. \\
& = &
0 \\
& = &
\idmap(\kappa_{1}n,\kappa_{2}m) \\
\lefteqn{\big(q^{\dag} \relafter q\big)(\kappa_{2}n, \kappa_{1}m)} \\
& = &
\sum_{x}q(\kappa_{2}n,x)\cdot q^{\dag}(x,\kappa_{1}m) \\
& = &
q(\kappa_{2}n,\kappa_{1}(n-1))\cdot q^{\dag}(\kappa_{1}(n-1),\kappa_{1}m) +
q(\kappa_{2}n,\kappa_{2}(n+1))\cdot q^{\dag}(\kappa_{2}(n+1),\kappa_{1}m) \\
& = &
\left\{\begin{array}{ll}
\frac{1}{\sqrt{2}}\cdot \frac{1}{\sqrt{2}} +
   -\frac{1}{\sqrt{2}} \cdot \frac{1}{\sqrt{2}} \quad & \mbox{if }n=m \\
0 & \mbox{otherwise} \end{array}\right. \\
& = &
0 \\
& = &
\idmap(\kappa_{2}n,\kappa_{1}m) \\
\lefteqn{\big(q^{\dag} \relafter q\big)(\kappa_{2}n, \kappa_{2}m)} \\
& = &
\sum_{x}q(\kappa_{1}n,x)\cdot q^{\dag}(x,\kappa_{2}m) \\
& = &
q(\kappa_{2}n,\kappa_{1}(n-1))\cdot q^{\dag}(\kappa_{1}(n-1),\kappa_{2}m) +
q(\kappa_{2}n,\kappa_{2}(n+1))\cdot q^{\dag}(\kappa_{2}(n+1),\kappa_{2}m) \\
& = &
\left\{\begin{array}{ll}
\frac{1}{\sqrt{2}}\cdot \frac{1}{\sqrt{2}} +
   -\frac{1}{\sqrt{2}} \cdot -\frac{1}{\sqrt{2}} \quad & \mbox{if }n=m \\
0 & \mbox{otherwise} \end{array}\right. \\
& = &
\idmap(\kappa_{2}n,\kappa_{2}m).
\end{array}$$
}

In the remainder of this section we investigate some of the
categorical structure of the category of bifinite multirelations.  It
will clarify, for instance, that the above ``walks'' map $q$ is an
endomap
$\mathbb{Z}\oplus\mathbb{Z}\rightarrow\mathbb{Z}\oplus\mathbb{Z}$,
where $\oplus$ is a biproduct.

%% The situation is thus fairly similar to Proposition~\ref{LocBifinProp}
%% where we have tame relations given by bifinite relations $X\times
%% Y\rightarrow 2$ factoring both as $X\rightarrow\Powfin(Y)$ and
%% $Y\rightarrow\Powfin(X)$.  In~\cite{Jacobs11a}, in the context of
%% quantum computing, a coalgebra
%% $X\rightarrow\Mlt_{\mathbb{C}}(Y)$---for the complex numbers
%% $\mathbb{C}$ as semiring---is called swappable if the associated
%% (swapped) map $Y\rightarrow \mathbb{C}^{X}$ factors through
%% $Y\rightarrow\Mlt_{\mathbb{C}}(X)$. Clearly, this means that one has a
%% tame relation. The resulting dagger structure provides the
%% reversibility of quantum computing.

\begin{proposition}
\label{BfMRCatProp}
For an involutive commutative semiring $S$, the category
$\BifMRel_{S}$ of sets and bifinite $S$-valued multirelations has
(symmetric) dagger tensors $(\times, 1)$ and dagger biproducts
$(+,0)$, where tensors distribute over biproducts. 

The set of scalars in $\BifMRel_{S}$ (endomaps of the tensor unit $1$)
is $S$. The induced additive structure on homsets is obtained pointwise
from $S$. The homsets are (Abelian) groups iff $S$ is a ring.

The objects $X$ in $\BifMRel_{S}$ that are finite (as a set) are
$S$-modules of the form $S^{n}$ that carry a compact structure.  The
induced monodial trace operation $\tr(s)\colon X \rightarrow Y$, for
$s\colon X \times A \rightarrow Y\times A$ is given by the sum of the
`diagonal' elements:
$$\begin{array}{rcl}
\tr(s)(x,y)
& = &
\displaystyle\sum_{a\in A} \frac{s(x,a,y,a)}{\#A\cdot 1}.
\end{array}$$

\noindent Here we assume that the number of elements $\#A\in\NNO$ is
not zero; otherwise, trivially, $A = 0$ is the zero-object and $s$ is
the zero-map.
\end{proposition}

In case $S$ is a field like $\mathbb{R}$ or $\mathbb{C}$, the latter
category of modules $S^{n}$ is of course the category of
finite-dimensional vector (or Hilbert) spaces.

\begin{myproof}
The tensor is given on objects by Cartesian product: $X_{1}\otimes
X_{2} = X_{1}\times X_{2}$. And if we have $r_{i}\colon
X_{i}\rightarrow Y_{i}$, then $r_{1}\otimes r_{2}\colon X_{1}\otimes
X_{2}\rightarrow Y_{1}\otimes Y_{2}$ is given by the function:
$$\xymatrix@R-2pc{
(X_{1}\times Y_{1}) \times (X_{2}\times Y_{2}) 
   \ar[rr]^-{r_{1}\otimes r_{2}} & & S \\
\tuple{(x_{1},y_{1}), \, (x_{2},y_{2})} \ar@{|->}[rr] & &
   r_{1}(x_{1},y_{1})\cdot r_{2}(x_{2},y_{2}).
}$$

\noindent The singleton set $1$, say $1=\{*\}$, is tensor unit.  The
monoidal (dagger) isomorphisms are given by equalities, such as:
$$\xymatrix@R-2pc@C-1.6pc{
(1\times X)\times X\ar[r]^-{\lambda} & S
& 
(X\times Y)\times (Y\times X)\ar[r]^-{\gamma} & S \\
\tuple{(*,x),x'}\ar@{|->}[r] &
   {\left\{\begin{array}{ll} 1 \; & \mbox{if }x=x' \\
      0 & \mbox{otherwise} \end{array}\right.}
& 
\tuple{(x,y),(y',x')}\ar@{|->}[r] &
   {\left\{\begin{array}{ll} 1 \; & \mbox{if }x=x',y=y' \\
      0 & \mbox{otherwise.} \end{array}\right.}
}$$

\noindent The endomaps on the tensor unit $1$ are maps $1\times
1\rightarrow S$, corresponding to elements of the semiring $S$.

The category $\BifMRel_{S}$ also has biproducts, given on objects by
finite coproducts on sets (whose coprojections we write as $\kappa_i$,
like above). The empty set $0$ is zero object in $\BifMRel_{S}$,
with empty multirelations $X\rightarrow 0$ and $0\rightarrow Y$. The
resulting zero map $0\colon X\rightarrow Y$ is the relation $0\colon
X\times Y\rightarrow S$ that is always $0$.  The coprojections and
projections $X_{i} \stackrel{\kappa_i}{\longrightarrow} X_{1}\oplus
X_{2} \stackrel{\pi_i}{\longrightarrow} X_{i}$ in $\BifMRel_{S}$ are
given by:
$$\xymatrix@R-2pc@C-1.5pc{
X_{i}\times (X_{1}+X_{2})\ar[r]^-{\kappa_i} & S
& 
(X_{1}+X_{2})\times X_{i}\ar[r]^-{\pi_i} & S \\
\tuple{x,u}\ar@{|->}[r] &
   {\left\{\begin{array}{ll} 1 \; & \mbox{if }u=\kappa_{i}x \\
      0 & \mbox{otherwise} \end{array}\right.}
& 
\tuple{u,x}\ar@{|->}[r] &
   {\left\{\begin{array}{ll} 1 \; & \mbox{if }u=\kappa_{i}x \\
      0 & \mbox{otherwise.} \end{array}\right.}
}$$

\noindent (Notice that two different coprojections $\kappa$ occur: in
$\BifMRel_{S}$ and in \Sets.)

We have $\pi_{i} = (\kappa_{i})^{\dag}$ in $\BifMRel_S$. Tuples and
cotuples, for $r_{i}\colon Z\rightarrow X_{i}$ and $t_{i}\colon X_{i}
\rightarrow Z$ are given by:
$$\xymatrix@R-2pc@C-1.8pc{
Z\times (X_{1}+X_{2})\ar[r]^-{\tuple{r_{1},r_{2}}} & S
& 
(X_{1}+X_{2})\times Y\ar[r]^-{[t_{1},t_{2}]} & S \\
\tuple{z,u}\ar@{|->}[r] & r_{i}(z,x), 
   \mbox{ for }u = \kappa_{i}x
& 
\tuple{u,z}\ar@{|->}[r] &
   t_{i}(x,z), \mbox{ for }u = \kappa_{i}x
}$$

\noindent It is not hard to see that $\tuple{r_{1},r_{2}}^{\dag} =
          [(r_{1})^{\dag}, (r_{2})^{\dag}]$.

\auxproof{
We check some equations:
$$\begin{array}{rcl}
\big(\pi_{i} \relafter \tuple{r_{1},r_{2}}\big)(z,x)
& = &
\sum_{u}\tuple{r_{1},r_{2}}(z,u)\cdot \pi_{i}(u,x) \\
& = &
\tuple{r_{i},r_{2}(z,\kappa_{i}x)} \\
& = &
r_{i}(z,x) \\
\big(\tuple{r_{1},r_{2}}\relafter s\big)(y,u)
& = &
\sum_{z}s(y,z)\cdot \tuple{r_{1},r_{2}}(z,u) \\
& = &
\sum_{z}s(y,z)\cdot r_{i}(z,x), \quad
   \mbox{ for }u = \kappa_{i}x \\
& = &
(r_{i}\relafter s)(y,x), \quad
   \mbox{ for }u = \kappa_{i}x \\
& = &
\tuple{r_{1}\relafter s, r_{2}\relafter s}(y,u) \\
\tuple{\pi_{1},\pi_{2}}(w,u)
& = &
\pi_{i}(w,x), \quad \mbox{ for }u = \kappa_{i}x \\
& = &
{\left\{\begin{array}{ll} 1 \; & \mbox{if }w=\kappa_{i}x \\
      0 & \mbox{otherwise.} \end{array}\right.} \qquad
   \mbox{ for }u = \kappa_{i}x \\
& = &
{\left\{\begin{array}{ll} 1 \; & \mbox{if }w=u \\
      0 & \mbox{otherwise.} \end{array}\right.} \\
& = &
\idmap(w,u) \\
\tuple{r_{1},r_{2}}^{\dag}(u,z)
& = &
\overline{\tuple{r_{1},r_{2}}(z,u)} \\
& = &
\overline{r_{i}(z,x)}, \quad \mbox{ for }u = \kappa_{i}x \\
& = &
(r_{i})^{\dag}(x,z), \quad \mbox{ for }u = \kappa_{i}x \\
& = &
[(r_{1})^{\dag}, (r_{2})^{\dag}](u,z).
\end{array}$$

Notice that $r_{1}\oplus r_{2}$ can be described in two ways as:
$$\begin{array}{rcl}
\tuple{r_{1}\relafter \pi_{1}, r_{2} \relafter \pi_{2}}
   (\kappa_{i}x, \kappa_{j}z) 
& = &
(r_{j} \relafter \pi_{j})(\kappa_{i}x, z) \\
& = &
\sum_{x'}\pi_{j}(\kappa_{i}x,x')\cdot r_{j}(x',z) \\
& = &
r_{i}(x,z), \quad \mbox{if }i=j \\
{[\kappa_{1}\relafter r_{1}, \kappa_{2}\relafter r_{2}]}
   (\kappa_{i}x, \kappa_{j}z)
& = &
(\kappa_{i}\relafter r_{i})(x, \kappa_{j}z) \\
& = &
\sum_{z'}r_{i}(x,z')\cdot \kappa_{i}(z',\kappa_{j}z) \\
& = &
r_{i}(x,z), \quad\mbox{if }i=j.
\end{array}$$

The induced additive structure on homsets of maps $X\rightarrow Y$ is
pointwise:
$$\begin{array}{rcl}
(r+s)(x,y)
& \smash{\stackrel{\textrm{def}}{=}} &
\big(\nabla \relafter (r\oplus s) \relafter \Delta\big)(x,y) \\
& = &
\big(\nabla \relafter (r\oplus s)\big)(\kappa_{1}x,y) +
   \big(\nabla \relafter (r\oplus s)\big)(\kappa_{2}x,y) \\
& = &
(r\oplus s)(\kappa_{1}x,\kappa_{1}y) + 
   (r\oplus s)(\kappa_{1}x,\kappa_{2}y)\; + \\
& & \qquad
   (r\oplus s)(\kappa_{2}x,\kappa_{1}y) + 
   (r\oplus s)(\kappa_{2}x,\kappa_{2}y) \\ 
& = &
r(x,y) + 0 + 0 + s(x,y) \\
& = &
r(x,y) + s(x,y).
\end{array}$$
}

There are distributivity (dagger) isomorphisms $X\otimes (Y_{1}\oplus
Y_{2}) \xrightarrow{\cong} (X\otimes Y_{i})\oplus (X\otimes Y_{2})$
given by:
$$\xymatrix@R-2pc@C-1.5pc{
\Big(X\times (Y_{1}+Y_{2})\Big)\times 
  \Big((X\times Y_{1}) + (X\times Y_{2})\Big)\ar[r] & S \\
\tuple{(x,u),v}\ar@{|->}[r] &
   {\left\{\begin{array}{ll} 1 \; & 
      \mbox{if }u=\kappa_{i}y\mbox{ and }v=\kappa_{i}(x,y) \\
      0 & \mbox{otherwise.} \end{array}\right.}
}$$

\auxproof{
We check naturality of this isomorphism (say with name $d$).
$$\begin{array}{rcl}
\lefteqn{\big(d \relafter (r\otimes (t_{1}\oplus t_{2}))\big)
   (\tuple{a,\kappa_{i}z},\kappa_{j}(x,y))} \\
& = &
\sum_{u}(r\otimes (t_{1}\oplus t_{2}))(\tuple{a,\kappa_{i}z},
   \tuple{x,u}) \cdot d(\tuple{x,u},\kappa_{j}(x,y)) \\
& = &
(r\otimes (t_{1}\oplus t_{2}))
   (\tuple{a,\kappa_{i}z}, \tuple{x, \kappa_{j}y}) \\
& = &
r(a,x) \cdot (t_{1}\oplus t_{2})(\kappa_{i}z, \kappa_{j}y) \\
& = &
r(a,x)\cdot t_{i}(z,y), \quad \mbox{if }i=j \\
& = &
(r\otimes t_{i})((a,z), (x,y)), \quad \mbox{if }i=j \\
& = &
\big((r\otimes t_{1})\oplus (r\otimes t_{2})\big)
   (\kappa_{i}(a,z), \kappa_{j}(x,y)) \\
& = &
\sum_{v} d(\tuple{a,\kappa_{i}z}, v) \cdot
   \big((r\otimes t_{1})\oplus (r\otimes t_{2})\big)(v, \kappa_{j}(x,y)) \\
& = &
\big((r\otimes t_{1})\oplus (r\otimes t_{2}) \relafter d\big)
   (\tuple{a,\kappa_{i}z},\kappa_{j}(x,y)).
\end{array}$$
}

Finally we note that if $X$ is a finite set, say with $n = |X|$
elements, then $\Mlt_{S}(X) \cong S^{n}$. Further, each multirelation
$r\colon X\times Y\rightarrow S$ is automatically bifinite, if $X,Y$
are finite. Such a morphism is thus determined by the associated map
$X\rightarrow \Mlt_{S}(Y) \cong S^{|Y|}$. The latter corresponds to a
linear map $S^{|X|}\rightarrow S^{|Y|}$. The compact structure on
a finite set $X$ has $X^{*}=X$ with unit $\eta\colon 1\rightarrow
X^{*}\otimes X$ and counit $\varepsilon \colon X\otimes
X^{*}\rightarrow 1$ given by:
$$\xymatrix@R-1.8pc@C-.5pc{
1\times (X\times X)\ar[r]^-{\eta} & S
&
(X\times X)\times 1\ar[r]^-{\varepsilon} & S \\
\tuple{*,(x,x')}\ar@{|->}[r] & {\left\{\begin{array}{ll}
   1 \; & \mbox{if }x=x' \\ 0 & \mbox{otherwise} \end{array}\right.}
&
\tuple{(x,x'),*}\ar@{|->}[r] & {\left\{\begin{array}{ll}
   1 \; & \mbox{if }x=x' \\ 0 & \mbox{otherwise.} \end{array}\right.}
}$$

\noindent These multirelations are bifinite because $X$ is finite.
The formula for traces is obtained in the standard manner from compact
structure $\eta, \varepsilon$, see~\cite{JoyalSV96}. \QED

\auxproof{
We check the trace requirements:
\begin{itemize}
\item Yanking: for the swap map $\gamma \colon A\times A \congrightarrow 
A\times A$, assuming $A$ is non-empty:
$$\begin{array}{rcl}
\tr(\gamma)(x,y)
& = &
\frac{\sum_{a}{\gamma(x,a,y,a)}}{\#A \cdot 1} \\
& = &
\sum_{a}{\left\{\begin{array}{ll}
   \frac{1}{\#A \cdot 1} \quad & \mbox{if }x=y \\ 0 \quad \mbox{otherwise} 
   \end{array}\right.} \\
& = &
\left\{\begin{array}{ll}
   1 \quad & \mbox{if }x=y \\ 0 \quad \mbox{otherwise} 
   \end{array}\right. \\
& = &
\idmap[A](x,y).
\end{array}$$

\noindent (If $A=\emptyset$, the zero object, then $\tr(s) 
= \idmap[A]$, trivially.

\item Vanishing I: for $s\colon X\times 1 \rightarrow Y\times 1$
we get:
$$\begin{array}{rcl}
\tr(s)(x,y)
& = &
\sum_{a\in 1}\frac{s(x,a,y,a)}{\#1\cdot 1} \\
& = &
s(x, *, y, *) \\
& = &
\big(\lambda \after s \after \lambda^{-1}\big)(x,y).
\end{array}$$

\item Vanishing II: for $s\colon (X\times A)\times B \rightarrow
(X\times A)\times B$, 
$$\begin{array}{rcl}
\tr(\tr(s))(x,y)
& = &
\sum_{a}\frac{\tr(s)(x, a, y, a)}{\#A\cdot 1} \\
& = &
\sum_{a}\sum_{b}\frac{s(x, a, b, y, a, b)}
   {(\#B \cdot 1) \cdot (\#A\cdot 1)} \\
& = &
\sum_{(a,b)}\frac{s(x, a, b, y, a, b)}
   {\#(A\times B) \cdot 1} \\
& = &
\tr(s)
\end{array}$$

\item Naturality: for $u\colon Z\rightarrow X$ and $t\colon Y
\rightarrow W$,
$$\begin{array}{rcl}
\lefteqn{\tr\big((t\times\idmap)\after s \after 
   (u\times\idmap)\big)(z,w)} \\
& = &
\sum_{a}\frac{\big((\idmap\times t)\after r \after (\idmap\times s)\big)
   (z,a,w,a)}{\#A\cdot 1} \\
& = &
\sum_{x,y,a}\frac{u(z,x) \cdot s(x,a,y,a) \cdot t(y,w)}{\#A\cdot 1} \\
& = &
\sum_{x,y}u(z,x) \cdot (\sum_{a}\frac{s(x,a,y,a)}{\#A\cdot 1}) 
   \cdot t(y,w) \\
& = &
\sum_{x,y} u(z,x) \cdot \tr(s)(x,y) \cdot t(y,w) \\
& = &
\big(u \after \tr(s) \after t\big)(z,w).
\end{array}$$

\item Superposition:
$$\begin{array}{rcl}
\tr(\idmap\times s)(z,x,w,y)
& = &
\sum_{a}\frac{(\idmap\times s)(z, x, a, w, y, a)}{\#A\cdot 1} \\
& = &
\sum_{a} \frac{\idmap(z,w) \cdot s(x,a,y,a)}{\#A\cdot 1} \\
& = &
\idmap(z,w) \cdot \sum_{a}\frac{s(x,a,y,a)}{\#A\cdot 1} \\
& = &
\idmap(z,w) \cdot \tr(s)(x,y) \\
& = &
(\idmap\times \tr(s))(z,x,w,y).
\end{array}$$

\item Exchange:
$$\begin{array}{rcl}
\lefteqn{\tr(\tr((\idmap\times \gamma) \after s \after 
   (\idmap\times \gamma)))(x,y)} \\
& = &
\sum_{a}\frac{(\tr((\idmap\times \gamma) \after s \after 
   (\idmap\times \gamma)))(x,a,y,a)}{\#A\cdot 1} \\
& = &
\sum_{a} \frac{\sum_{b}
   \frac{((\idmap\times \gamma) \after s \after (\idmap\times \gamma))
   (x,a,b,y,a,b)}{\#B\cdot 1}}{\#A\cdot 1} \\
& = &
\sum_{a}\sum_{b}\frac{s(x,b,a,y,b,a)}{(\#B\cdot 1)\cdot (\#A\cdot 1)} \\
& = &
\sum_{b}\sum_{a}\frac{s(x,b,a,y,b,a)}{(\#A\cdot 1)\cdot (\#B\cdot 1)} \\
& = &
\sum_{b} \frac{\sum_{a}
   \frac{s(x,a,b,y,a,b)}{\#A\cdot 1}}{\#B\cdot 1} \\
& = &
\sum_{b}\frac{\tr(s)(x,b,y,b)}{\#B\cdot 1} \\
& = &
\tr(\tr(s)).
\end{array}$$

\end{itemize}
}

\end{myproof}

%% Interestingly, if one adapts the multiset monad to allow countable
%% support---as also mentioned at the end of
%% Subsection~\ref{MltDstSubsec}---countable sets/objects in \BifMRel
%% carry a compact structure.

%% \begin{example}
%% \label{BitEx}
%% Consider a map $q\colon 1 \rightarrow 1\oplus 1 = 2$ in $\BifMRel_{S}$
%% that is a dagger mono---\textit{i.e.}~satisfies $q^{\dag} \after q =
%% \idmap$. This $q\colon 1\times 2 \rightarrow S$ corresponds to two
%% scalars $a = q(*,0)\in S$ and $b = q(*,1)\in S$, where we use $1 =
%% \{*\}$ and $2 = \{0,1\}$. The dagger mono requirement yields:
%% $$\textstyle 1 
%% = 
%% (q^{\dag} \relafter q)(*,*) 
%% = 
%% \sum_{x\in 2}q(*,x)\cdot \overline{q(*,x)}
%% = 
%% a\cdot\overline{a} + b\cdot\overline{b}.$$

%% \noindent This means
%% that such a dagger mono $q\colon 1\rightarrow 2$ is:
%% \begin{itemize}
%% \item a \emph{quantum} bit for $S=\mathbb{C}$, since
%% $q = (a,b)\in\mathbb{C}^{2}$ satisfies $\|a\|^{2} + \|b\|^{2} = 1$;
%% \item a \emph{classical} bit for $S = 2$, since $q = (a,b)\in 2^{2}$
%%   satisfies $a^{2}+b^{2} = a + b = 1$, and thus $a = 1,b=0$ or
%%   $a=0,b=1$.
%% \end{itemize}

%% \noindent Hence these categories $\BifMRel_{S}$ of bifinite
%% multirelations provide a uniform setting for both classical and
%% quantum bits in (reversible) computations.
%% \end{example}

The logical structure of the category $\BifMRel_{S}$ will be described
in terms of its dagger kernels, following~\cite{HeunenJ10a}.  We first
borrow some more terminology from linear algebra. Two multisets
$\varphi,\psi\in\Mlt_{S}(X)$ will be called orthogonal, written as
$\varphi\orthogonal\psi$, if $\cpMlt(\varphi,\psi) = 0$. Recall that
this corresponds to the usual condition
$\sum_{x}\overline{\varphi(x)}\cdot\psi(x) = 0$, for the ``dot'' inner
product. A subset $V\subseteq\Mlt_{S}(X)$ will be called orthogonal if
all its pairs of (different) elements are orthogonal; it will be
called orthonormal if additionally each $\varphi\in V$ satisfies
$\smash{\|\varphi\|^{2} \stackrel{\textrm{\tiny def}}{=}
  \cpMlt(\varphi,\varphi) = 1}$.

%% Given a subset $V\subseteq\Mlt_{S}(X)$, the span
%% $\tuple{V}\subseteq\Mlt_{S}(X)$ is the least subspace containing $V$;
%% it is given by all finite linear combinations of elements of $V$.

Since dagger kernels are both kernels and dagger monos, the following
characterisation sheds light on the situation.

\begin{lemma}
\label{BfMRelDagMonoLem}
A morphism $r\colon X\rightarrow Y$ in $\BifMRel_{S}$ is a dagger
mono iff the set of multisets $\set{r(x,-)}{x\in X}\subseteq\Mlt_{S}(Y)$
is orthonormal.
\end{lemma}

\begin{myproof}
The crucial point is:
$$\textstyle (r^{\dag} \relafter r)(x,x')
=
\sum_{y}r(x,y)\cdot\overline{r(x',y)}
=
\cpMlt\big(r(x',-), r(x,-)\big).$$

\noindent Thus:
$$\begin{array}[b]{rcl}
r \mbox{ is dagger mono}
& \Longleftrightarrow &
r^{\dag} \relafter r = \idmap \\
& \Longleftrightarrow &
\all{x,x'}{(r^{\dag} \relafter r)(x,x') =
  \left\{\begin{array}{ll} 1 \; & \mbox{if }x=x' \\
     0 & \mbox{otherwise} \end{array}\right.} \\
& \Longleftrightarrow &
\all{x,x'}{\cpMlt\big(r(x',-), r(x,-)\big) = 
  \left\{\begin{array}{ll} 1 \; & \mbox{if }x=x' \\
     0 & \mbox{otherwise} \end{array}\right.} \\
& \Longleftrightarrow &
\set{r(x,-)}{x\in X} \mbox{ is orthonormal.}
\end{array}\eqno{\QEDbox}$$
\end{myproof}

Before giving the general construction of dagger kernels, it may be
helpful to see an illustration first.

\begin{example}
\label{KernelEx}
We use $S=\mathbb{R}$ as semiring (actually as field) and start from a
morphism $r\colon\NNO\rightarrow\NNO$ in $\BifMRel_{S}$, described
as the following multirelation $r\colon\NNO\times\NNO \rightarrow S$.
$$\begin{array}{rcl}
r(x,y)
& = &
\left\{\begin{array}{ll}
1 \quad & \mbox{if } x = 2y \\
-1 \quad & \mbox{if } x=2y+1 \\
0 & \mbox{otherwise.}
\end{array}\right.
\end{array}$$

\noindent This $r$ involves an infinite number of multisets:
$$r(-,0) = 1\cdot 0 + (-1)\cdot 1, 
\quad
r(-,1) = 1\cdot 2 + (-1)\cdot 3,
\quad
r(-,2) = \ldots\quad\textit{etc.}$$

\noindent We illustrate how to interpret them as infinitely many
linear equations:
$$v_{0} = v_{1}
\qquad
v_{2} = v_{3}
\qquad
v_{4} = v_{5}
\qquad \textit{etc.}$$

\noindent Assume we have a map $f\colon 1\rightarrow \NNO$ with
$r\relafter f = 0$. Then, for each $y\in\NNO$,
$$\textstyle 0 = (r\relafter f)(*,y) = \sum_{x}f(x)\cdot r(x,y) =
   f(2y) - f(2y+1).$$

\noindent Thus, this $f$, as function $f\colon \NNO\rightarrow S$ with
finite support, satisfies $f(2y) = f(2y+1)$. It thus provides a
``solution'' $f(0) = f(1), f(2) = f(3), \ldots$ to the ``equations''
$r(-,y)=0$.

We wish to describe the dagger kernel of $r$ as the solution space for
these equations $r(-,y)=0$. Lemma~\ref{BfMRelDagMonoLem} tells that we
have to look for an orthonormal basis for this space. An obvious
choice for such a basis is the infinite set of multisets:
$$\begin{array}{rclcrcl}
B
& = &
\set{\varphi_{i}}{i\in\NNO}
& \quad\mbox{where}\quad &
\varphi_{i}
& = &
\frac{1}{\sqrt{2}}\cdot (2i) + \frac{1}{\sqrt{2}}\cdot (2i+1)
   \,\in\,\Mlt_{\mathbb{R}}(X).
\end{array}$$

\noindent We take as kernel the map $\ker(r)\colon B\rightarrow\NNO$,
given as function $\ker(r)\colon B\times\NNO\rightarrow S$ simply by:
$$\begin{array}{rcl}
k(\varphi_{i},x)
& = &
\varphi_{i}(x).
\end{array}$$

\noindent Clearly, this is well-defined, in the sense that $\ker(r)$ is
bifinite, as multirelation. Further, $\ker(r)$ satisfies the appropriate
properties:
$$\begin{array}{rcl}
(r\relafter \ker(r))(\varphi_{i}, y)
& = &
\sum_{x}\ker(r)(\varphi_{i},x)\cdot r(x,y) \\
& = &
\varphi_{i}(2y)\cdot 1 + \varphi_{i}(2y+1)\cdot -1 \\
& = &
\left\{\begin{array}{ll}
\frac{1}{\sqrt{2}} + -\frac{1}{\sqrt{2}} \quad
  & \mbox{if }i=y \\
0 & \mbox{otherwise} \end{array}\right. \\
& = &
0 \\
(\ker(r)^{\dag} \after \ker(r))(\varphi_{i},\varphi_{j})
& = &
\sum_{x}\ker(r)(\varphi_{i},x)\cdot \overline{\ker(r)(\varphi_{j},x)} \\
& = &
\left\{\begin{array}{ll}
(\frac{1}{\sqrt{2}})^{2} + (\frac{1}{\sqrt{2}})^{2} \quad
   & \mbox{if }i=j \\
0 & \mbox{otherwise} \end{array}\right. \\
& = &
\left\{\begin{array}{ll}
1 \quad & \mbox{if }\varphi_{i} = \varphi_{j} \\
0 & \mbox{otherwise} \end{array}\right. \\
& = &
\idmap(\varphi_{i},\varphi_{j}).
\end{array}$$

Next assume we have a map $t\colon Z\rightarrow\NNO$ in $\BifMRel_{S}$
satisfying $r\relafter t = 0$. We have to show that $t$ factors
through the kernel $\ker(r)$. For each $z\in Z$ and $y\in\NNO$ we have
$0 = (r\relafter t)(z,y) = \sum_{x}t(z,x)\cdot r(x,y) = t(z,2y)\cdot 1
+ t(z,2y+1)\cdot -1$. Hence $t(z,2y)=t(z,2y+1)$, so that $t$ solves
the equations $r(-,y)=0$. Since $t$ is bifinite, there are for a fixed
$z\in Z$, only finitely many $y$ with $t(z,y)\neq 0$. Hence we can
express $t(z,-)\in\Mlt_{\mathbb{R}}(\NNO)$ in terms of the base
vectors in $B$, say as:
$$\begin{array}{rcl}
t(z,-)
& = &
a_{1}\cdot \varphi_{y_1} + \cdots + a_{n}\cdot \varphi_{y_n},
\qquad\mbox{where}\qquad
a_{i} = t(z,2y_{i})\cdot\sqrt{2} \in \mathbb{R},
\end{array}$$

\noindent for certain $y_{1},\ldots,y_{n}\in\NNO$ (depending on $z$).
We thus define the required map $t'\colon Z\rightarrow B$ by
$t'(z,\varphi_{y_i}) = a_{i}$, for these $y_{1},\ldots,y_{n}$ (and $0$
elsewhere). Then:
$$\begin{array}{rcl}
(\ker(r) \relafter t')(z,x)
& = &
\sum_{i}t'(z,\varphi_{i})\cdot \ker(r)(\varphi_{i},x) \\
& = &
\sum_{i}a_{i}\cdot\varphi_{y_i}(x) \\
& = &
t(z,x).
\end{array}$$

\auxproof{
We also have an example with a finite number of equations.
We illustrate how to compute a kernel for a specific map
$r\colon\NNO\rightarrow\NNO$ in $\BifMRel_{S}$, where for instance
$S$ is the field $\mathbb{R}$ or $\mathbb{C}$. We assume some rather
arbitrary values:
$$\begin{array}{rclcrclcrclcrcl}
r(0,0) & = & -1
& \qquad &
r(4,0) & = & 2
& \qquad &
r(0,1) & = & 1
& \qquad &
r(1,1) & = & 1  \\
r(3,1) & = & -1
& &
r(0,2) & = & 2
& &
r(1,2) & = & 1
& &
r(2,2) & = & 1
\end{array}$$

\noindent Alternatively, we may describe this bifinite multirelation $r$
via the following three multisets.
$$\begin{array}{rcl}
r(-,0)
& = &
(-1)\cdot 0 + 2\cdot 4 \\
r(-,1)
& = &
1\cdot 0 + 1\cdot 1 + (-1)\cdot 3 \\
r(-,2)
& = &
2\cdot 0 + 1\cdot 1 + 1\cdot 2.
\end{array}$$

\noindent In order to obtain the kernel of $r$ we shall view these
multisets as describing relations in variables $x_{0}, x_{1}, x_{2},
x_{3}, x_{4}$, corresponding to the elements of the multisets. The
equations that we thus consider are:
\begin{equation}
\label{KernelExEqn}
\begin{array}{rcl}
0
& = &
-x_{0} + 2 x_{4} \\
0
& = &
x_{0} + x_{1} - x_{3} \\
0
& = &
2 x_{0} + x_{1} + x_{2}.
\end{array}
\end{equation}

Assume a morphism $t\colon Z\rightarrow \NNO$ with $r\relafter t = 0$
in $\BifMRel_{S}$; we need to construct a kernel $\ker(r)$ of $r$ so
that each such $t$ factors through $\ker(r)$. For each $x\in X$ and
$y=0,1,2$ we have $0 = (r\relafter t)(z,y) = \sum_{x} t(z,x)\cdot
r(x,y)$. Hence the multisets $t(z,-)$ also solves the
equations~\eqref{KernelExEqn}.

\noindent The first two equations give us $x_{0} = 2x_{4}$ and $x_{0}
= -x_{1}+x_{3}$. Substituting the latter in the third equation gives
$-2x_{1} + 2x_{3} + x_{1} + x_{2} = 0$, and thus $x_{1} = x_{2} +
2x_{3}$. The solution space thus has 2 dimensions. We choose an
orthonormal basis, such as:
$$\frac{1}{4}(-2 , 3, 1, 1, -1)
\qquad\qquad
\frac{1}{4\sqrt{39}}(-6,-7,19,-13,-3).$$

\noindent These vectors may also be understood as multisets over
$\NNO$, as in:
$$\begin{array}{rcl}
\varphi_{1}
& = &
(-\frac{1}{2})\cdot 0 + \frac{3}{4}\cdot 1 + \frac{1}{4}\cdot 2 + 
   \frac{1}{4}\cdot 3 + (-\frac{1}{4})\cdot 4 \\
\varphi_{2}
& = &
(-\frac{3}{2\sqrt{39}})\cdot 0 + (-\frac{7}{4\sqrt{39}})\cdot 1 +
  \frac{19}{4\sqrt{39}}\cdot 2 + (-\frac{13}{4\sqrt{39}})\cdot 3 +
  (-\frac{3}{4\sqrt{39}})\cdot 4.
\end{array}$$

\noindent Orthonormality now amounts to:
$$\cpMlt(\varphi_{1},\varphi_{1}) = 1
\qquad
\cpMlt(\varphi_{1},\varphi_{2}) = 0
\qquad
\cpMlt(\varphi_{2},\varphi_{2}) = 1$$

\noindent where $\cpMlt$ is as introduced in~\eqref{MltComparEqn}.

We are finally in a position to describe the kernel of $r\colon
\NNO\rightarrow\NNO$. The relevant object is:
$$\begin{array}{rcl}
\Ker(r)
& = &
\Big(\NNO - \{0,1,2,3,4\}\Big) \cup \{\varphi_{1},\varphi_{2}\},
\end{array}$$

\noindent with map $\ker(r)\colon\Ker(r)\rightarrow\NNO$, as
function $\ker(r)\colon \Ker(r)\times\NNO\rightarrow S$, given by:
$$\begin{array}{rclcrcl}
\ker(r)(n,m)
& = &
\left\{\begin{array}{ll}
1 \; & \mbox{if }n>4\mbox{ and }n=m \\
0 \; & \mbox{otherwise} \end{array}\right.
& \qquad\quad &
\ker(\varphi_{i},m)
& = &
\varphi_{i}(m)
\end{array}$$

\noindent Then, for $n>4$,
$$\begin{array}{rcl}
(r \relafter \ker(r))(n,k)
& = &
\sum_{m}\ker(r)(n,m)\cdot r(m,n) \\
& = &
r(n,k) \\
& = &
0 \qquad \mbox{since $r(n,k)=0$ for $n>4$} \\
(r \relafter \ker(r))(\varphi_{i},k) 
& = &
\sum_{m} \ker(r)(\varphi_{i},m)\cdot r(m,k) \\
& = &
\sum_{m} \varphi_{i}(m) \cdot r(m,k) \\
& = &
0 \qquad \mbox{since $\varphi_i$ solves~\eqref{KernelExEqn}.}
\end{array}$$

\noindent It is not hard to see that $\ker(r)$ is a dagger mono,
\textit{i.e.}~satisfies $\ker(r)^{\dag} \after \ker(r) = \idmap$,
precisely because the $\varphi_{i}$ form an orthonormal basis. The
crucial case is:
$$\begin{array}{rcl}
(\ker(r)^{\dag} \relafter \ker(r))(\varphi_{i},\varphi_{j})
& = &
\sum_{m} \ker(r)(\varphi_{i},m)\cdot \ker(r)(\varphi_{j},m) \\
& = &
\sum_{m} \varphi_{i}(m) \cdot \varphi_{j}(m) \\
& = &
\cpMlt(\varphi_{i}, \varphi_{j}) \\
& = &
\left\{\begin{array}{ll} 1 \; & \mbox{if }i=j \\
   0 & \mbox{otherwise} \end{array}\right. \\
& = &
\idmap[\Ker(r)](\varphi_{i},\varphi_{j}).
\end{array}$$

\auxproof{
The other cases: for $n,n'>4$,
$$\begin{array}{rcl}
(\ker(r)^{\dag} \relafter \ker(r))(n,n')
& = &
\sum_{m} \ker(r)(n,m)\cdot \ker(r)(n',m) \\
& = &
\left\{\begin{array}{ll} 1 \; & \mbox{if }n=n' \\
   0 & \mbox{otherwise} \end{array}\right. \\
& = &
\idmap[\Ker(r)](n,n') \\
(\ker(r)^{\dag} \relafter \ker(r))(n,\varphi_{i})
& = &
\sum_{m} \ker(r)(n,m)\cdot \ker(r)(\varphi_{i},m) \\
& = &
\ker(r)(\varphi_{i},n) \\
& = &
\varphi_{i}(n) \\
& = &
0 \\
& = &
\idmap[\Ker(r)](n,\varphi_{i}).
\end{array}$$
}

Now suppose we have $t\colon X\rightarrow \NNO$ with $r\relafter t =
0$. Then for each $x\in X$ and $k=0,1,2$ we have $0 = (r\relafter
t)(x,k) = \sum_{m} t(x,m)\cdot r(m,k)$. Hence the multiset $t(x,-)$
also solves~\eqref{KernelExEqn}. When restricted to the subset
$\{0,1,2,3,4\}\subseteq \NNO$, these $t(x,-)$ are thus elements of
the solution space spanned by $\varphi_{1},\varphi_{2}$. Hence we can
write, for each $x\in X$, an equation of multisets:
$$\begin{array}{rcl}
\sum_{n\leq 4}t(x,n)\cdot n
& = &
f_{1}(x)\cdot \varphi_{1} + f_{2}(x)\cdot \varphi_{2},
\end{array}$$

\noindent where $f_{i}\colon X\rightarrow S$. Thus we define
$t'\colon X\rightarrow \Ker(r)$ as:
$$\begin{array}{rclcrcl}
t'(x,n)
& = &
t(x,n)
& \qquad &
t'(x,\varphi_{i})
& = &
f_{i}(x).
\end{array}$$

\noindent This $t'$ is the required, necessarily unique (because
$\ker(r)$ is mono) map, since:
$$\begin{array}{rcl}
\lefteqn{\big(\ker(r) \relafter t'\big)(x,m)} \\
& = &
\Big(\sum_{n>4}t'(x,n)\cdot \ker(r)(n,m)\Big) + 
   \Big(\sum_{i}t'(x,\varphi_{i})\cdot \ker(r)(\varphi_{i},m)\Big) \\
& = &
\left\{\begin{array}{ll} t'(x,m) & \mbox{if }m>4 \\
   \sum_{i}f_{i}(x)\cdot \varphi_{i}(m)\quad 
      & \mbox{if } m\leq 4 \end{array}\right. \\
& = &
\left\{\begin{array}{ll} t(x,m) & \mbox{if }m>4 \\
   t(x,m)\quad & \mbox{if } m\leq 4 \end{array}\right. \\
& = &
t(x,m).
\end{array}$$
}
\end{example}

\begin{proposition}[AC]
\label{BfMRelKernelProp}
The category $\BifMRel_{S}$, restricted to countable objects, has
dagger kernels, assuming $S=\mathbb{R}$ or $S=\mathbb{C}$.
\end{proposition}

\begin{myproof}
For an arbitrary map $r\colon X\rightarrow Y$ we consider, like in
Example~\ref{KernelEx}, the multisets $r(-,y)\in\Mlt_{S}(X)$ as
equations, whose solutions, also in $\Mlt_{S}(X)$, give rise to
kernels. The support $\support(r(-,y)) = \set{x}{r(x,y)\neq 0}$ of
such an equation captures the variables that occur. We first collect
all such variables in a subset $X_{r}\subseteq X$, and then describe
the set of solutions in terms of multisets over these variables.
$$\begin{array}{rcl}
X_{r}
& = &
\bigcup_{y\in Y}\support\big(r(-,y)\big) \\
\mathrm{Sol}_{r}
& = &
\setin{\varphi}{\Mlt_{S}(X_{r})}{\allin{y}{Y}{
   \sum_{x}\varphi(x)\cdot r(x,y)=0}}.
\end{array}$$

\noindent Clearly, $\mathrm{Sol}_{r} \subseteq \Mlt_{S}(X)$ is a
linear subspace. Hence, using the Axiom of Choice, we can choose a
basis $B_{r} \subseteq \mathrm{Sol}_{r}$, of linearly independent,
with norm 1. Since the set $\set{1x}{x\in X}$ is a countable basis
for $\Mlt_{S}(X)$, $B_r$ has at most countably many elements.

We claim: for each $z\in X$, the set
$\setin{\varphi}{B_r}{\varphi(z)\neq 0}$ is finite. Suppose not,
\textit{i.e.}~suppose there are infinitely many $\varphi_{i}\in B_{r}$
with $\varphi_{i}(z)\neq 0$. Since $\varphi_{i}\in\Mlt_{S}(X_{r})$ and
$z\in\support(\varphi_{i})\subseteq X_{r}$, there must be an $y_{i}\in
Y$ with $r(z,y_{i}) \neq 0$. Because $r$ is bifinite there can only be
finitely many such $y_{i}$, say $y_{1}, \ldots, y_{n}$. Since the
$\varphi_{i}$ are in $B_{r} \subseteq \mathrm{Sol}_{r}$, we have
$\sum_{x}\varphi_{i}(x) \cdot r(x,y_{j})=0$ for each $i$ and $j\leq
n$. The solution space of these $n$ equations $r(-,y_{j})$ has finite
dimension. Hence it cannot contain infinitely many linearly
independent $\varphi_i$.

We now define a kernel object $\Ker(r) = \big(X-X_{r}\big)\cup B_{r}$,
with kernel map $\ker(r)\colon \Ker(r)\rightarrow X$ given by:
$$\begin{array}{rclcrcl}
\ker(r)(x,x')
& = &
\left\{\begin{array}{ll}
1 \; & \mbox{if }x\in X-X_{r}\mbox{ and }x=x' \\
0 \; & \mbox{if }x\in X-X_{r}\mbox{ and }x\neq x' \end{array}\right.
& \qquad &
\ker(\varphi,x)
& = &
\varphi(x).
\end{array}$$

\noindent This gives a bifinite multirelation by the claim above.
We check that this $\ker(r)$ is a dagger kernel in three steps.
\begin{itemize}
\item[-] In order to obtain that $\ker(r)$ is a dagger mono by
  applying Lemma~\ref{BfMRelDagMonoLem} we need to transform the set
  of base vectors $B_{r} \subseteq \mathrm{Sol}_{r}$ into an
  orthonormal basis. This can be done in a standard way, via
  Gram-Schmidt, using the inner product $\cpMlt$
  from~\eqref{MltComparEqn}. Because $B_r$ is countable, we can write
  $B_{r} = \set{\varphi_n}{n\in\NNO}$ and replace each $\varphi_n$ by
  $\varphi'_n$ obtained as: $\varphi'_{0} = \varphi_{0}$ and:
$$\begin{array}{rcl}
\varphi'_{n+1}
& = &
\varphi_{n+1} - \cpMlt(\varphi'_{0}, \varphi_{n+1})\varphi'_{0} - \cdots
    - \cpMlt(\varphi'_{n}, \varphi_{n+1})\varphi_{n}.
\end{array}$$

\noindent By construction, $\cpMlt(\varphi'_{k}, \varphi'_{n}) = 0$,
for $k<n$. Hence we may assume that $B_r$ is orthonormal.

\auxproof{
Then, for $k\leq n+1$,
$$\begin{array}{rcl}
\lefteqn{\cpMlt(\varphi'_{k}, \varphi'_{n+1})} \\
& = &
\cpMlt(\varphi'_{k}, \varphi_{n+1} - 
   \cpMlt(\varphi'_{0}, \varphi_{n+1})\varphi'_{0} - \cdots
    - \cpMlt(\varphi'_{n}, \varphi_{n+1})\varphi_{n}) \\
& = &
\cpMlt(\varphi'_{k}, \varphi_{n+1}) - 
   \cpMlt(\varphi'_{k}, 
   \cpMlt(\varphi'_{0}, \varphi_{n+1})\varphi'_{0}) - \cdots \\
& & \qquad
    - \; \cpMlt(\varphi'_{k}, 
   \cpMlt(\varphi'_{n}, \varphi_{n+1})\varphi_{n}) \\
& = &
\cpMlt(\varphi'_{k}, \varphi_{n+1}) - 
   \cpMlt(\varphi'_{0}, \varphi_{n+1})\cpMlt(\varphi'_{k}, 
   \varphi'_{0}) - \cdots \\
& & \qquad
    -\; \cpMlt(\varphi'_{n}, \varphi_{n+1})\cpMlt(\varphi'_{k}, 
   \varphi_{n}) \\
& \smash{\stackrel{\text{(IH)}}{=}} &
\cpMlt(\varphi'_{k}, \varphi_{n+1}) - 
   \cpMlt(\varphi'_{k}, \varphi_{n+1})\cpMlt(\varphi'_{k}, \varphi'_{k}) \\
& = &
\cpMlt(\varphi'_{k}, \varphi_{n+1}) - 
   \cpMlt(\varphi'_{k}, \varphi_{n+1})\|\varphi'_{k}\|^{2} \\
& = &
\cpMlt(\varphi'_{k}, \varphi_{n+1}) - 
   \cpMlt(\varphi'_{k}, \varphi_{n+1}) \\
& = &
0.
\end{array}$$
}

\item[-] We have $r \relafter \ker(r) = 0$, since for $x\in X-X_{r}$
  and $\varphi\in B_{r}$,
$$\begin{array}{rcl}
\big(r \relafter \ker(r)\big)(x,y)
& = &
\sum_{x'\in X}\ker(r)(x,x')\cdot r(x',y) \\
& = &
r(x,y) \\
& = &
0 \qquad \mbox{since }x\not\in X_{r} \\
\big(r \relafter \ker(r)\big)(\varphi,y)
& = &
\sum_{x'\in X}\ker(r)(\varphi,x')\cdot r(x',y) \\
& = &
\sum_{x'\in X}\varphi(x')\cdot r(x',y) \\
& = &
0, \qquad \mbox{since }\varphi\in B_{r}\subseteq \mathrm{Sol}_{r}.
\end{array}$$

\auxproof{
One has for $x,x'\in X-X_{r}$ and $\varphi,\varphi'\in B_{r}$,
$$\begin{array}{rcl}
\big(\ker(r)^{\dag} \relafter \ker(r)\big)(x,x')
& = &
\sum_{z\in X}\ker(r)(x,z)\cdot \overline{\ker(r)(x',z)} \\
& = &
\left\{\begin{array}{ll}
   1 \quad & \mbox{if }x=x' \\
   0 & \mbox{otherwise} \end{array}\right. \\
& = &
\idmap[\Ker(r)](x,x') \\
\big(\ker(r)^{\dag} \relafter \ker(r)\big)(x,\varphi)
& = &
\sum_{z\in X}\ker(r)(x,z)\cdot \overline{\ker(r)(\varphi,z)} \\
& = &
\overline{\varphi(x)} \\
& = &
0 \qquad \mbox{since }x\not\in X_{r}\mbox{ and }
   \varphi\in\Mlt_{S}(X_{r}) \\
\big(\ker(r)^{\dag} \relafter \ker(r)\big)(\varphi,\varphi')
& = &
\sum_{z\in X}\ker(r)(\varphi,z)\cdot \overline{\ker(r)(\varphi',z)} \\
& = &
\sum_{z\in X}\varphi(z)\cdot \overline{\varphi'(z)} \\
& = &
\cp(\varphi, \varphi') \\
& = &
\left\{\begin{array}{ll}
   1 \quad & \mbox{if }\varphi=\varphi' \\
   0 & \mbox{otherwise} \end{array}\right. \\
& = &
\idmap[\Ker(r)](\varphi,\varphi') \\
\end{array}$$
}

\item[-] We also check the universal property of $\ker(r)$. Let $t\colon
  Z\rightarrow X$ satisfy $r\relafter t = 0$. We split each multiset
  $t(z,-)\in\Mlt_{S}(X)$ in two parts:
$$\qquad t(z,-) = t_{1}(z, -) + t_{2}(z, -)
\quad\mbox{where}\quad
\left\{\begin{array}{l}
   \support(t_{1}(z,-))\subseteq X_{r} \\
   \support(t_{2}(z, -)) \cap X_{r} = \emptyset. \end{array}\right.$$

\noindent Then $t_{1}(z,-)\in\mathrm{Sol}_{r}$, since for each $y\in Y$,
$$\begin{array}{rcl}
0
\hspace*{\arraycolsep} = \hspace*{\arraycolsep}
(r \relafter t)(z,y)
& = &
\sum_{x}t(z,x)\cdot r(x,y) \\
& = &
\sum_{x}t_{1}(z,x)\cdot r(x,y) + t_{2}(z,x)\cdot r(x,y) \\
& = &
\sum_{x}t_{1}(z,x)\cdot r(x,y).
\end{array}$$

\noindent Since $t$ is bifinite there are only finitely many $y$ with
$t_{1}(z,y) \neq 0$. Hence each $t_{1}(z,-)\in\Mlt_{S}(X_r)$ can be
expressed in terms of finitely many base vectors from $B_r$, say as:
$$\begin{array}{rcl}
t_{1}(z, -)
& = &
a_{1}^{z}\cdot \varphi_{1}^{z} + \cdots + a_{n_z}^{z} \cdot \varphi_{n_z}^{z}.
\end{array}$$

\noindent We then define the required mediating map $t'\colon Z
\rightarrow \Ker(r)$ as function $t'\colon Z \times \big((X-X_{r})
\cup B_{r}\big) \rightarrow S$, given on $x\in X-X_{r}$ and
$\varphi\in B_{r}$ by:
$$\qquad\begin{array}{rccclcrcl}
t'(z,x)
& = &
t(z,x)
& = &
t_{2}(z,x)
& \quad\mbox{and}\quad &
t'(z,\varphi)
& = &
\left\{\begin{array}{ll}
a_{i}^{z} \quad & \mbox{if }\varphi=\varphi_{i}^{z} \\
0 & \mbox{otherwise.} \end{array}\right.
\end{array}$$

\auxproof{
For bifiniteness, the interesting part is: for each $\varphi\in
\Mlt_{S}(X_{r})$ there are only finitely many $z$ with $\varphi(z)
\neq 0$. Fix $\varphi$, say with $\support(\varphi) = \{x_{1}, \ldots,
x_{n}\}$. Write $Z_{i} = \support(t(-,x_{i})) \subseteq Z$, which is
finite. We have:
$$t'(z,\varphi) \neq 0
\Longrightarrow \ex{x}{t(z,x) \conjun \varphi(x)\neq 0}.$$

\noindent As a result $\support'(-, \varphi) \subseteq \bigcup_{i}Z_{i}$,
making it finite. The above implication holds because if 
$\all{x}{t(z,x) \neq 0 \Rightarrow \varphi(x) = 0}$, then $\varphi$
can be removed from the equation describing $t(z, -)$ in terms of
base vectors.
}

\noindent This $t'$ is bifinite, and is the right map, since:
$$\hspace*{1.5em}\begin{array}[b]{rcl}
\lefteqn{\big(\ker(r) \relafter t'\big)(z,x)} \\
& = &
\displaystyle\sum_{k\in\Ker(r)}t'(z,k)\cdot \ker(r)(k,x) \\
& = &
\Big(\displaystyle\sum_{x'\in X-X_{r}}t'(z,x) \cdot \ker(r)(x',x)\Big) +
   \Big(\sum_{\varphi\in B_{r}}t'(z,\varphi) \cdot \ker(r)(\varphi,x)\Big) \\
& = &
\left\{\begin{array}{ll}
   t'(z,x) \quad & \mbox{if }x\not\in X_{r} \\
   \sum_{i}\,a_{i}^{z}\cdot\varphi_{i}^{z}(x) \quad
      & \mbox{if }x\in X_{r} \end{array}\right. \\
& = &
\left\{\begin{array}{ll}
   t_{2}(z,x) \quad & \mbox{if }x\not\in X_{r} \\
   t_{1}(z,x) & \mbox{if }x\in X_{r} \end{array}\right. \\
& = &
t(z,x).
\end{array}\eqno{\QEDbox}$$
\end{itemize}
\end{myproof}

The category $\BifMRel_{S}$ is thus a dagger kernel category.  The
kernel subobjects of an object then form an orthomodular lattice,
see~\cite{HeunenJ10a}. Further investigation is needed to see if
$\BifMRel_{S}$ can really be seen as a ``light'' version of the
category of Hilbert spaces, suitable for discrete quantum
computations. Further steps in more logic-oriented investigations,
including measurement, can be found in~\cite{Jacobs12b}.

\subsubsection*{Acknowledgments}
Thanks to Jorik Mandemaker and Chris Heunen for helpful discussions.

\bibliographystyle{plain} 
\bibliography{../../Macros/bib}

\end{document}